# Nonlinear Fourier Analysis

**Terence Tao, Christoph Thiele**

# Contents

# Nonlinear Fourier Analysis

**Terence Tao**

**Christoph Thiele**

## 1. Introduction

These are lecture notes for a short course presented at the IAS Park City Summer School in July 2003 by the second author. The material of these lectures has been developed in cooperation by both authors.

The aim of the course was to give an introduction to nonlinear Fourier analysis from a harmonic analyst's point of view. Indeed, even the choice of the name for the subject reflects the harmonic analyst's taste, since the subject goes by many names such as for example scattering theory, orthogonal polynomials, operator theory, logarithmic integrals, continued fractions, integrable systems, Riemann Hilbert problems, stationary Gaussian processes, bounded holomorphic functions, etc.

We present only one basic model for the nonlinear Fourier transform among a large family of generalizations of our model. The focus then is to study analogues of classical questions in harmonic analysis about the linear Fourier transform in the setting of the nonlinear Fourier transform. These questions concern for example the definition of the Fourier transform in classical function spaces, continuity properties, invertibility properties, and a priori estimates. There is an abundance of analytical questions one can ask about the nonlinear Fourier transform, and we only scratch the surface of the subject.

The second half of the lecture series is devoted to showing how the nonlinear Fourier transform appears naturally in several fields of mathematics. We only present a few of the many applications that are suggested by the above (incomplete) list of names for the subject.

[1]Department of Mathematics, UCLA, Los Angeles, CA 90095-1555.
**E-mail address**: tao@math.ucla.edu.
[2]Department of Mathematics, UCLA, Los Angeles, CA 90095-1555.
**E-mail address**: thiele@math.ucla.edu.
Received by the editors July, 2003.
The first author is a Clay fellow and was partially supported by a grant form the Packard foundation. The second author was partially supported by NSF grants DMS 9985572 and DMS 9970469.

There is a vast literature on the subject of this course, in part generated by research groups with few cross-references to each other. Unfortunately we are not sufficiently expert to turn these lecture notes into anything near a survey of the existing literature. In the bibliography, we present only a small number of fairly randomly chosen entrance points to the vast literature.

We would like to thank the Park City Math Institute, its staff, and the conference organizers for organizing a stimulating and enjoyable summer school. We would like to thank R. Killip and S. Klein for carefully reading earlier versions of the manuscript and making many suggestions to improve the text. Finally, we thank J. Garnett for teaching us bounded analytic functions.

# LECTURE 1
# The nonlinear Fourier transform on $l_0$, $l^1$ and $l^p$

## 1. The nonlinear Fourier transform

In this lecture series, we study a special case of a wide class of nonlinear Fourier transforms which can be formulated at least as general as in the framework of generalized AKNS-ZS systems in the sense of ([**1**]). For simplicity we refer to the special case of a nonlinear Fourier transform in this lecture series as "the nonlinear Fourier transform", but the possibility of a more general setting should be kept in mind.

More precisely, we discuss (briefly) a nonlinear Fourier transform of functions on the real line, and (at length) a nonlinear Fourier series of coefficient sequences, i.e., functions on the integer lattice **Z**. Fourier series can be regarded as abstract Fourier transform on the circle group **T** or dually as abstract Fourier transform on the group **Z** of integers, while ordinary Fourier transform is the abstract Fourier transform of the group **R** of real numbers. We shall therefore use the word *Fourier transform* for both models which we discuss. Indeed, to the extend that we discuss the general theory here, it is mostly parallel in both models, with the possible exception of the general existence result for an inverse Fourier transform in Lecture 3 which the authors have not been able to verify in the model of the nonlinear Fourier transform of functions on the real line.

For a sequence $F = (F_n)$ of complex numbers parameterized by $n \in \mathbf{Z}$, we define the Fourier transform as

$$\widehat{F}(\theta) = \sum_{n \in \mathbf{Z}} F_n e^{-2\pi i \theta n} \tag{1}$$

and one has the inversion formula

$$F_n = \int_0^1 \widehat{F}(\theta) e^{2\pi i \theta n}\, d\theta$$

A natural limiting process takes this Fourier transform to the usual Fourier transform of functions on the real line. We have made the choice of signs in the exponents so that this limit process is consistent with the definition of the Fourier transform in [**19**].

We shall pass to a complex variable

$$z = e^{-2\pi i \theta}$$

so that (1) becomes

$$\widehat{F}(z) = \sum_{n \in \mathbf{Z}} F_n z^n$$

after identifying 1-periodic functions in $\theta$ with functions in $z \in \mathbf{T}$. The choice of sign in the exponent here is the one most convenient for us.

The discrete nonlinear Fourier transform acts on sequences $F_n$ parameterized by the integers, $n \in \mathbf{Z}$, such that each $F_n$ is a complex number in the unit disc $D$. To begin with we shall assume these sequences are compactly supported. That is, $F_n = 0$ for all but finitely many values of $n$.

For a complex parameter $z$ consider the following formally infinite recursion:

$$\begin{pmatrix} a_n & b_n \end{pmatrix} = \frac{1}{\sqrt{1-|F_n|^2}} \begin{pmatrix} a_{n-1} & b_{n-1} \end{pmatrix} \begin{pmatrix} 1 & F_n z^n \\ \overline{F_n} z^{-n} & 1 \end{pmatrix}$$

$$a_{-\infty} = 1, \quad b_{-\infty} = 0$$

Here $a_{-\infty} = 1$ and $b_{-\infty} = 0$ is to be interpreted as $a_n = 1$ and $b_n = 0$ for sufficiently small $n$, which is consistent with the recursion formula since the transfer matrix

$$\frac{1}{\sqrt{1-|F_n|^2}} \begin{pmatrix} 1 & F_n z^n \\ \overline{F_n} z^{-n} & 1 \end{pmatrix} \tag{2}$$

is the identity matrix for sufficiently small $n$ by the assumption that $F_n$ is compactly supported.

The nonlinear Fourier transform of the sequence $F_n$ is the pair of functions $(a_\infty, b_\infty)$ in the parameter $z \in \mathbf{T}$, where $a_\infty$ and $b_\infty$ are equal to $a_n$ and $b_n$ for sufficiently large $n$. We write

$$\widehat{F}(z) = (a_\infty(z), b_\infty(z))$$

We will momentarily identify the pair of functions $(a_\infty, b_\infty)$ with an $SU(1,1)$ valued function on $\mathbf{T}$.

While evidently $a_\infty$ and $b_\infty$ are finite Laurent polynomials in $z$ (rational functions with possible poles only at 0 and $\infty$), we regard the nonlinear Fourier transform as functions on the unit circle $\mathbf{T}$. Later, when we consider properly infinite sequences $F_n$, restriction to $\mathbf{T}$ as domain will be a necessity.

Observe that for $z \in \mathbf{T}$ the transfer matrices are all in $SU(1,1)$. Hence we can write equivalently for the above recursion

$$\begin{pmatrix} a_n & b_n \\ \overline{b_n} & \overline{a_n} \end{pmatrix} = \frac{1}{\sqrt{1-|F_n|^2}} \begin{pmatrix} a_{n-1} & b_{n-1} \\ \overline{b_{n-1}} & \overline{a_{n-1}} \end{pmatrix} \begin{pmatrix} 1 & F_n z^n \\ \overline{F_n} z^{-n} & 1 \end{pmatrix}$$

with

$$\begin{pmatrix} a_{-\infty} & b_{-\infty} \\ \overline{b_{-\infty}} & \overline{a_{-\infty}} \end{pmatrix} = \begin{pmatrix} 1 & 0 \\ 0 & 1 \end{pmatrix}$$

and all matrices

$$\begin{pmatrix} a_n & b_n \\ \overline{b_n} & \overline{a_n} \end{pmatrix}$$

are in $SU(1,1)$, and in particular $|a_n|^2 = 1 + |b_n|^2$.

Thus the Fourier transform can be regarded as a map

$$l_0(\mathbf{Z}, D) \to C(\mathbf{T}, SU(1,1))$$

where $l_0(\mathbf{Z}, D)$ are the compactly supported sequences with values in $D$, and $C(\mathbf{T}, SU(1,1))$ are the continuous functions on $\mathbf{T}$ with values in $SU(1,1)$.

While we shall not do this here, one can naturally define similar nonlinear Fourier transforms for a variety of Lie groups in place of $SU(1,1)$. The group $SU(2)$ leads to an interesting example. We remark that here we define Fourier transforms using Lie groups in a quite different manner from the way it is done in representation theory. There one defines Fourier transforms of complex valued functions on groups, and one remains in the realm of linear function spaces. Here we end up with group valued functions, a much more nonlinear construction.

If $E$ is an open set in the Riemann sphere, define $E^*$ to be the set reflected across the unit circle, i.e,

$$E^* = \{z : \overline{z}^{-1} \in E\}$$

The operation $*$ is the identity map on $E \cap \mathbf{T}$.

If $c$ is a function on $E$, define

$$c^*(z) = \overline{c(\overline{z}^{-1})}$$

as a function on $E^*$. This operation preserves analyticity. On the circle $T$, this operation coincides with complex conjugation:

$$c^*(z) = \overline{c(z)}$$

for all $z \in \mathbf{T} \cap E$.

We then observe the recursion

$$(3) \qquad \begin{pmatrix} a_n & b_n \\ b_n^* & a_n^* \end{pmatrix} = \frac{1}{\sqrt{1-|F_n|^2}} \begin{pmatrix} a_{n-1} & b_{n-1} \\ b_{n-1}^* & a_{n-1}^* \end{pmatrix} \begin{pmatrix} 1 & F_n z^n \\ \overline{F_n} z^{-n} & 1 \end{pmatrix}$$

with

$$\begin{pmatrix} a_{-\infty} & b_{-\infty} \\ b_{-\infty}^* & a_{-\infty}^* \end{pmatrix} = \begin{pmatrix} 1 & 0 \\ 0 & 1 \end{pmatrix}$$

All entries in these matrices are meromorphic functions on the entire Riemann sphere. Namely, these recursions hold on $\mathbf{T}$ and thus hold on the entire sphere by meromorphic continuation of $a_n, b_n, a_n^*, b_n^*$.

Observe that the matrix

$$\begin{pmatrix} a_n(z) & b_n(z) \\ b_n^*(z) & a_n^*(z) \end{pmatrix}$$

is not necessarily in $SU(1,1)$ for $z$ outside the circle $\mathbf{T}$. However,

$$a_n a_n^* = 1 + b_n b_n^*$$

continues to hold on the complex plane since it holds on the circle $\mathbf{T}$.

Thinking of the pair $(c, d)$ as the first row of an element of a function which takes values in $SU(1,1)$ on the circle $\mathbf{T}$, we shall use the convention to write

$$(a, b)(c, d) = (ac + bd^*, ad + bc^*)$$

For small values of $F_n$ the nonlinear Fourier transform is approximated the linear inverse Fourier transform. This can be seen by linearizing in $F$. The factor $(1-|F_n|^2)^{-1/2}$ is quadratic and we disregard it. The remaining formula for $a_\infty$ and $b_\infty$ is polynomial in $F$ and $\overline{F}$. If we only collect the constant and the linear term, we obtain

$$(a_\infty, b_\infty) = (1, \sum_{n \in \mathbf{Z}} F_n z^n)$$

Thus $a_\infty$ is constant equal to 1 in linear approximation and $b_\infty$ is the Fourier transform

$$\sum_{n \in \mathbf{Z}} F_n z^n$$

in linear approximation.

The following lemma summarizes a few algebraic properties of the nonlinear Fourier transform.

**Lemma 1.** *If $F_n = 0$ for $n \neq m$, then*

$$\widehat{(F_n)} = (1 - |F_m|^2)^{-1/2}(1, F_m z^m) \tag{4}$$

*If $\widehat{(F_n)} = (a, b)$, then we have for the shifted sequence whose $n$-th entry is $F_{n+1}$*

$$\widehat{(F_{n+1})} = (a, bz^{-1}) \tag{5}$$

*If the support of $F$ is entirely to the left of the support of $G$, then*

$$\widehat{(F+G)} = \widehat{F}\,\widehat{G} \tag{6}$$

*If $|c| = 1$ then*

$$\widehat{(cF_n)} = (a, cb) \tag{7}$$

*For the reflected sequence whose $n$-th entry is $F_{-n}$*

$$\widehat{(F_{-n})}(z) = (a^*(z^{-1}), b(z^{-1})) \tag{8}$$

*Finally, for the complex conjugate of a sequence, we have*

$$\widehat{(\overline{F_n})}(z) = (a^*(z^{-1}), b^*(z^{-1})) \tag{9}$$

Observe that statements (5),(7), (8) and (9) are exactly the behaviour of the linearization $a \sim 1$ and $b \sim \sum F_n z^n$. Statements (5) and (7) are most easily proved by conjugation with diagonal elements in $SU(1,1)$. Concerning statements (8) and (9), observe that under the reflection $n \to -n$ or under complex conjugation the Laurent expansion of the diagonal element $a$ turns into the expansion with complex conjugate coefficients.

## 2. The image of finite sequences

Our next concern is the space of functions $a_\infty, b_\infty$ obtained as Fourier transforms of finite $D$-valued sequences $F_n$.

It is immediately clear that $a$ and $b$ are finite Laurent polynomials. The following lemma describes the degree of these Laurent polynomials. Define the upper degree of a Laurent polynomial to be the largest $N$ such that the $N$-th coefficient is nonzero, and define the lower degree to be the least $N$ such that the $N$-th coefficient is nonzero.

**Lemma 2.** *Let $F_n$ be a nonzero finite sequence with NLFT $(a, b)$. Let $N_-$ be the smallest integer such that $F_{-N} \neq 0$. and let $N_+$ be the largest integer such that $F_{N_+} \neq 0$, Then $a$ is a Laurent polynomial*

$$a = \sum_{n=N_- - N_+}^{0} \check{a}(n) z^k$$

*with exact lowest degree $N_- - N_+$ and exact highest degree $0$. The constant term of this Laurent polynomial is*

$$\check{a}(0) = \prod_k (1 - |F_k|^2)^{-1/2}$$

*Moreover, $b$ is of the form*

$$b = \sum_{k=N_-}^{N_+} \check{b}(n) z^n$$

*with exact highest degree $N_+$ and exact lowest degree $N_-$.*

A particular consequence of this lemma is that the order of the highest and lowest nonvanishing coefficients for the Laurent polynomial of $b$ are the same as for the sequence $F$.

**Proof.** The Lemma can be proved by induction on the length $1 + N_+ - N_-$ of the sequence $F_n$. If the length is 1, then $(a, b)$ is equal to a transfer matrix (set $N = N_- = N_+$),

$$(a, b) = \left((1 - |F_N|^2)^{-1/2}, (1 - |F_N|^2)^{-1/2} F_N z^N\right)$$

This proves the lemma for sequences of length one.

Now let $l > 1$ and assume the theorem is true for lengths less than $l$ and assume $F_n$ has length $l$. Set $N = N_+$, then we have

$$a = (1 - |F_N|^2)^{-1/2} a' + (1 - |F_N|)^{-1/2} z^{-N} \overline{F_N} b'$$

where $(a', b')$ is the NLFT of the truncated sequence $F'$ which coincides with $F$ everywhere except for the $N$-th entry where we have $F'_N = 0$. By induction, $z^{-N} b'$ has highest degree at most $-1$ (since $b'$ has highest degree at most $N - 1$), so the constant term of $a$ is $(1 - |F_n)^{-1/2}$ times the constant term of $a'$ and there are no terms of positive order of $a$. The lowest order term of $a'$ is at least $N_- - N_+ + 1$, while the lowest order term of $z^{-N} b'$ is exactly $N_- - N_+$. Thus for $a$ we have the lowest and highest order coefficients claimed in the lemma.

Similarly, we have

$$b = (1 - |F_N|^2)^{-1/2} a' z^N F_N + (1 - |F_N|)^{-1/2} b'$$

Only the second summand produces a nonzero coefficient of degree $N_-$ and this is the lowest degree of $b$. Only the first summand produces a nonzero coefficient or order $N_+$ and this is the highest order coefficient of $b$. □

The meromorphic extensions of $a$ and $b$ to the Riemann sphere satisfy the recursion (3). On the Riemann sphere, we shall be interested in the open unit disc $D = \{z : |z| < 1\}$ and the unit disc at infinity $D^* = \{1/z : |z| < 1\}$. Observe that $a$ is holomorphic on $D^*$.

**Lemma 3.** *Let $(a, b) = \widehat{F}$ for some finite sequence $F$. Then $a$ has no zeros in the disc at infinity $D^*$.*

**Proof.** It suffices to prove the lemma under the assumption $F_n = 0$ for $n < 0$, because we can translate $F$ and use Lemma 1.

The constant term of $a$ is nonzero, and therefore $a$ is not zero at $\infty$.

For $|z| > 1$ we rewrite the recursion

$$\begin{pmatrix} a_n & b_n \end{pmatrix} = \frac{1}{\sqrt{1 - |F_n|^2}} \begin{pmatrix} a_{n-1} & b_{n-1} \end{pmatrix} \begin{pmatrix} 1 & F_n z^n \\ \overline{F_n} z^{-n} & 1 \end{pmatrix}$$

as

$$\begin{pmatrix} |z|^n a_n & b_n \end{pmatrix} = \frac{1}{\sqrt{1 - |F_n|^2}} \begin{pmatrix} |z|^n a_{n-1} & b_{n-1} \end{pmatrix} \begin{pmatrix} 1 & F_n (z/|z|)^n \\ \overline{F_n} (z/|z|)^{-n} & 1 \end{pmatrix}$$

Therefore

$$|z^n a_n|^2 - |b_n|^2 = |z^n a_{n-1}|^2 - |b_{n-1}|^2 \geq |z^{n-1} a_{n-1}|^2 - |b_n|^2$$

because $|z| > 1$. Now it follows by induction that

$$|z^n a_n|^2 - |b_n|^2 > 0$$

because this is true for $n$ near $-\infty$. Thus $a$ is not zero at $z$. □

**Corollary 1.** *We have $|a_n(z)| \geq 1$ for $|z| > 1$. Moreover, $a_n(\infty)$ is positive and greater than or equal to* 1.

**Proof.** Since $a$ has no zeros, the function $\log |a|$ is harmonic on $D^*$. On the boundary it is non-negative, and by the maximum principle it is non-negative on $D^*$. It remains to see positivity of $a(\infty)$, but this is clear since $a(\infty)$ is the constant term in the Laurent polynomial of $a$. □

Next we observe that if $(a, b)$ is the NLFT of a finite sequence, then $a$ is already determined by $b$.

**Lemma 4.** *Let $b$ be a Laurent polynomial. Then there exists a unique Laurent polynomial $a$ such that $aa^* = 1 + bb^*$, $a$ has no zeros in $D^*$, and $a(\infty) > 0$.*

**Proof.** Observe that $P = 1 + bb^*$ is a nonzero Laurent polynomial. It can only have poles at 0 and $\infty$, and by symmetry the order of these two poles have to be equal. Assume $P$ has pole of oder $n$ at 0, clearly $n \geq 0$. By symmetry $P$ has pole of order $n$ at $\infty$. Since $P$ is clearly has no zeros on $\mathbf{T}$, by symmetry it has exactly $n$ zeros in $D$ and $n$ zeros in $D^*$.

Uniqueness: In order for $aa^*$ to have the correct order of pole at 0, the Laurent polynomial $a$ has to be a polynomial in $z^{-1}$ of order $n$. In order for $aa^*$ to have the same zeros as $P$, it has to have the same zeros as $P$ in $D$ and no other zeros. This determines $a$ up to a scalar factor. The condition $a(\infty) > 0$ determines $a$ up to a positive factor. Since $1 + bb^*$ is nonzero, this positive factor is determined by $1 + bb^* = aa^*$.

Existence: Let $a$ be a polynomial of degree $n$ in $1/z$ whose zeros are exactly the zeros of $P$ in $D$. Then $a(\infty) \neq 0$. By multiplying by a phase factor we may assume $a(\infty) > 0$. By multiplying by a positive factor we may assume that $a^*a$ coincides with $P$ on at least one point of $\mathbf{T}$. Then $1 + bb^* - aa^*$ has $2n + 1$ zeros, but at most two poles of order $n$. Thus $1 + bb^* - aa^* = 0$. □

We are now ready to characterize the target space of the Fourier transform of finite sequences.

**Theorem 1.** *The nonlinear Fourier transform is a bijection from the set of all finite sequences $(F_n)$ in the unit disc into the space of all pairs $(a, b)$ with $b$ an arbitrary Laurent polynomial and $a$ the unique Laurent polynomial which satisfies $aa^* = 1 + bb^*$, $a(\infty) > 0$, and has no zeros in $D^*$.*

Remark: While not stated explicitly in the theorem, all pairs $(a, b)$ described in the theorem satisfy not only $a(\infty) > 0$ but also $a(\infty) \geq 1$. Moreover, $a$ and $b$ have the same length. This follows from the theorem and the previous discussion.

**Proof.** Clearly the Fourier transform maps into the described space by the previous discussion.

We know that the upper and lower degree of $F$ are the same as the upper and lower degree of $b$. Thus it suffices to prove bijectivity under the assumption of fixed upper and lower degree of $F$ and $b$. By shifting $F$ and $b$ we may assume both have lower degree 0, and we can use induction on the common upper degree $N$.

In the case $N = 0$, meaning $F \equiv 0$, and the in case $N = 1$ we have

$$b = F_0(1 - |F_0|^2)^{-1/2}$$

and the map $F_0 \to b$ is clearly a bijection from $D$ to $\mathbf{C}$. Now assume we have proved bijectivity up to upper degree $N - 1$.

We first prove injectivity, i.e., $F$ of upper degree $N$ can be recovered from $(a, b)$. It suffices to show that $F_0$ can be recovered from $(a, b)$. Then we can by induction recover the truncated sequence $F'$, which coincides with $F$ except for $F'_0 = 0$ from the Fourier transform of $F'$ which can be calculated as

$$(a', b') = (1 - |F_0|^2)^{-1/2}(1, -F_0)(a, b)$$

However, this last identity also implies that

$$0 = b'(0) = (1 - |F_0|^2)^{1/2} b'(0) = b(0) - F_0 a^*(0)$$

Hence

$$F_0 = \frac{b(0)}{a^*(0)}$$

and this quotient is well defined since $a^*(0) \neq 0$. Thus $F_0$ is determined by $(a, b)$ and we have proved injectivity for upper degree $N$.

Next, we prove surjectivity. Let $b$ have upper degree $N$ and let $a$ be the unique Laurent polynomial which satisfies $aa^* = 1 + bb^*$, $a(\infty) > 0$, and has no zeros in $D^*$. We set formally

$$F_0 = \frac{b(0)}{a^*(0)}$$

Observe that $|F_0| < 1$ since $b/a^*$ is holomorphic in $D$, continuous up to $\mathbf{T}$, and bounded by 1 on $\mathbf{T}$. We calculate formally the truncated Fourier transform data

$$(a', b') = (1 - |F_0|^2)^{-1/2}(1, -F_0)(a, b)$$

Then $b'$ is a polynomial of upper degree at most $N$ and lower degree at least 1.

It now suffices to prove that $(a', b')$ is the nonlinear Fourier transform of a sequence of length $N - 1$. For this it suffices to show that $a'$ is the unique Laurent polynomial such that $1 + b'b'^* = a'a'^*$, $a(\infty) > 0$, and $a'$ has no zeros in $D^*$.

However,

$$1 + b'b'^* = a'a'^*$$

holds on $\mathbf{T}$ and therefore everywhere since the determinant of the matrix $(a, b)$ coincides with that of $(a', b')$ on $\mathbf{T}$. The recursion

$$(a, b) = (1 - |F_0|^2)^{-1/2}(1, F_0)(a', b')$$

implies

$$a(\infty) = (1 - |F_0|^2)^{-1/2} a'(\infty)$$

and thus $a'(\infty) > 0$. Finally, we observe that

$$a'(z) = (1 - |F_0|^2)^{-1/2}(a - F_0 b^*)$$

does not vanish in $D^*$ since $b^*/a$ is bounded by 1 in $D^*$ and thus $a$ strictly dominates $F_0 b^*$ in $D^*$. $\square$

## 3. Extension to $l^1$ sequences

We have defined the Fourier transform for finite sequences. Now, we would like to extend the definition to infinite sequences.

As in the case of the linear Fourier transform, the defining formula actually extends to sequences in $l^1(\mathbf{Z}, D)$, i.e., summable sequences of elements in $D$.

Define a metric on the space $SU(1,1)$ by

$$\mathrm{dist}(G, G') = \|G - G'\|_{op}$$

This clearly makes $SU(1,1)$ a complete metric space, since $\mathbf{C}^4$ is a complete metric space and $SU(1,1)$ is a closed subset of $\mathbf{C}^4$ with the inherited topology.

Define $L^\infty(\mathbf{T}, SU(1,1))$ to be the metric space of all essentially bounded functions $G : \mathbf{T} \to SU(1,1)$

$$\sup_z \mathrm{dist}(\mathrm{id}, G(z)) < \infty$$

(in the usual sense of the essential supremum) with the distance

$$\mathrm{dist}(G, G') = \sup_z \mathrm{dist}(G(z), G'(z))$$

On the space of all summable sequences in $D$ define the distance

$$\mathrm{dist}(F, F') = \sum_n \|T_n - T'_n\|_{op}$$

where $T_n$ denotes the transfer matrix defined in (2). This makes $l^1(\mathbf{Z}, D)$ a complete metric space. We claim that on sets

$$B_\epsilon = \{F_n : \sup_n |F_n| < 1 - \epsilon\}$$

with $\epsilon > 0$ (every element in $l^1(\mathbf{Z}, D)$ is in such a set, and also every Cauchy sequence in $l^1(\mathbf{Z}, D)$ is inside one of these sets) this distance is bi-Lipschitz to

$$\mathrm{dist}'(F, F') = \sum_n |F_n - F'_n|$$

Namely, if $F_n$ and $F'_n$ are in $B_\epsilon$, then

$$\|T_n - T'_n\|_{op} =$$

$$= \left|(1 - |F_n|^2)^{-1/2} - (1 - |F_n|^2)^{-1/2}\right| + \left|(1 - |F_n|^2)^{-1/2} F_n - (1 - |F'_n|^2)^{-1/2} F'_n\right|$$

This is bounded by a constant depending on $\epsilon$. Thus we only need to show equivalence to $|F_n - F'_n|$ if the latter is smaller than a constant depending on $\epsilon$. This however follow easily by Taylor expansion of the nonlinear terms in the expression for $\|T_n - T'_n\|_{op}$.

In particular, we observe that the finite sequences are dense in $l^1$.

**Lemma 5.** *With the above metrics, the NLFT on $l_0(\mathbf{Z}, D)$ extends uniquely to a locally Lipschitz map from $l^1(\mathbf{Z}, D)$ to $L^\infty(\mathbf{T}, SU(1,1))$. The NLFT of such sequences can be written as the convergent infinite ordered product of the transfer matrices.*

**Proof.** To prove existence and uniqueness of the extension, it suffices to prove the Lipschitz estimate on bounded sets for finite sequences.

Using $\|T_n\| \geq 1$ we have by Trotter's formula:

$$\|\prod_n T_n - \prod_n T'_n\|_{op} \leq [\sum_n \|T_n - T'_n\|_{op}][\prod_n \|T_n\|_{op}][\prod_n \|T'_n\|_{op}]$$

Moreover, we have

$$\prod_n \|T_n\|_{op} \leq \exp[\sum_n [\|T_n\|_{op} - 1]] \leq \exp[\sum_n \|T_n - \mathrm{id}\|_{op}]$$

hence the right-hand side remains bounded on bounded sets of $l^1(\mathbf{Z}, D)$. Thus, on a bounded set we have

$$\|\prod_n T_n - \prod_n T'_n\|_{op} \leq C[\sum_n \|T_n - T'_n\|_{op}]$$

with $C$ depending on the set. This proves the Lipschitz estimate on bounded sets.

By the abstract theory of metric spaces, the NLFT extends to a locally Lipshitz map on $l^1(\mathbf{Z}, D)$. Given any sequence $F$, the truncations of the sequence to the interval $[-N, N]$ converge in $l^1(\mathbf{Z}, D)$ to $F$, and thus the sequence of nonlinear Fourier transforms converges to the NLFT of $F$.

Convegence in the target space is not only in $L^\infty(\mathbf{T}, SU(1,1))$, but also in the space $C(\mathbf{T}, SU(1,1))$ of continuous functions, which is a closed subspace. This implies that the NLFT of the truncated sequences converge pointwise and uniformly to the NLFT of $F$.

Observe that if $F \in B_\epsilon$, then the truncations of $F_n$ to intervals $[-N, N]$ remain in $B_\epsilon$, and converge to $F$ in $l^1(\mathbf{Z}, D)$. Thus the products of the transfer matrices converge to the NLFT of $F$. □

We observe that for $F \in l^1(\mathbf{Z}, D)$ we have

$$\sup_{z \in \mathbf{T}} \|(a(z), b(z))\|_{op} \leq \prod_n \|(1 - |F_n|^2)^{-1/2}(1, F_n)\|$$

or, applying the logarithm to both sides and using Lemma 33:

$$\sup_z \operatorname{arccosh}|a(z)| \leq \sum_n \operatorname{arccosh}((1 - |F_n|^2)^{-1/2})$$

Define $g(y) = (\log(\cosh(y)))^{1/2}$. Then $g$ vanishes at 0 and is concave on the positive half axis, and therefore $g(x) + g(y) \geq g(x + y)$ for all $0 \leq x, y$, and the analogue inequality holds for any countable number of summands. Applying $g$ to the last display we thus obtain

$$\sup_z (\log|a(z)|)^{1/2} \leq \sum_n (\log((1 - |F_n|^2)^{-1/2}))^{1/2} \tag{10}$$

In the following section, this estimate will be compared to estimates for sequences in spaces $l^p(\mathbf{Z}, D)$ for various $p$.

## 4. Extension to $l^p$ sequences, $1 < p < 2$

In this section we define the nonlinear Fourier transform of $l^p$ sequences with $1 < p < 2$. The discussion in this section is an extraction from the work of Christ and Kiselev. Mainly we rely on [**4**] and we state and use but not prove theorems from that paper.

Let $F \in l_0(\mathbf{Z}, D)$ be a finite sequence. By the distributive law, we can write

$$\widehat{F}(z) = \prod_{n \in \mathbf{Z}} (1 - |F_n|^2)^{-1/2}(1, F_n z^n)$$

$$= \left(\prod_{n \in \mathbf{Z}} (1 - |F_n|^2)^{-1/2}\right)\left(\sum_{n=0}^{\infty}\left[\sum_{i_1 < \cdots < i_n} \prod_{k=1}^{n} (0, F_{i_k} z^{i_k})\right]\right)$$

Here all formally infinite products and sums are actually finite since almost all factors and summands are trivial, and where the summand for $n = 0$ on the right-hand side is equal to $(1, 0)$.

Observe that the first factor in the last display is independent of $z$ and is a convergent product under the assumption $F \in l^p(\mathbf{Z}, D)$ with $1 < p < 2$. The second factor in the last display is a multilinear expansion, i.e., a Taylor expansion in the sequence $F$ near the trivial sequence $F = 0$.

We shall see that for $F \in l^p(\mathbf{Z}, D)$, each term in this multilinear expansion is well defined as a measurable function in $z$ and that the multilinear expansion is absolutely summable for almost all $z \in \mathbf{T}$. This allows us to define the nonlinear Fourier transform for $l^p$ sequences as a measurable function on $\mathbf{T}$.

**Theorem 2.** *Let $1 \le p < 2$ and let $p'$ be the dual exponent $p/(p-1)$. Let $F \in l^p(\mathbf{Z}, D)$. Then the multilinear term*

$$\sum_{i_1 < \cdots < i_n} \left[ \prod_{k=1}^{n} (0, F_{i_k} z^{i_k}) \right] \tag{11}$$

*is a well defined element of the quasi-metric space $L^{p'/n}(\mathbf{T}, M_{2\times 2})$ and depends continuously on the sequence $F \in l^p(\mathbf{Z}, D)$. The multilinear expansion*

$$\sum_{n=0}^{\infty} \left[ \sum_{i_1 < \cdots < i_n} \prod_{k=1}^{n} (0, F_{i_k} z^{i_k}) \right] \tag{12}$$

*is absolutely summable for almost every $z$. Defining*

$$(a, b)(z) := \left( \prod_{n \in \mathbf{Z}} (1 - |F_n|^2)^{-1/2} \right) \left( \sum_{n=0}^{\infty} \left[ \sum_{i_1 < \cdots < i_n} \prod_{k=1}^{n} (0, F_{i_k} z^{i_k}) \right] \right)$$

*we have $|a|^2 = 1 + |b|^2$ almost everywhere, the function $a$ has an outer extension to $D^*$ with $a(\infty) > 0$, and we have the estimate*

$$\|(\log |a|)^{1/2}\|_{L^{p'}(\mathbf{T})} \le C_p \left\| |\log(1 - |F_.|^2)|^{1/2} \right\|_{l^p(\mathbf{Z})} \tag{13}$$

The case $p = 1$ of inequality (13) has been observed in (10). Indeed, the case $p = 1$ of this theorem is considerably easier than the case $p > 1$.

The multilinear expansion described in this theorem fails to converge in general if $p = 2$, see [**15**]. However, inequality (13) remains true for $p = 2$ if the nonlinear Fourier transform is defined properly. We will discuss this in subsequent sections. It is an interesting open problem whether $C_p$ in inequality (13) can be chosen uniformly in $p$ as $p$ approaches 2.

**Proof.** The quasi-metric of the space $L^q(\mathbf{T}, M_{2\times 2})$ is defined as

$$\left( \int_{\mathbf{T}} \|G(z)\|_{op}^q \right)^{1/q}$$

To prove that each multilinear term (11) is a well defined element in $L^{p'/n}(\mathbf{T}, M_{2\times 2})$, it suffices to prove that the multilinear map $T_n$, originally defined on finite sequences by

$$T_n(F^{(1)}, \dots, F^{(n)})(z) = \sum_{i_1 < \cdots < i_n} \left[ \prod_{k=1}^{n} (0, F_{i_k}^{(k)} z^{i_k}) \right]$$

satisfies the a priori estimate

$$||T_n(F^{(1)}, \dots, F^{(n)})||_{p'/n} \le C \prod_{k=1}^{n} ||F^{(k)}||_p$$

By the general theory of multilinear maps on topological vector spaces, the map $T_n$ extends then uniquely to a continuous map

$$l^p(\mathbf{Z}) \times \cdots \times l^p(\mathbf{Z}) \to L^{p'/n}(\mathbf{T})$$

To prove the above a priori estimate, we use the following theorem formulated slightly differently in [**3**]:

**Theorem 3.** *Let $1 < p < 2$. Let $k_j$ for $j = 1, \dots, n$ be locally integrable functions on $\mathbf{R} \times \mathbf{R}$ such that the map*

$$K_j f(y) := \int k_j(y, x) f(x)\, dx \tag{14}$$

*which is defined on bounded compactly supported $f$ satisfies the a priori bound*

$$||K_j f||_{p'} \le C||f||_p$$

*Then the operator $T_n$ defined by*

$$T_n(f_1, \dots, f_n) = \int_{x_1 < \cdots < x_n} \prod_{j=1}^{n} k_j(y, x_j) f_j(x_j)\, dx_j$$

*satisfies the a priori bound*

$$||T_n(f_1, \dots, f_n)||_{p'/n} \le C \prod_{j=1}^{n} ||f_j||_p$$

For the proof of this theorem we refer to [**3**] with improvements on the constant $C$ in [**4**].

To apply the theorem to the case at hand we need to convert the integral in (14) to a sum. This is easily done by considering functions that are constant on each interval $[l, l+1)$ for $l \in \mathbf{Z}$, but some care is to be taken so that the integration in the definition of $T_n$ can be turned into a summation over a discrete set. Specifically, define $k_j(y, x)$ to be zero for $y \notin [0, 2\pi]$ and

$$k_j(x, y) = e^{\pm iy[x/n]}$$

where $[x]$ denotes the largest integer smaller than $x$ and the sign $\pm$ is positive or negative depending on whether $j$ is odd or even. Further we define $f_j$ such that for every integer $m$, the restriction of $f_j$ to the interval $[nm, n(m+1))$ is equal to $F_m 1_{[n(m+1)-j, n(m+1)-j+1)}$ if $j$ is odd and equal to the complex conjugate of this expression if $j$ is even. With $z = e^{iy}$ one observes that for odd $n$

$$T_n(f_1, \dots, f_n) = \sum_{i_1 < \cdots < i_n} F_{i_1} z^{i_1} \overline{F_{i_2} z^{i_2}} \cdots F_{i_n} z^{i_n}$$

and for even $n$

$$T_n(f_1, \dots, f_n) = \sum_{i_1 < \cdots < i_n} F_{i_1} z^{i_1} \overline{F_{i_2} z^{i_2}} \dots \overline{F_{i_n} z^{i_n}}$$

Thus the desired bound for (11) follows from Theorem 3.

To obtain good bounds of the multilinear terms and conclude that the expansion (12) is absolutely summable almost everywhere we invoke the next theorem from [**4**].

Define a martingale structure on $\mathbf{R}$ to be a collection of intervals $E_j^m$ with $m \geq 0$ and $1 \leq j \leq 2^m$ such that the following conditions are satisfied modulo endpoints

1. The union $\cup_j E_j^m$ is equal to $\mathbf{R}$ for every $m$
2. The intervals $E_j^m$ and $E_{j'}^m$ are disjoint for $j \neq j'$
3. If $j < j'$, $x \in E_j^m$ and $x' \in E_{j'}^m$, then $x < x'$
4. For every $j, m$ we have $E_j^m = E_{2j-1}^{m+1} \cup E_{2j}^{m+1}$

Given such a martingale structure, to each locally integrable function $f$ we associate $g_f \in [0, \infty]$ by

$$g_f = \sum_{m=1}^{\infty} \left( \sum_{j=1}^{2^m} | \int_{E_j^m} f|^2 \right)^{1/2}$$

**Theorem 4.** *There is a constant $B$ such that the following holds. Define*

$$T_n(f_1, \dots, f_n) := \int_{x_1 < \dots < x_n} \prod_{i=1}^{n} f_i(x_i)\, dx_i \tag{15}$$

*and let $X$ be a finite set of locally integrable functions on $\mathbf{R}$. Assume that we are given a fixed martingale structure and define*

$$g := \max_{f \in X} g_f$$

*Then for every $n > 1$ and every $f_1, \dots, f_n \in X$,*

$$|T_n(f_1, \dots, f_n)| \leq (n!)^{-1/2} B^n g^n$$

For each parameter $y$, we apply this theorem with the same $k_i$ and $f_i$ as in the application of the previous theorem. Thus the number $g$ depends on the parameter $y$. Writing again $z = e^{iy}$ we obtain for even $n$

$$\sum_{i_1 < \dots < i_n} F_{i_1} z^{i_1} \overline{F_{i_2} z^{i_2}} \dots \overline{F_{i_n} z^{i_n}} \leq (n!)^{-1/2} B^n g^n(z)$$

and similarly for odd $n$. If we can show that for a proper choice of the martingale structure the function $g(z)$ is finite for almost every $z$, then this implies immediately that the multilinear expansion (12) converges for almost all $z$.

We choose the martingale structure adapted to $f_i$ in the sense of [**4**] (observe that all the $f_i$ are identical up to complex conjugation). As in the remark to Theorem 1.1 of [**4**] one checks that

$$\|g\|_{p'} \leq \|F\|_p$$

Thus in particular the expansion (12) converges for almost every $z$.

The diagonal entry $a$ of $(a, b)$ is only affected by the even terms in the multilinear expansion (12). Thus we obtain with the same martingale structure as before

$$|a| \leq [\prod_n (1 - |F_n|^2)^{-1/2}] \sum_n ((2n)!)^{-1/2} B^{2n} g^{2n}$$

$$\leq [\prod_n (1 - |F_n|^2)^{-1/2}] \sum_n (n!)^{-1} B^{2n} g^{2n} = [\prod_n (1 - |F_n|^2)^{-1/2}] \exp(B^2 g^2)$$

$$(\log|a|)^{1/2} \le C[\sum_n |\log(1-|F_n|^2)|]^{1/2} + Cg$$

This easily proves (13).

To complete the proof of Theorem 2 it remains to prove that the nonlinear Fourier transform $(a,b)$ of an $l^p$ sequence, which we have now defined using the multilinear expansion, satisfies $|a|^2 = 1+|b|^2$ and $a$ has an outer extension to $D^*$. We shall, for simplicity, only argue for sequences supported on $\mathbf{Z}_{\ge 0}$. The case of sequences on the full line $\mathbf{Z}$ is then a slight variation using truncations at both ends of the sequence.

Assume $F \in l^p(\mathbf{Z}_{\ge 0}, D)$ and consider the nonlinear Fourier transforms

$$(a^{(\le N)}, b^{(\le N)})$$

of the truncations $F^{(\le N)}$. If we can show that $(a^{(\le N)}, b^{(\le N)})$ converge to $(a,b)$ almost everywhere, then we obtain immediately $|a|^2 = 1+|b|^2$. Moreover, with the additional a priori estimate (16) and Lebesgue's dominated convergence theorem we have that $\log|a^{(\le N)}|$ converges to $\log|a|$ in $L^1$ and one easily concludes that $a$ is outer. □

**Theorem 5.** *Let $F \in l^p(\mathbf{Z}_{\ge 0}, D)$. With the notation as above, the sequence $(a^{(\le N)}, b^{(\le N)})$ converges for almost every $z$ to $(a,b)$. Moreover, we have the a priori estimate*

$$\|\sup_N \log|a^{(\le N)}|\|_L^{p'} \le C_p \||\log((1-|F_n|^2)|^{1/2}\|_{l^p} \tag{16}$$

**Proof.** We first show that almost everywhere convergence follows from the a priori estimate (16).

Let $M > N$ an write

$$(a^{(\le M)}, b^{(\le M)}) = (a^{(\le N)}, b^{(\le N)})(a', b')$$

By continuity of multiplication in $SU(1,1)$, we have to show smallness of $(a',b')$ for large $N$, independently of $M$ and outside a set of prescribed small measure.

This however is precisely what (16) provides if applied to the tail $F^{(>N)}$.

Thus it remains to prove (16). This a priori estimate follows by arguments similar to those given before with the following theorem taken from [**4**], which is a modification of Theorem 4.

For a given martingale structure $E_j^m$ we define

$$\tilde{g}(f) = \sum_{m=1}^{\infty} m \left( \sum_{j=1}^{2^m} |\int_{E_j^m} f|^2 \right)^{1/2}$$

**Theorem 6.** *Define*

$$M_n(f_1,\dots,f_n) := \sup_{y,y'} \left| \int_{y<x_1<\dots<x_n<y'} \prod_{i=1}^n f_i(x_i)\, dx_i \right| \tag{17}$$

*and let $X$ be a finite set of locally integrable functions on $\mathbf{R}$. Assume we are given a fixed martingale structure and define*

$$\tilde{g} := \max_{f\in X} \tilde{g}(f)$$

*Then for every $n > 1$ and every $f_1,\dots,f_n \in X$*

$$|M_n(f_1,\dots,f_n)| \le (n!)^{-1/2} B^n \tilde{g}^n$$

Applying this theorem as before proves Theorem 5. This also completes the proof of Theorem 2. □

# LECTURE 2
## The nonlinear Fourier transform on $l^2(\mathbf{Z}_{\geq 0})$

### 1. Extension to half-infinite $l^2$ sequences

Most of this section is an adaption from an article by Sylvester and Winebrenner [**20**].

Assume $F$ is a square summable sequence with values in the open unit disc, i.e., an element of $l^2(\mathbf{Z}, D)$.

As for the linear Fourier transform, the defining equation for the nonlinear Fourier transform of $F$ (infinite product of transfer matrices) does not necessarily converge for given $z \in \mathbf{T}$. Indeed, almost everywhere convergence of the partial products - a nonlinear version of a theorem of Carleson - is an interesting open problem, see the discussion in [**16**].

It however converges in a certain $L^2$ sense. As in the linear theory, the main ingredient to prove this is a Plancherel type identity.

**Lemma 6.** *Let $F_n$ be a finite sequence of elements in the unit disc. Then with $(a, b) = \widehat{F}$ we have*

$$\int_0^1 \log |a(e^{2\pi i\theta})| \, d\theta = \int_{\mathbf{T}} \log |a(z)| = -\frac{1}{2} \sum_n \log(1 - |F_n|^2)$$

Remark 1: Observe that the integrand $\log |a(z)|$ on the left - hand side is positive, as is each summand $-\log(1 - |F_n|^2)$ on the right - hand side. Thus this equation has the flavor of a norm identity.

Exercise: prove that in lowest order approximation (quadratic) this becomes the usual Plancherel identity.

Remark 2: This formula appears at least as early as in a 1936 paper by Verblunsky [**24**] (page 291).

**Proof.** Since $F_n$ is a finite sequence, $a = a_\infty$ is a polynomial in $z^{-1}$ with constant term

$$\prod_n (1 - |F_n|^2)^{-1/2}$$

and non-vanishing in $D^*$ by Lemma 3.

Thus we have

$$\int_{\mathbf{T}} \log(a(z)) = \log(a(\infty)) = -\frac{1}{2}\sum_n \log(1-|F_n|^2)$$

Since the right-hand side is real, we also have

$$\int_{\mathbf{T}} \log|a(z)| = -\frac{1}{2}\sum_n \log(1-|F_n|^2)$$

□

Let $l^2(\mathbf{Z}_{\geq 0}, D)$ be the space of sequences supported on the nonnegative integers ("right half-line") with values in $D$. We now proceed to describe the space $\mathbf{H}$ which we will show is the range of the nonlinear Fourier transform on $l^2(\mathbf{Z}_{\geq 0}, D)$.

Consider the space $\mathbf{K}$ of all measurable $SU(1,1)$ functions on the circle with

$$\int_{\mathbf{T}} \log|a(z)| < \infty \tag{18}$$

We can embed this space into the space $L^1(\mathbf{T}) \times L^2(\mathbf{T}) \times L^2(\mathbf{T})$ mapping the function $(a,b)$ to the function $(\log|a|, b/|a|, a/|a|)$. Clearly by our assumption on the space $\mathbf{K}$, $\log|a|$ is in $L^1(\mathbf{T})$, while $b/|a|$ is an essentially bounded measurable function because $(b,a)$ is in $SU(1,1)$ almost everywhere, and $a/|a|$ is also an essentially bounded measurable function.

This embedding is indeed injective, since we can recover the modulus of $a$ almost everywhere from $\log|a|$, we can recover the phase of $a$ almost everywhere from $a/|a|$, thus we can recover $a$ almost everywhere. Then we can recover $b$ almost everywhere from $b/|a|$.

We endow the space $\mathbf{K}$ with the inherited metric. Thus the distance between two functions $(a,b)$ and $(a',b')$ is given by

$$\int_{\mathbf{T}} \log|a| - \log|a'| + \left(\int_{\mathbf{T}} |b/|a| - b'/|a'||^2\right)^{1/2} + \left(\int_{\mathbf{T}} |a/|a| - a'/|a'||^2\right)^{1/2}$$

Indeed, this makes $\mathbf{K}$ a complete metric space. To see this, it is enough to show that the image of the embedding is closed, because the space $L^1 \times L^2 \times L^2$ is complete. However, the image is the subspace of all functions $(f,g,h)$ such that $f$ is real and nonnegative almost everywhere, $g$ satisfies $|ge^f|^2 + 1 = |e^f|^2$ almost everywhere, and $h$ has values in $\mathbf{T}$ almost everywhere. Any limit of a sequence in this subspace satisfies the same constraints almost everywhere, and thus the subspace is closed.

We observe that in the above definition of the metric we could have used quotients of the type $b/a - b'/a'$ instead of $b/|a| - b'/|a'|$ and obtained an equivalent metric. This is because

$$|b/|a| - b'/|a'|| = |(b/a)(a/|a|) - (b/a)(a'/|a'|) + (b/a)(a'/|a'|) - (b'/a')(a'/|a'|)|$$
$$\leq |a/|a| - a'/|a'|| + |(b/a) - (b'/a')|$$

and similarly

$$|b/a - b'/a'| \leq |(b/|a|) - (b'/|a'|)| + |a/|a| - a'/|a'||$$

Here we have used $|b|/|a| < 1$.

Likewise, we could have used quotients $b/a^*$ in the definition of the distance. If $G$ denotes the group $SU(1,1)$ and $K$ is the compact subgroup of diagonal elements, then $b/a$ parameterizes the left residue classes $K\backslash G$ while $b/a^*$ parameterizes the right residue classes $G/K$.

Some calculations to follow will be slightly simplified if we note that we can work with quasi-metrics rather than metrics. A quasi-metric on a space is a distance function with definiteness, $d(x,y)=0$ implies $x=y$, symmetry, $d(x,y)=d(y,x)$, and a modified triangle inequality

$$d(x,z) \leq C \max(d(x,y), d(y,z))$$

Just like a metric, a quasi-metric defines a topology through the notion of open balls. This topology is completely determined by its convergent sequences. Also, Cauchy sequences are defined, and one can talk about completeness of quasi-metric spaces.

Two quasi-metrics on the same space are called equivalent if there is some strictly monotone continuous function $C$ vanishing at 0 such that

$$d(x,y) \leq C(d'(x,y))$$

$$d'(x,y) \leq C(d(x,y))$$

Two equivalent quasi-metrics produce the same convergent sequences and the same Cauchy sequences.

A quasi-metric equivalent to the above metric is given by

$$\mathrm{dist}((a,b),(a',b')) =$$

$$\int_{\mathbf{T}} |\log|a| - \log|a'|| + \int_T |b/|a| - b'/|a'||^2 + \int |a/|a| - a'/|a'||^2$$

The space $\mathbf{K}$ is too large to be the image of the NLFT. As we shall see, the image of the NLFT lies in a subspace of $\mathbf{K}$ on which the phase of $a$ does not contain any information other than that already contained in $\log|a|$.

More precisely, let $\mathbf{L}$ be the subspace of $\mathbf{K}$ consisting of all pairs $(a,b)$ such that $a$ is the boundary value of an outer function — also denoted by $a$ — on $D$ that is positive at $\infty$.

The outerness condition together with positivity of $a$ at $\infty$ can be rephrased as

$$a/|a| = e^{-ig}$$

where

$$g = p.v. \int_{\mathbf{T}} \log|a(\zeta)| \mathrm{Im}\,(\frac{\zeta+z}{\zeta-z})\, d\zeta$$

i.e., $g$ is the Hilbert transform of $\log|a|$. Recall that the harmonic extension of the Hilbert transform to $D$ vanishes at 0.

**Lemma 7.** *Let $F_n$ be a finite sequence of elements in $D$ and $(a,b) = \widehat{F}$. Then $(a,b) \in \mathbf{L}$.*

**Proof.** Clearly $(a,b) \in \mathbf{K}$. The function $a$ is holomorphic in a neighborhood of the closure of $D^*$ and has no zeros there. Therefore $a$ and $a^{-1}$ are in $H^\infty(D^*)$ and by Lemma 32 in the appendix the function $a$ is outer on $D^*$. □

**Lemma 8.** *The space $\mathbf{L}$ is closed in $\mathbf{K}$. The restriction of the quasi-metric of $\mathbf{K}$ to $\mathbf{L}$ is equivalent to the following quasi-metric on $\mathbf{L}$:*

$$\mathrm{dist}((a,b),(a',b')) = \int_{\mathbf{T}} |\log|a| - \log|a'|| + \int_T |b/|a| - b'/|a'||^2$$

If $(a, b)$ and $(a', b')$ are in $\mathbf{L}$, then

$$\int_{\mathbf{T}} |a/|a| - a'/|a'||^2 \leq C \int_T |\log|a| - \log|a'||$$

Namely,

$$\int |a/|a| - a'/|a'||^2$$
$$\leq \int_0^2 \lambda |\{|a/|a| - a'/|a'|| > \lambda\}| \, d\lambda$$
$$\leq \int_0^2 \lambda |\{|\mathrm{Im}\,(\log(a) - \log(a'))| > \lambda\}| \, d\lambda$$

where $\log(a)$ and $\log(a')$ are the boundary values of the branches of the logarithm which are real at $\infty$. Using the weak type 1 bound for the Hilbert transform, we can estimate the last display by

$$\leq \int_0^2 \lambda (C \| \log|a| - \log|a'| \|_1 / \lambda) \, d\lambda$$
$$\leq C \| \log|a| - \log|a'| \|_1$$

This estimate shows that the quasi-metric defined in the lemma is equivalent on $\mathbf{L}$ to the distance defined on $\mathbf{K}$. Moreover, it shows that a sequence in $\mathbf{L}$ which is convergent in $\mathbf{K}$ converges to an element in $\mathbf{L}$ and thus $\mathbf{L}$ is closed. □

Observe that on $\mathbf{L}$ we have

$$\mathrm{dist}(\mathrm{id}, (a, b)) \leq 3 \int_{\mathbf{T}} \log|a| \tag{19}$$

because $aa^* = 1 + bb^*$ and $1 - x \leq \log(1/x)$ imply

$$\int_{\mathbf{T}} |b/|a||^2 \leq 2 \int_{\mathbf{T}} \log|a|$$

Also observe that for $(a, b) \in \mathbf{L}$, knowledge of the quotient $b/|a|$ (or $b/a^*$ or $b/a$) is sufficient to recover $(a, b)$. Namely, we can recover $|a|$ almost everywhere using the formula $|a|^2 = 1 + |b|^2$. Then we can recover the argument of $a$ almost everywhere as the Hilbert transform of $\log|a|$. Then we can recover $b$ almost everywhere from $a$ and the quotient $b/|a|$ (or $b/a^*$ or $b/a$).

Define the space $\mathbf{H}$ to be the subspace of $\mathbf{L}$ of all functions such that $b/a^*$ is the boundary value of an analytic function on $D$ (also denoted by $b/a^*$) that is in the Hardy space $H^2$. Since the Hardy space $H^2$ (identified as space of functions on $\mathbf{T}$) is a closed subspace of $L^2$, we have that $\mathbf{H}$ is a closed subspace of $\mathbf{L}$.

If $F$ is a finite sequence supported on the right half-line, i.e., $F_n = 0$ for $n < 0$, then clearly $(a, b) \in \mathbf{H}$.

The following string of lemmas will prove that the nonlinear Fourier transform is a homeomorphism from $l^2(\mathbf{Z}_{\geq 0}, D)$ onto $\mathbf{H}$.

**Lemma 9.** *Let $F$ be a sequence in $l^2(\mathbf{Z}_{\geq 0}, D)$ and let $F^{(\leq N)}$ denote the truncations to $[0, N]$. Then $(a_N, b_N) = \widehat{F^{(\leq N)}}$ is a Cauchy sequence in $\mathbf{H}$.*

Remark: Once this lemma has been established, it is possible to define $\widehat{F}$ to be the limit of this Cauchy sequence.

**Proof.** We need the following auxiliary lemma:

**Lemma 10.** *For $G, G' \in \mathbf{L}$ we have*

$$\operatorname{dist}(GG', G) \leq C\operatorname{dist}(G', \mathrm{id}) + C\left[\operatorname{dist}(G, \mathrm{id})\operatorname{dist}(G', \mathrm{id})\right]^{1/2}$$

**Proof.** We have

$$\operatorname{dist}(GG', G) \sim \int_{\mathbf{T}} |\log|aa' + b\overline{b'}| - \log|a|| + \int_{\mathbf{T}} |\frac{ab' + b\overline{a'}}{aa' + b\overline{b'}} - \frac{b}{a}|^2 \tag{20}$$

Consider the first summand. We have

$$|\log|aa' + b\overline{b'}| - \log|a||$$

$$= |\log|a'| + \log|1 + (b/a)\overline{(b'/a')}(\overline{a'}/a')||$$

$$\leq |\log|a'|| + |\log|1 + (b/a)\overline{(b'/a')}(\overline{a'}/a')||$$

Upon integration over $\mathbf{T}$, the first summand is bounded by $\operatorname{dist}(G', \mathrm{id})$. On the set of all $z$ such that $(b'/a')(z) \geq 1/10$, we estimate

$$|\log|1 + (b/a)\overline{(b'/a')}(\overline{a'}/a')||$$

$$\leq |\log(1 - |b'/a'|)| \leq C|\log|1 - |b'/a'|^2| = 2C|\log|a'||$$

which again upon integration is bounded by $\operatorname{dist}(G', \mathrm{id})$. On the set of all $z$ with $(b'/a')(z) \leq 1/10$, we estimate

$$|\log|1 + (b/a)\overline{(b'/a')}(\overline{a'}/a')|| \leq C|b/a||b'/a'|$$

Upon integration over $\mathbf{T}$ and application of Cauchy-Schwarz, this is bounded by the square root of $\operatorname{dist}(G, \mathrm{id})\operatorname{dist}(G', \mathrm{id})$.

We consider the second summand on the right-hand side of (20). We claim that

$$\left|\frac{ab' + b\overline{a'}}{aa' + b\overline{b'}} - \frac{b}{a}\right| \leq C\frac{|b'|}{|a'|} + C\left|1 - \frac{a'}{|a'|}\right| \tag{21}$$

This will finish the proof, since upon taking the square and integrating, the right-hand side is bounded by $\operatorname{dist}(G', \mathrm{id})$ (here we use that the distance functions on $\mathbf{L}$ and $\mathbf{K}$ are equivalent). The claim is evident if $|b'|/|a'|$ is greater than $1/10$, since the left-hand side of (21) is bounded by 2. Assume

$$\frac{1}{10} > \frac{|b'|}{|a'|} > \frac{1}{10}\frac{|b|}{|a|}$$

Then we use triangle inequality on the left-hand side of (21) and obtain

$$\left|\frac{ab' + b\overline{a'}}{aa' + b\overline{b'}} - \frac{b}{a}\right| \leq C\left|\frac{ab' + b\overline{a'}}{aa'}\right| + C\frac{|b|}{|a|} \leq C\frac{|b'|}{|a'|}$$

Assume $|b'|/|a'| < \frac{1}{10}|b|/|a|$. Then

$$\left|\frac{ab' + b\overline{a'}}{aa' + b\overline{b'}} - \frac{b}{a}\right| \leq \left|\frac{ab'}{aa' + b\overline{b'}}\right| + \left|(\frac{b\overline{a'}}{aa' + b\overline{b'}} - \frac{b\overline{a'}}{aa'})\right| + \left|\frac{\overline{a'}}{a'} - 1\right|\frac{|b|}{|a|}$$

$$\leq C\frac{|b'|}{|a'|} + C\frac{|b'|}{|a'|} + C\left|1 - \frac{a'}{|a'|}\right|$$

This proves Lemma 10. □

We continue the proof of Lemma 9. Consider

$$\operatorname{dist}(\prod_{n=0}^{M} T_n, \prod_{n=0}^{N} T_n) = \operatorname{dist}(GG', G)$$

where

$$G = \prod_{n=0}^{N} T_n \quad \text{and} \quad G' = \prod_{n=N+1}^{M} T_n$$

By the Plancherel identity and (19) we have

$$\operatorname{dist}(G, \mathrm{id}) \leq \int \log |a_N| \leq C \sum_{n=1}^{N} \left|\log |1 - |F_n|^2|\right| \leq C$$

(since $F \in l^2$) and similarly,

$$\operatorname{dist}(G', \mathrm{id}) \leq C \sum_{n=N+1}^{M} \left|\log |1 - |F_n|^2|\right| \leq \epsilon$$

if $N > N(\epsilon)$ is chosen large enough depending on the choice of $\epsilon$. Thus, for $M > N > N(\epsilon)$, we have by Lemma 10

$$\operatorname{dist}(\prod_{n=0}^{M} T_n, \prod_{n=0}^{N} T_n) \leq C\epsilon^{1/2}$$

This shows that $\widehat{F^{(\leq N)}}$ is Cauchy in $\mathbf{H}$.

□

Thus we can define the NLFT on $l^2(\mathbf{Z}_{\geq 0}, D)$ as the limit of the NLFT of the truncated sequences. We have not yet shown any genuine continuity of the NLFT, but we will do that further below. Using Theorem 5, one can show that this definition of the NLFT coincides with the old definition on the subset $l^p(\mathbf{Z}_{\geq 0}, D)$ of $l^2(\mathbf{Z}_{\geq 0}, D)$ for $1 \leq p < 2$.

As the distance between the NLFT of the truncated sequence and the NLFT of the full sequence converges to 0, the Plancherel identity continues to hold on all of $l^2(\mathbf{Z}_{\geq 0}, D)$.

**Lemma 11.** *The NLFT is injective on* $l^2(\mathbf{Z}_{\geq 0})$.

**Proof.** We know for the finite truncations that

$$F_0 = b^{(\leq N)}(0)/a^{(\leq N)*}(0) = \int_{\mathbf{T}} b^{(\leq N)}/a^{(\leq N)*}$$

Where we used that $b/a^*$ has holomorphic extension to a neighborhood of $D$.

By convergence of the data $(a^{(\leq N)}, b^{(\leq N)})$ in $\mathbf{H}$ we see that $b^{(\leq N)}/a^{(\leq N)*}$ converges in $L^2(\mathbf{T})$ to $b/a^*$, where $(a, b)$ is the NLFT of $F$. Thus

$$F_0 = \int_{\mathbf{T}} b/a^*$$

Observe that the quotient $b/a^*$ is sufficient to determine $F_0$. This is consistent with the earlier observation that this quotient contains the full information of $(a, b) \in \mathbf{H}$.

To proceed iteratively, we need to determine the NLFT $(\tilde{a}, \tilde{b})$ of the "layer stripped" sequence $\tilde{F}$ in $l^2(\mathbf{Z}_{\geq 0}, D)$; this is defined by $\tilde{F}_n = F_n$ for $n > 0$ and $\tilde{F}_0 = 0$. More precisely, we will determine the quotient $\tilde{b}/\tilde{a}^*$ using only the quotient $b/a^*$. Write $r = b/a^*$, $\tilde{r} = \tilde{b}/\tilde{a}^*$ etc.

Using the product formula for finite sequences we calculate

$$(\tilde{a}^{(\leq N)}, \tilde{b}^{\leq N}) = (1 - |F_0|^2)^{-1/2}(1, -F_0)(a^{(\leq N)}, b^{(\leq N)})$$

$$\tilde{r}^{(\leq N)} = \frac{r^{(\leq N)} - F_0}{-\overline{F_0} r^{(\leq N)} + 1}$$

As $N$ tends to $\infty$, the equation tends in $L^2$ norm to

$$\tilde{r} = \frac{r - F_0}{-\overline{F_0} r + 1}$$

For the left-hand side this follows directly from the definitions. For the right-hand side this follows from the fact that the map

$$s \to \frac{s - F_0}{-\overline{F_0} s + 1}$$

has bounded derivative on the closure of $D$ and thus turns $L^2$ convergence of $r$ into $L^2$ convergence of $\frac{r - F_0}{-\overline{F_0} r + 1}$ (recall $F_0$ is fixed and $|F_0| < 1$.)

Thus we can calculate $\tilde{b}/\tilde{a}^*$.

By an inductive procedure (conjugate the new problem by a shift to reduce it to the old problem for sequences starting at 0) we can calculate all $F_n$. This proves injectivity.

□

The layer stripping method in the proof of this lemma can be used to obtain the following result: If $F$ is in $l^2$ and if $F^{(\leq N)}$ and $F^{(>N)}$ are the truncations to $[0, N]$ and $[N + 1, \infty)$, then

$$(a, b) = (a^{(\leq N)}, b^{(\leq N)})(a^{(>N)}, b^{(\leq N)})$$

This is clear if $N = 0$. If $N = 1$, this has been observe in the proof of the previous lemma. Then one can use induction to prove this for all $N$.

For later reference we note that

$$r_{>N} := b^{(>N)}/(a^{(>N)})^*$$

has a holomorphic extension to $D$ which vanishes to order $N + 1$ at 0.

**Lemma 12.** *The NLFT is surjective from $l^2$ onto the space* $\mathbf{H}$

**Proof.** Let $(a, b) \in \mathbf{H}$. Then $r = b/a^*$ is an analytic function in $D$ bounded by 1. Moreover, $r(0) < 1$ since the extension of $r$ to $\mathbf{T}$ is strictly less than 1 almost everywhere. Following the calculations in the proof of the previous lemma, we set formally

$$F_0 = r(0)$$

and

$$z\tilde{r} = \frac{r - F_0}{-F_0 r - 1}$$

Being the Möbius transform of a bounded analytic function, the right-hand side is still an analytic function in $D$ bounded by 1 and it vanishes at 0 by construction. Thus, by the Lemma of Schwarz, we can divide by $z$ and calculate formally $\tilde{r}$ which is again an analytic function in $D$ bounded by 1. This procedure can be iterated and gives a sequence $F_n$.

This iteration process can be applied to any analytic function bounded by 1, regardless of any further regularity of this function. This process is called Schur's algorithm after [**17**].

We show that $\tilde{r}$ as constructed above is of the form $\tilde{b}/\tilde{a}^*$ for some $(a,b) \in \mathbf{H}$. To this end, it suffices to show that $1-|\tilde{r}|^2$ (which is formally $|\tilde{a}|^{-2}$) is log integrable over $\mathbf{T}$. Then we can determine the outer function $\tilde{a}$ and calculate $\tilde{b}$. Observe that as bounded analytic function, $\tilde{r}$ is automatically in $H^2(D)$.

We have

$$1-|\tilde{r}|^2 = \frac{|\overline{F_0}r-1|^2 - |r-F_0|^2}{|\overline{F_0}r-1|^2}$$

$$= \frac{|\overline{F_0}r|^2 - |r|^2 - |F_0|^2 + 1}{|\overline{F_0}r-1|^2}$$

$$= \frac{(1-|F_0|^2)(1-|r|^2)}{|\overline{F_0}r-1|^2}$$

Observe that $\overline{F_0}r-1$ is bounded and bounded away from 0 and so is its holomorphic extension to $D$, thus its logarithm is integrable over $\mathbf{T}$ and equal to the value of the logarithm at 0:

$$\int \log(1-\overline{F_0}r) = \log(1-|F_0|^2)$$

Taking logarithms and integrating gives

$$\int \log(1-|\tilde{r}|^2) = \log(1-|F_0|^2) + \int (1-|r|^2) - 2\log(1-|F_0|^2)$$

$$\int \log(1-|\tilde{r}|^2) = -\log(1-|F_0|^2) + \int \log(1-|r|^2)$$

Discussing the signs we obtain

$$\int |\log(1-|\tilde{r}|^2)| = -|\log(1-|F_0|^2)| + \int |\log(1-|r|^2)| \tag{22}$$

Thus $\log(1-|\tilde{r}|^2)$ is integrable.

We can iterate to calculate formally $F_n$. Using (22) inductively, we obtain

$$\sum_{n=0}^{\infty} |\log(1-|F_n|^2)| \le \int |\log(1-|r|^2)| \tag{23}$$

Thus the sequence $F_n$ we calculated is in $l^2(\mathbf{Z}_{\ge 0}, D)$ and it is a candidate for the preimage of $(a,b)$ under the NLFT.

Let $(\tilde{a},\tilde{b})$ denote the NLFT of $F_n$. We will show that indeed, $(a,b) = (\tilde{a},\tilde{b})$ and we have equality in (23).

As noted earlier,

$$(\tilde{a},\tilde{b}) = (\tilde{a}^{(\le N)}, \tilde{b}^{(\le N)})(\tilde{a}^{(>N)}, \tilde{b}^{(>N)})$$

where the factors on the right-hand side are defined by the usual truncations.

Observe that by unwinding Schur's algorithm introduced above, we obtain

$$(a,b) = (\tilde{a}^{(\le N)}, \tilde{b}^{(\le N)})(a^{(>N)}, b^{(>N)})$$

where $(a^{(>N)}, b^{(>N)})$ is the unique element in $\mathbf{H}$ such that $b^{(>N)})/a^{(>N)}$ is equal to the $N$-th function in Schur's algorithm.

Thus in the expression

$$(\tilde{a},\tilde{b})^{-1}(a,b)$$

we can cancel factors to obtain

$$(\tilde{a},\tilde{b})^{-1}(a,b) = (\tilde{a}^{(>N)},\tilde{b}^{(>N)})^{-1}(a^{(>N)},b^{(>N)}) \tag{24}$$

Consider an off- diagonal entry of the right-hand side of (24):

$$(\tilde{a}^{(>N)})^* b^{(>N)} - \tilde{b}^{(>N)}(a^{(>N)})^*$$

This is a Nevanlinna function on $D$. The Taylor coefficients of this function at 0 vanish up to order $N$. By (24), this function does not actually depend on $N$, therefore all Taylor coefficients at 0 vanish. Thus the function vanishes on $D$, and so do its radial limits on $\mathbf{T}$ almost everywhere. Thus the off diagonal coefficients of

$$(\tilde{a},\tilde{b})^{-1}(a,b)$$

vanish and we have for some function $c$

$$(a,b) = (\tilde{a},\tilde{b})(c,0)$$

This function $c$ has to have modulus 1 almost everywhere on $\mathbf{T}$, since the last display can be read as identity between $SU(1,1)$ valued measurable functions on $\mathbf{T}$. Calculating a diagonal entry in the last display, we obtain

$$a = \tilde{a}c$$

Since $a$ and $\tilde{a}$ are outer, so is $c$. However, any outer function with constant modulus on $\mathbf{T}$ is constant. Moreover, $c$ is positive since $a$ and $\tilde{a}$ are positive at $\infty$. This proves $c = 1$. □

**Lemma 13.** *The NLFT is a continuous map from $l^2(\mathbf{Z}_{\geq 0}, D)$ to $\mathbf{H}$.*

**Proof.** Fix $F \in l^2(\mathbf{Z}_{\geq 0}, D)$ and choose $\epsilon > 0$. Choose $N$ very large depending on $\epsilon$.

For $F'$ close to $F$ depending on $N$, $\epsilon$ we write

$$\mathrm{dist}((a,b),(a',b')) \leq \mathrm{dist}((a,b),(a^{(\leq N)},b^{(\leq N)})) +$$
$$+\mathrm{dist}((a^{(\leq N)},b^{(\leq N)}),(a'^{(\leq N)},b'^{(\leq N)})) + \mathrm{dist}((a'^{(\leq N)},b'^{(\leq N)}),(a',b'))$$

We intend to argue that all the terms on the right-hand side are less than $\epsilon/3$.

By the definition of $(a,b)$, the first term can be made small by choosing $N$ large enough. Likewise the third term can be made small, since the distance between the truncation and the full Fourier transform depends only on the $l^2$ norm of the sequence $F'$ and the $l^2$ norm of the tail of this sequence, which can be both controlled by choosing $N$ large enough and $F'$ close enough to $F$. Thus it remains to control the middle term.

Consider the space of $D$-valued sequences on $[0,N]$ with the $l^2$ norm. Since the space is finite dimensional, the $l^2$ norm is equivalent to the $l^1$ norm.

Observe that for two matrices $(a,b)$ and $(a',b')$ in $SU(1,1)$ we have

$$|\log|a| - \log|a'|| \leq |a-a'| \leq \|(a,b)-(a',b')\|_{op}$$

and

$$|b/|a| - b'/|a'|| \leq |b-b'| \leq \|(a,b)-(a',b')\|_{op}$$

Thus, if $F'^{(\leq N)}$ is sufficiently close to $F^{(\leq N)}$ w.r.t. $l^2$ and thus $l^1$, then

$$\sup_z \left|\log|a^{(\leq N)}| - \log|a'^{(\leq N)}|\right| + \sup_z \left|b^{(\leq N)}/|a^{(\leq N)}| - b'^{(\leq N)}/|a'^{(\leq N)}|\right|$$

is small, and thus

$$\mathrm{dist}((a',b'),(a^{(\leq N)},a^{(\leq N)}))$$

is small. □

We remark that this proof does not give any uniform continuity. The weak point in the argument is the comparison of the $l^2$ with the $l^1$ norm of a finite sequence without any good control over the length of the sequence.

**Lemma 14.** *The inverse of the NLFT is a continuous map from* $\mathbf{H}$ *to* $l^2(\mathbf{Z}_{\geq 0}, D)$.

**Proof.** We first prove that all $F_n$ depend continuously on $(a, b)$. This is clear for $F_0$ since

$$F_0 = \int b/a^*$$

and the integral is continuous in the $L^2$ norm of $b/a^*$.

To use induction, we need to show that (see the proof of injectivity)

$$r \to \frac{r - F_0}{-\overline{F_0} r + 1}$$

is a jointly continuous mapping of $F_0 \in D$ and $r \in H^2$ into $H^2$. This follows from the fact that the Möbius transform with $F_0$ provides a Lipschitz distortion on the closed unit disc, and the distortion depends continuously on $F_0$.

Now let $(a, b)$ be given and let $F$ be its inverse NLFT. Given $\epsilon$, we can find $N$ very large so that the $[N, \infty)$ tail of $F$ is very small. Let $(a', b')$ be close to $(a, b)$ and let $F'$ be the inverse NLFT. Then we can assume for all $n \leq N$

$$\log(1 - |F_n|^2) - \log(1 - |F_n'|^2)$$

is much smaller than $\epsilon$, by continuity of all $F_n$ individually.

Next we have

$$\sum_{n>N} |\log(1 - |F_n'|^2)| = \sum_n |\log(1 - |F_n'|^2)| - \sum_{n \leq N} |\log(1 - |F_n'|^2)|$$

$$\leq \sum_n |\log(1 - |F_n|^2)| - \sum_{n \leq N} |\log(1 - |F_n|^2)| + \epsilon \leq 2\epsilon$$

Here we have used that $(a, b)$ and $(a', b')$ are close and therefore, by the Plancherel identity, the quantities

$$\sum_n \log(1 - |F_n|^2)$$

$$\sum_n \log(1 - |F'_n|^2)$$

are close. Also we have used the previously observed continuity for individual $F_n$.

Now it is straight forward to obtain

$$\sum_n |\log(1 - |F_n|^2) - \log(1 - |F'_n|^2)| < 4\epsilon$$

by considering separately $n > N$ and $n \leq N$. □

## 2. Higher order variants of the Plancherel identity

The main ingredient in the $l^2$ theory of the nonlinear Fourier transform described in the previous section is the Plancherel identy

$$\int_{\mathbf{T}} \log|a| = -\frac{1}{2}\sum_{n\in\mathbf{Z}} \log(1-|F_n|^2)$$

Both sides of the identity are equal to $a(\infty)$, which on the left-hand side is expressed as a Cauchy integral and on the right-hand side is expressed in terms of the sequence $F$ by solving explicitly the recursion for $a(\infty)$.

There are higher order identities of this type, which arise from calculating higher derivatives of $\log(a)$ at $\infty$. These identities are nonlinear analogues of Sobolev identities, i.e., identities between expressions for Sobolev norm of a function in terms of the function itself and in terms of its Fourier transform.

We discuss $a^*$ instead of $a$. The contour integral for $\log(a^*)^{(k)}(0)$ can easily be written in closed form:

$$\log(a^*)(z) = \int_{\mathbf{T}} \frac{\zeta+z}{\zeta-z}\log|a|(\zeta) = \int_T \left[\frac{2\zeta}{\zeta-z}-1\right]\log|a|(\zeta)$$

Taking derivatives in $z$ we obtain for $k>0$

$$\log(a^*)^{(k)}(z) = \int_T \frac{2\zeta k!}{(\zeta-z)^{k+1}}\log|a|(\zeta)$$

$$\log(a^*)^{(k)}(0) = \int_T \frac{2k!}{\zeta^k}\log|a|(\zeta)$$

Solving the recursion in terms of the $F_n$ is harder to do in closed form. Such formulae are stated in Case's paper [**2**], see also [**13**].

We shall calculate only the cases $k=1,2$. For an application to the theory of orthogonal polynomials, see [**7**].

**Lemma 15.** *For $F$ a square summable sequence we have*

$$2\int_{\mathbf{T}} z^{-1}\log|a|(z) = \sum_n \overline{F_n}F_{n+1}$$

$$4\int_{\mathbf{T}} z^{-2}\log|a|(z) = -\sum_n (\overline{F_n}F_{n+1})^2 + 2\sum_n \overline{F_n}(1-|F_{n+1}|^2)F_{n+2}$$

**Proof.** We first reduce to the case of compactly supported seqeunces $F$ by approximating an arbitrary sequence by its truncations $F_n$. At least for half infinite sequences we have already seen that $\log|a_n|$ converges to $\log|a|$ in $L^1$ norm, and in the next section we will discuss and establish the same fact for general sequences in $l^2(D,\mathbf{Z})$. Thus the left-hand side of each of the identities in the lemma is well approximated by the truncations. Also the right-hand side clearly converges if $F$ is in $l^2$. Thus it suffices to show the identities for compactly supported $F$.

Assume $F$ is compactly supported. We expand the product

$$\prod_{n\in\mathbf{Z}} (1-|F_n|^2)^{-1/2}[(1,0)+(0,F_nz^n)]$$

Only the terms of even order in $F$ contribute to the diagonal elements and thus

$$a^*(z)\prod_n (1-|F_n|^2)^{1/2}$$

$$= 1 + \sum_{n_1 < n_2} \overline{F_{n_1}} F_{n_2} z^{n_2 - n_1} + \sum_{n_1 < n_2 < n_3 < n_4} \overline{F_{n_1}} F_{n_2} \overline{F_{n_3}} F_{n_4} z^{n_2 - n_1} z^{n_4 - n_3} + \dots$$

Since $n_1 < n_2$ and $n_3 < n_4$ etc., we see that the bilinear term in $F$ is has lowest order $z$ and the four-linear term has lowest order $z^2$, while all other terms have order at least $z^3$. Thus for the purpose of calculating the first two derivatives at $\infty$ we only need to consider the terms that are explicitly written.

Indeed, we have

$$a^*(z) \prod_n (1 - |F_n|^2)^{1/2} = 1 + \sum_n \overline{F_n} F_{n+1} z$$

$$+ \sum_n \overline{F_n} F_{n+2} z^2 + \sum_{n_1 + 1 < n_2} \overline{F_{n_1}} F_{n+1} \overline{F_{n_2}} F_{n_2+1} z^2 + O(z^3)$$

Considering the case $k = 1$ we have

$$\log(a^*)'(0) = \frac{(a^*)'(0)}{a^*(0)} = \sum_n \overline{F_n} F_{n+1}$$

This proves the first identity of Lemma 15.

Considering the case $k = 2$ we have

$$\log(a^*)''(0) = \frac{(a^*)''(0)}{a(0)} - \frac{(a^*)'(0)^2}{a^*(0)^2}$$

$$= 2 \sum_n \overline{F_n} F_{n+2} + 2 \sum_{n_1 + 1 < n_2} \overline{F_{n_1}} F_{n+1} \overline{F_{n_2}} F_{n_2+1} - [\sum_n \overline{F_n} F_{n+1}]^2$$

$$= - \sum_n (\overline{F_n} F_{n+1})^2 + 2 \sum_n \overline{F_n} (1 - |F_{n+1}|^2) F_{n+2}$$

This proves the second identity of Lemma 15. □

# LECTURE 3
# The nonlinear Fourier transform on $l^2(\mathbf{Z})$

## 1. The forward NLFT on $l^2(\mathbf{Z})$

We have defined the nonlinear Fourier transform for square summable sequences supported on the nonnegative integers. Indeed, we have shown that it is a homeomorphism onto $\mathbf{H}$, the space of all measurable $SU(1,1)$ valued function $(a,b)$ such that $a$ has an outer extension to $D^*$, $a(\infty)>0$, and $b/a^*$ has a holomorphic extension to $D$ which is in the Hardy space $H^2(D)$.

Demanding that property (8) of Lemma 1 continues to hold for infinite sequences, we define for $F$ supported on $n \leq 0$:

$$\overbrace{F}(z) := (a^*(z^{-1}), b(z^{-1}))$$

where $(a,b)$ is the Fourier transform of the reflected sequence $\tilde{F}$ with $\tilde{F}_n = F_{-n}$.

It is clear that the nonlinear Fourier transform thus defined is a homeomorphism from $l^2(\mathbf{Z}_{\leq 0}, D)$ to $\mathbf{H}^*$, the latter denoting the space of all $SU(1,1)$ valued measurable functions $(a,b)$ such that $a$ has an outer extension to $D^*$, $a(\infty)>0$, and $b/a$ has a holomorphic extension to $D^*$ which is in the Hardy space $H^2(D^*)$.

Let $\mathbf{H}_0^*$ be the space of all elements in $\mathbf{H}^*$ such that $b(\infty)=0$. By the shifting property (5) of Lemma 1, it is easy to deduce that this space is the homeomorphic image of $l^2(\mathbf{Z}_{\leq -1}, D)$.

If $F_n$ is any square summable sequence in $l^2(\mathbf{Z}, D)$, then we can cut it as

$$F_n = F_n^{(\leq -1)} + F_n^{(\geq 0)}$$

where $F_n^{(\leq -1)} = 0$ for $n \geq 0$ and $F_n^{(\geq 0)} = 0$ for $n \leq -1$. Then we define a measurable $SU(1,1)$ valued function on $\mathbf{T}$ by

$$\overbrace{F} := \overbrace{F^{(\leq -1)}}\,\overbrace{F^{(\geq 0)}} \tag{25}$$

in accordance with property (6) of Lemma 1. We shall use the suggestive notation

$$(a,b) = (a_-, b_-)(a_+, b_+)$$

for (25).

It is easy to verify that the NLFT defined by (25) on $l^2(\mathbf{Z}, D)$ satisfies the properties of Lemma 1. The properties of Lemma 1 imply that the exact location of the cut we used in (25) is not relevant for the definition.

As noted previously, the definition of the NLFT on $l^2(\mathbf{Z}_{\geq 0}, D)$ was consistent with earlier definitions on the subset $l^p(\mathbf{Z}_{\geq 0}, D)$ for $p < 2$. Passing from the half-line to the full line, since definition (25) and the old definition of the NLFT on $l^p(\mathbf{Z}, D)$ are consistent with Lemma 1, the two definitions coincide on $l^p(\mathbf{Z}, D)$.

**Lemma 16.** *The NLFT is continuous from $l^2(\mathbf{Z}, D)$ to $\mathbf{L}$. The Plancherel identity*

$$\int_{\mathbf{T}} \log |a(z)| = -\frac{1}{2} \sum_n \log(1 - |F_n|^2)$$

*holds.*

**Proof.** Recall that $\mathbf{L}$ is the space of all $SU(1,1)$ valued measurable functions $(a, b)$ such that $a$ has an outer extension to $D^*$ and $a(\infty) > 0$.

First we check that the nonlinear Fourier transform indeed maps to $\mathbf{L}$. We need to verify that $a$ has an outer extension to $D^*$ and that $a(\infty) > 0$. But we have

$$a = a_- a_+ + b_- b_+^*$$

and all functions on the right-hand side extend holomorphically to $D^*$ with

$$b_-(\infty) = 0$$

Therefore

$$a(\infty) = a_-(\infty) a_+(\infty) > 0$$

Moreover,

$$a = a_- a_+ \left(1 + \frac{b_-}{a_-} \frac{b_+^*}{a_+}\right)$$

and the first two factors on the right-hand side are outer. The last factor has positive real part on on $D^*$ because the extensions of $b_-/a_-$ and $b_+^*/a_+$ to $D^*$ are bounded by 1. Thus the Herglotz representation theorem applies to the last factor, which then can be seen to be in $H^p(D^*)$ for all $p < 1$. The reciprocal of the last factor has also positive real part and is also in $H^p(D^*)$ for all $p < 1$. Thus the last factor is an outer function.

The proof of continuity invokes Lemma 10. Given $F$, we can use the independence of the cut in Definition (25) to cut $F$ and any nearby $F'$ at a very large integer $N$ (depending only on $F$), so that the tail to the right of $N$ of both $F$ and $F'$ is negligible by the Plancherel identity and Lemma 10. Then we can apply continuity of the nonlinear Fourier transform on the (shifted) half line to show that the parts of $F$ and $F'$ to the left of $N$ have nearby nonlinear Fourier transforms.

The same argument also proves the Plancherel identity. □

We now observe

**Lemma 17.** *The nonlinear Fourier transform is not injective on $l^2(\mathbf{Z}, D)$.*

**Proof.** We claim that

$$(a, b) = \left(\frac{2z}{z-1}, \frac{z+1}{z-1}\right)$$

is in $\mathbf{H} \cap \mathbf{H}^*$. Therefore it has nonzero preimages in $l^2(\mathbf{Z}_{\geq 0})$ and in $l^2(\mathbf{Z}_{\leq 0})$, and since these preimages are not finite sequences ($a$ is not a Laurent polynomial), these two preimages are necessarily distinct members of $l^2(\mathbf{Z}, D)$.

It remains to prove the claim. We observe

$$a(z)a^*(z) - b(z)b^*(z)$$

$$= \frac{(2z)(2z^{-1})}{(z-1)(z^{-1}-1)} - \frac{(z+1)(z^{-1}+1)}{(z-1)(z^{-1}-1)}$$

$$= \frac{4}{-z+2-z^{-1}} - \frac{z+2+z^{-1}}{-z+2-z^{-1}} = 1$$

The function $a$ is outer on $D^*$ since it is in $H^p(D^*)$ for all $p < 1$ and its reciprocal is in $H^\infty(D^*)$. We also have $a(\infty) = 2 > 0$.

Moreover, both

$$\frac{b(z)}{a(z)} = \frac{z+1}{2z}$$

and

$$\frac{b^*(z)}{a(z)} = -\frac{z+1}{2z}$$

are holomorphic in $D^*$ and in $H^2(D^*)$. This proves the claim. □

We now discuss the inverse problem, i.e., finding a (sometimes not unique) sequence $F$ whose nonlinear Fourier transform is a given $(a, b)$.

Given data $(a, b) \in \mathbf{L}$, we need to factorize it

$$(a_-, b_-)(a_+, b_+) = (a, b)$$

with $(a_-, b_-) \in \mathbf{H}_0^*$ and $(a_+, b_+) \in \mathbf{H}$. Any such factorization is in bijective correspondence to a sequence $F \in l^2(\mathbf{Z}, D)$ whose truncations satisfy

$$\widehat{F_{\leq -1}} = (a_-, b_-)$$

$$\widehat{F_{\geq 0}} = (a_+, b_+)$$

Thus the inverse problem for the nonlinear Fourier transform is a matrix factorization problem with (mainly but not exclusively) holomorphicity conditions on the matrix factors.

Observe that the corresponding linear problem is the decomposition of a function $f \in L^2(\mathbf{T})$ as the sum of a function in the Hardy space $H^2$ and a function in the conjugate Hardy space $\overline{H}_0^2$ (where the index 0 stands for functions with mean zero).

Finding a factorization of a matrix function on $\mathbf{T}$ into a product of two matrix functions, one extending to $D$ and one extending to $D^*$ is called a Riemann-Hilbert problem.

Our factorization is a somewhat twisted Riemann-Hilbert problem, because the matrices both have entries which extend to $D$ and $D^*$. Moreover, the factorization problem is constrained in that there are algebraic relations between the matrix entries and there is an outerness condition on $a_-, a_+$ and a normalization condition at $\infty$.

However, the factorization problem can be reduced to a more genuine Riemann-Hilbert problem by the following algebraic manipulations. The factorization equation together with the determinant condition can be rewritten as

$$\begin{pmatrix} a_+{}^* & -b_+ \\ b_-^* & a_-^* \end{pmatrix} \begin{pmatrix} a_+ & b_+ \\ b_+^* & a_+^* \end{pmatrix} = \begin{pmatrix} 1 & 0 \\ b^* & a^* \end{pmatrix}$$

Here the second row comes from the factorization problem while the first row comes from

$$(a_+, b_+)^{-1}(a_+, b_+) = (1, 0)$$

Similarly we obtain

$$\begin{pmatrix} a_+^* & -b_+ \\ -b_+^* & a_+ \end{pmatrix} \begin{pmatrix} a_+ & -b_- \\ b_+^* & a_- \end{pmatrix} = \begin{pmatrix} 1 & -b \\ 0 & a \end{pmatrix}$$

Multiplying the two equations and using the determinant condition on $(a, b)$ gives

$$\begin{pmatrix} a_+{}^* & -b_+ \\ b_-{}^* & a_-{}^* \end{pmatrix} \begin{pmatrix} a_+ & -b_- \\ b_+{}^* & a_- \end{pmatrix} = \begin{pmatrix} 1 & -b \\ b^* & 1 \end{pmatrix}$$

In the last equation, all entries of the first factor on the left-hand side extend to $D$ while all entries of the second factor on the left-hand side extend to $D^*$ (the function $b_-$ is in addition required to vanish at $\infty$). The Riemann-Hilbert problem is still constrained in that the entries of the two matrices on the left are dependent, however the constraints can be subsumed in the statement that the factorization should be invariant under the map

$$T : \begin{pmatrix} a & b \\ c & d \end{pmatrix} \to \begin{pmatrix} a^* & -c^* \\ -b^* & d^* \end{pmatrix}$$

This map reverses the order of multiplication, $T(G)T(G') = T(G'G)$, and one can easily check that this symmetry produces all algebraic constraints between the two factors.

Observe that while $a$ no longer appears explicitly in the factorization problem, $a^*$ and $a$ formally coincide with the determinants of the two factors and thus taking determinants everywhere we formally obtain the equation

$$a^* a = 1 + bb^*$$

Observe that for any solution of this Riemann-Hilbert problem we can produce another solution by multipying the first factor by a function $c$ of modulus 1 on $\mathbf{T}$ with holomorphic extension to $D$ and the second factor by $c^*$. To obtain any hope for uniqueness, one has to make the additional analytic assumption that say all functions of the first matrix are in $N^+(D)$ (as defined in the appendix) and all entries of the second matrix are in $N^+(D^*)$ and that the determinants of the two matrices are outer on $D$ and $D^*$ respectively. A sharper constraint is to require that the diagonal entries of the two factors are outer functions on $D$ and $D^*$ respectively.

Then the only obvious ambiguity left is a scalar factor, which can be normalized by requiring the determinants of both matrices to be positive at 0 and $\infty$ respectively.

## 2. Existence and uniqueness of an inverse NLFT for bounded $a$

We shall now prove existence and uniqueness of the solution to the factorization problem under the additional assumption that $a$ is bounded. We shall use the original formulation rather than the more genuine Riemann-Hilbert problem.

**Lemma 18.** *If $(a, b) \in \mathbf{L}$ and in addition $a$ is a bounded function, then there is a unique $F_n \in l^2(\mathbf{Z}, D)$ such that*

$$\widehat{F_n} = (a, b)$$

**Proof.** By the half-line theory it suffices to find and show uniqueness of a decomposition

$$(a_-, b_-)(a_+, b_+) = (a, b) \tag{26}$$

such that the factors on the left-hand side are in $\mathbf{H}_0^*$ and $\mathbf{H}$ respectively.

We first prove the following, which does not require $a$ to be bounded.

**Lemma 19.** *Let $(a, b) \in \mathbf{L}$. For any factorization of the Riemann-Hilbert problem*

$$(a_-, b_-)(a_+, b_+) = (a, b)$$

*with $(a-, b-) \in \mathbf{H}_0^*$ and $(a_+, b_+) \in \mathbf{H}$, we have that $a_-/a$ and $a_+/a$ are functions in $H^2(D^*)$.*

**Proof.** As $a_-/a$ and $a_+/a$ are outer, it suffices to show that the boundary values of these functions on $\mathbf{T}$ are in $L^2$.

The Riemann-Hilbert problem gives

$$a = a_- a_+[1 - (b_-/a_-)(b_+^*/a_+)]$$

or equivalently,

$$a_- a_+/a = [1 - (b_-/a_-)(b_+^*/a_+)]^{-1} \tag{27}$$

We first show that the real part of the right-hand side is in $L^1(\mathbf{T})$.

This function extends to $D^*$ with positive real part, because the quotients $b_-/a_-$ and $b_+^*/a_+$ are strictly bounded by 1 on $D^*$. By the Theorem 20 of Herglotz discussed in the appendix, the real part of (27) is the harmonic extension of a positive measure. Almost everywhere on $\mathbf{T}$, the real part of the function $a_- a_+/a$ coincides with the density of the absolutely continuous part of this measure and is thus in $L^1(\mathbf{T})$.

As

$$\mathrm{Re}\,(\frac{a_- a_+}{a}) + \mathrm{Re}\,(\frac{b_- b_+^*}{a}) = 1$$

from the Riemann-Hilbert factorization, we also have that

$$\mathrm{Re}\left[\frac{a_- a_+}{a} - \frac{b_- b_+^*}{a}\right] \tag{28}$$

is absolutely integrable.

The Riemann-Hilbert factorization can be rewritten in terms of the equations

$$(a_-, b_-) = (a, b)(a_+^*, -b_+)$$
$$(a_+, b_+) = (a_-^*, -b_-)(a, b)$$

These give

$$a_- = a a_+^* - b b_+^*$$
$$b_+ = a_-^* b - b_- a^*$$

Which in turn give

$$a_+^* = \frac{a_-}{a} + \frac{b b_+^*}{a}$$
$$b_- = -\frac{b_+}{a^*} + \frac{a_-^* b}{a^*}$$

Thus we can write for (28) on $\mathbf{T}$:

$$\mathrm{Re}\left[\frac{a_- a_-^*}{a a^*} + \frac{a_- b^* b_+}{a a^*} + \frac{b_+^* b_+}{a a^*} - \frac{b_+^* a_-^* b}{a a^*}\right]$$
$$= \mathrm{Re}\left[\frac{a_- a_-^*}{a a^*} + \frac{b_+^* b_+}{a a^*}\right]$$

In the last line we have cancelled two terms inside the brackets which added up to a purely imaginary quantity on $\mathbf{T}$. From integrability of the last line, we observe that

$$a_-/a \in L^2(\mathbf{T})$$

and also, since $a_+a_+^* = 1 + b_+b_+^*$ and $|a| > 1$ on $\mathbf{T}$,

$$a_+/a \in L^2(\mathbf{T})$$

This proves the lemma. □

We now prove the uniqueness part of Lemma 18.

By Lemma 19, it suffices to prove uniqueness under the additional assumption that $a_+, b_+^*, a_-, b_-$ are in $H^2(D^*)$.

We rewrite the Riemann-Hilbert problem as

$$(a_-, b_-) = (a, b)(a_+{}^*, -b_+)$$

Since $a$ is nonvanishing on $\mathbf{T}$, we can rewrite the second equation as

$$\frac{b_-}{a} = -b_+ + \frac{b}{a}a_+$$

Let $P_D$ be the orthogonal projection from $L^2(T)$ to $H^2(D)$. Then the previous display implies

$$b_+ = P_D\left(\frac{b}{a}a_+\right)$$

since $b_+$ is already in $H^2(D)$ and $b_-/a$ is in $H^2(D^*)$ with vanishing constant term.

Next, we have again from the Riemann-Hilbert factorization

$$\frac{a_-^*}{a^*} = a_+ - \frac{b^*}{a^*}b_+$$

Applying the orthogonal projection $P_{D^*}$ from $L^2(T)$ to $H^2(D^*)$, we obtain

$$a_+ = \frac{a_-(\infty)}{a(\infty)} + P_{D^*}\left(\frac{b^*}{a^*}b_+\right)$$

Here we have used that $a_+$ is already in $H^2(D^*)$ and the quotient $a_-^*/a^*$ is in $H^2(D)$ and thus its $P_D$ projection is equal to its constant term, which is real and equal to $a_-(\infty)/a(\infty)$.

Observe that evaluating the extension of

$$a = a_-a_+ + b_-b_+^*$$

at $\infty$ gives

$$a(\infty) = a_-(\infty)a_+(\infty) + 0$$

Thus we can can rewrite the expression for $a_+$ as

$$a_+ = \frac{1}{a_+(\infty)} + P_{D^*}\left(\frac{b^*}{a^*}b_+\right)$$

For any constant $c$, the affine linear map

$$(A, B) \mapsto \left(c + P_{D^*}\left(\frac{b^*}{a^*}B\right), P_D\left(\frac{b}{a}A\right)\right)$$

is a contraction in $L^2(\mathbf{T}) \oplus L^2(\mathbf{T})$ (Hilbert space sum). Namely, $P_D$ and $P_{D^*}$ have norm 1 in $L^2(\mathbf{T})$, while multiplication by $b/a$ or $b^*/a^*$ have norm strictly less than 1 in $L^2(\mathbf{T})$. Here we use that $a$ is bounded and thus

$$\left|\frac{b}{a}\right| \leq \left|1 - \frac{1}{|a|^2}\right|^{1/2} \leq 1 - \epsilon$$

Therefore, by the contraction mapping principle, this map has a unique fixed point $(A_c, B_c)$. Indeed, by linearity, this fixed point is $(cA_1, cB_1)$.

From the above it follows that $(a_+, b_+)$ is equal to this unique fixed point for some constant $c$. To prove uniqueness of $(a_+, b_+)$, it therefore suffices to show that we can determine $c$ uniquely.

The phase of the constant $c$ is determined by the requirement

$$a_+(\infty) = cA_1(\infty) > 0$$

The modulus of $c$ can then be determined by

$$cA_1(\infty) = \frac{1}{cA_1(\infty)} + P_D\left(\frac{b}{a}b_+^*\right)(\infty)$$

and thus

$$cA_1(\infty)^2 = 1 + P_D(\frac{b}{a}b_+^*)(\infty)cA_1(\infty)$$

The second summand on the right-hand side is necessarily real and positive since the left-hand side is larger than 1. This gives a quadratic equation for $cA_1(\infty)$ with a positive and a negative solution. Since $cA_1(\infty)$ is necessarily positive, it is therefore uniquely determined. Thus we can recover $(a_+, b_+) = (cA_1, cA_1)$ completely from $(a, b)$. By matrix division, we also obtain $(a_-, b_-)$. Thus the solution to the Riemann Hilbert problem is unique.

It remains to prove existence of the solution to the Riemann-Hilbert problem. In the next section, we will prove existence without assuming boundedness of $a$. However, the proof of existence for bounded $a$ is much easier. Therefore we choose to present it here.

Consider again the above fixed point equation and let $(A, B) \in H^2(D^*) \oplus H^2(D)$ be the unique solution for $c = 1$.

Observe that by interpolation, the linear map is also a contraction on $H^{2+\epsilon} \oplus H^{2+\epsilon}$ for small $\epsilon$. Namely, the map is bounded in any space $H^p \oplus H^p$ with $2 < p < \infty$. The operator norm may be large for any fixed $p$, but interpolating this estimate with the estimate for $H^2 \oplus H^2$, where the operator norm is less than 1, gives for sufficiently small $\epsilon$ an operator norm on $H^{2+\epsilon} \oplus H^{2+\epsilon}$ which is still less than 1. Hence the unique solution in $H^2 \oplus H^2$ is actually in the subspace $H^{2+\epsilon} \oplus H^{2+\epsilon}$, since this subspace also contains a solution by the contraction mapping principle. Hence the regularity of the solution to the Riemann-Hilbert problem will be slightly better than Lemma 19 suggests. We shall not need this extra regularity in the current proof.

We claim that the function

$$AA^* - BB^*$$

is constant on $\mathbf{T}$. Since it is manifestly real, it suffices to show that it is in the Hardy space $H^1(D)$ (being in a Hardy space $H^p$ with $p \geq 1$ makes the linear Fourier coefficients supported on a half-line, while being real makes the moduli of the Fourier coefficients symmetric about 0). Use the fixed point equation to write the function as

$$= A[1 + P_D(\frac{b}{a}B^*)] + B^* P_D(\frac{b}{a}A)$$

$$= A[1 - (\mathrm{id} - P_D)(\frac{b}{a}B^*)] + B^*(\mathrm{id} - P_D)(\frac{b}{a}A)$$

Here the two terms involving the identity operators that have artificially been inserted are negatives of each other. Observe that $\mathrm{id} - P_D$ is projecting onto the space $H^2_0(D^*)$ of functions in the Hardy space $H^2(D^*)$ with vanishing constant coefficient.

Thus the entire last displayed expression is an element in $H^1(D^*)$, since it is the sum of products of functions in $H^2(D^*)$.

Moreover, we observe that the constant coefficient of this expression is that of $A$:

$$\int_{\mathbf{T}} AA^* - BB^* = \int_{\mathbf{T}} A = A(\infty)$$

Thus the constant coefficient of $A$ is real. Indeed, it is positive, as we see from the following calculation:

$$\int_{\mathbf{T}} AA^* + BB^*$$

$$= \int_{\mathbf{T}} A(1 + P_D(\frac{b}{a}B^*)) + B^* P_D(\frac{b}{a}A)$$

$$= \int_{\mathbf{T}} A(1 + \frac{b}{a}B^*) + B^*\frac{b}{a}A$$

In the last line we have dropped the projection operators, because the operands are integrated against functions in the Hardy space $H^2(D^*)$. However, estimating the last display using $|b/a| \leq 1$ on $\mathbf{T}$, we obtain:

$$\int_{\mathbf{T}} AA^* + BB^* \leq \int_{\mathbf{T}} A + 2\int_{\mathbf{T}} |A||B|$$

or

$$\int_{\mathbf{T}} (|A| - |B|)^2 \leq A(\infty)$$

Thus the constant coefficient of $A$ is nonnegative.

Indeed, in the above string of inequlities, identity holds only if $|A| = |B| = 0$ almost everywhere on $\mathbf{T}$, since $|b/a|$ is strictly less than 1 almost everywhere. However, $A = B = 0$ is inconsistent with the fixed point equation. Therefore, we have strict inequality and the constant coefficient of $A$ is strictly positive.

Set

$$a_+(z) := A(z)[A(\infty)]^{-1/2}$$

$$b_+(z) := B(z)[A(\infty)]^{-1/2}$$

Then

$$|a_+|^2 - |b_+|^2 = 1$$

almost everywhere on $\mathbf{T}$ and $a_+ \in H^2(D^*)$ and $b_+ \in H^2(D)$.

Now we can define $a_-$ and $b_-$ by

$$(a_-, b_-) := (a, b)(a_+{}^*, -b_+)$$

but to complete the proof we need to show that $(a_-, b_-) \in \mathbf{H}_0^*$. We also need to show that $a_+$ is outer.

Clearly we have $a_-a_-^* = 1 + b_-b_-^*$ almost everywhere on $\mathbf{T}$ since the other matrices in the equation are $SU(1,1)$ almost everywhere.

Next we check that $a_-$ and $b_-$ have the correct holomorphicity properties. From the fixed point equations,

$$a_- = aa_+^* - bb_+^*$$

$$= a\frac{1}{a_+(\infty)} + aP_D(\frac{b}{a}b_+{}^*) - bb_+^*$$

$$= a\frac{1}{a_+(\infty)} - a(\mathrm{id} - P_D)(\frac{b}{a}b_+{}^*)$$

Clearly this is an element of $H^2(D^*)$. Moreover, the constant term obeys

$$a_-(\infty) = \frac{a(\infty)}{a_+(\infty)}$$

and so is positive as required.

Similarly, using the fixed point equation for $b_+$,

$$b_- = -ab_+ + ba_+$$

$$= -aP_D(\frac{b}{a}a_+) + ba_+$$

$$= (1 - P_D)(\frac{b}{a}a_+)$$

Thus $b_-$ is in $H^2(D^*)$

To prove that we have indeed solved the Riemann-Hilbert problem, we have to verify that $a_-$ and $a_+$ are outer.

Consider the equation

$$a = a_-a_+ + b_-b_+^*$$

Every function in this equation is holomorphic in $D^*$. We divide by the outer function $a$ to obtain

$$1 = a_-a_+a^{-1} + b_-b_+{}^*a^{-1}$$

since the first summand on the right is larger in modulus on $\mathbf{T}$ than the second, we conclude that

$$\mathrm{Re}\,(a_-a_+a^{-1}) \geq 1/2$$

almost everywhere on $\mathbf{T}$. This implies that the function $a_-a_+a^{-1}$, which is in $H^1(D^*)$ and thus equal to its Poisson integral on $D^*$, has real part larger than $1/2$ on $D^*$. Then the reciprocal function $\frac{a}{a_-a_+}$ is in $H^\infty(D^*)$ and

$$\frac{1}{a_+} = \left(\frac{a}{a_-a_+}\right)\left(\frac{a_-}{a_+}\right)$$

is in $H^2(D^*)$ by Lemma 19. By Lemma 32, $a_+$ is outer. Likewise one concludes that $a_-$ is outer.

This completes the proof of Lemma 18. □

## 3. Existence of an inverse NLFT for unbounded $a$

In this section we prove that for every $(a, b) \in \mathbf{L}$, there exists a factorization

$$(a_-, b_-)(a_+, b_+) = (a, b)$$

with $(a_-, b_-) \in \mathbf{H}_0^*$ and $(a_+, b_+) \in \mathbf{H}$. We shall call such a factorization a Riemann-Hilbert factorization. This factorization is not necessarily unique. If it is unique, then we say that $(a, b)$ has a unique Riemann-Hilbert factorization.

In the previous section we used the Banach fixed point theorem to produce a Riemann-Hilbert factorization when $a \in H^\infty(D^*)$. This same approach does not work in the general case; we shall instead use the Riesz representation theorem for linear functionals on a Hilbert space. In general there will be several choices of Hilbert space to work with, which will cause non-uniqueness of the Riemann-Hilbert factorization.

We introduce two examples of such Hilbert spaces, which in general may be different and will turn out to be extremal examples. They are vector spaces over the real (not complex) numbers. Given $(a, b) \in \mathbf{L}$, we consider the following inner product on pairs of measurable functions on $\mathbf{T}$:

$$\langle (A', B'), (A, B) \rangle := \int_{\mathbf{T}} \mathrm{Re}\,[A'(A^* - \frac{b}{a}B^*) + (B')^*(B - \frac{b}{a}A)] \tag{29}$$

whenever the integral on the right-hand side is absolutely integrable. We emphasize that absolute integrability is only required for the real part of the algebraic expression in the integrand.

This inner product is positive definite, as we see from the following calculation:

$$\begin{aligned} \|(A, B)\| := \langle (A, B), (A, B) \rangle &\geq \int_{\mathbf{T}} |A|^2 + |B|^2 - 2|b/a||A||B| \\ &\geq \int_{\mathbf{T}} |b/a|(|A| - |B|)^2 + \int (1 - |b/a|)(|A|^2 + |B|^2) \\ &\geq \frac{1}{2} \int_{\mathbf{T}} (1 - |b/a|^2)(|A|^2 + |B|^2) \\ &= \frac{1}{2} \int_{\mathbf{T}} (|A|^2 + |B|^2)|a|^{-2} \geq 0 \end{aligned}$$

Equality holds in the last estimate if and only if $A$ and $B$ vanish almost everywhere on $\mathbf{T}$. Therefore the inner product is positive definite. Moreover, we have seen that the integrand in (29) is nonnegative almost everywhere if $(A, B) = (A', B')$.

The above calculation also shows that a necessary condition for for the inner product (29) to be defined is $A/a \in L^2(\mathbf{T})$ and $B/a \in L^2(\mathbf{T})$.

Define $H_{\max}$ to be the space of all pairs $(A, B)$ such that $A/a \in H^2(D^*)$ and $B/a^* \in H^2(D)$ and $\|(A, B)\| < \infty$ with respect to the inner product (29). This space is evidently a pre-Hilbert space with inner product (29). It is indeed a Hilbert space, because for any Cauchy sequence, the boundary values of $A/a$ and $B/a^*$ converge in $L^2(\mathbf{T})$ and thus remain in $H^2(D^*)$ and $H^2(D)$ respectively. By an application of Fatou's lemma, the limit has again finite norm and thus is in $H_{\max}$. The previous display shows that $H_{\max}$ is continuously embedded in $aH^2(D^*) \times a^* H^2(D)$.

The space $H^2(D^*) \times H^2(D)$ is contained in $H_{\max}$, because

$$\langle (A, B), (A, B) \rangle \leq 2 \int_{\mathbf{T}} |A|^2 + |B|^2$$

Define $H_{\min}$ to be the closure of $H^2(D^*) \times H^2(D)$ in $H_{\max}$. As we will see, there are examples of data $(a, b)$ for which the space $H_{\min}$ is strictly contained in $H_{\max}$.

We now introduce the real linear functional to which the Riesz representation theorem will be applied. It takes the same form on $H_{\max}$ and $H_{\min}$ and is given by

$$\lambda : (A, B) \to \mathrm{Re}\,[A(\infty)]$$

Observe that this linear functional is indeed continuous on $H_{\max}$ and thus also $H_{\min}$, since it is even continuous on the larger space $aH^2(D^*) \times a^*H^2(D)$, as can be seen immediately from

$$\mathrm{Re}\,A(\infty) = a(\infty)\mathrm{Re} \int_{\mathbf{T}} A/a$$

Let $(A_{\min}, B_{\min})$ be the unique element in $H_{\min}$ which produces this linear functional in $H_{\min}$ and is guaranteed to exist by the Riesz representation theorem:

$$\langle (A_{\min}, B_{\min}), (A, B)\rangle = \lambda(A, B)$$

for all $(A, B) \in H_{\min}$. Let $(A_{\max}, B_{\max})$ be the unique element in $H_{\max}$ which produces this linear functional in $H_{\max}$.

**Theorem 7.** *Let $(a, b) \in \mathbf{L}$. Then there exists a factorization*

$$(a, b) = (a_-, b_-)(a_+, b_+)$$

*with $(a_-, b_-) \in \mathbf{H}_0^*$ and $(a_+, b_+) \in \mathbf{H}$. Moreover, with $(A_{\max}, B_{\max})$ and $(A_{\min}, B_{\min})$ defined as above, two possible choices of such a Riemann-Hilbert factorization are given by*

$$(a_+, b_+) := (A_{\min}, B_{\min})A_{\min}(\infty)^{-1/2} \tag{30}$$

*and*

$$(a_+, b_+) := (A_{\max}, B_{\max})A_{\max}(\infty)^{-1/2} \tag{31}$$

*with the corresponding $(a_-, b_-)$, which are easily determined by matrix division.*

**Proof.** The proofs that (30) and (31) give Riemann Hilbert factorizations are very similar. We shall formulate the proof for (30) and comment on the changes needed to prove (31).

Define $L$ to be the space of all pairs $(A, B)$ of measurable functions such that $A/a \in H^2(D^*)$, $B/a^* \in L^2(\mathbf{T})$, and $\|(A, B)\| < \infty$.

We claim that $H^2(D^*) \times L^2(\mathbf{T})$ is dense in $L$. Indeed, let $(A, B) \in L$ be orthogonal to all elements of $H^2(D^*) \times L^2(\mathbf{T})$. Choosing $A' = 0$ and $B'$ of modulus one such that $(B')^*(B - \frac{b}{a}A)$ is nonnegative real, we have

$$0 = \langle (A', B'), (A, B)\rangle = \int |B - \frac{b}{a}A|$$

and thus

$$B - \frac{b}{a}A = 0 \tag{32}$$

almost everywhere. Now choosing $A' \in H^2$ arbitrary and $B' = 0$ we obtain

$$0 = \langle (A', B'), (A, B)\rangle$$

$$= \int_{\mathbf{T}} A'(A^* - \frac{b}{a}B^*) = \int_{\mathbf{T}} A'A^*(1 - |b/a|^2) = \int_{\mathbf{T}} \frac{A'}{a}\frac{A^*}{a^*}$$

By Beurling's theorem (see [**9**]), since $a$ is outer, the set of all $A'/a$ with $A' \in H^2(D^*)$ is dense in $H^2(D^*)$. Thus, redefining $A'$,

$$\int_{\mathbf{T}} A' \frac{A^*}{a^*} = 0$$

for all $A' \in H^2(D^*)$. But $A^*/a^* \in H^2(D)$, thus $A^*/a^* = 0$. Hence $A = B = 0$ by (32) and we have shown that the orthogonal complement of $H^2(D^*) \times L^2(\mathbf{T})$ is trivial, thus proving our claim.

Define $H_n$ to be the closure of $H^2(D^*) \times z^n H^2(D)$ in $L$, in particular $H_0 = H_{\min}$. (Here is the main difference in proving the theorem for $(A_{\min}, B_{\min})$ and $(A_{\max}, B_{\max})$. To prove the theorem for $(A_{\max}, B_{\max})$, one would need to define $H_n$ to be the space of all $(A, B) \in L$ such that $B/a^* \in z^n H^2(D)$.)

Since evaluation of $z^{-n}B$ at 0 is a continuous functional on $H_n$,

$$z^{-n}B = a^*(0) \int_{\mathbf{T}} z^{-n}B/a^*$$

we see that $H_{n+1}$ is precisely the subspace of $H_n$ of all $(A, B)$ such that $z^{-n}B$ vanishes at 0. Thus $H_{n+1}$ has real co-dimension two in $H_n$.

Let $H_\infty$ be the intersection of all $H_n$ for $n \in \mathbf{Z}$, then it is clear that $H_\infty$ consists of pairs $(A, B)$ such that $B$ vanishes to infinite order at 0 and thus is identically equal to 0. Finiteness of the norm of $(A, B)$ is then equivalent to $A \in H^2(D^*)$ and thus evidently $H_\infty = H^2(D^*) \times \{0\}$.

Let $H_{-\infty}$ be the closure of the union of all $H_n$ for $n \in \mathbf{Z}$. Since every element of $H^2(D^*) \times L^2(\mathbf{T})$ can be approximated by a sequence of elements in spaces $H_n$ with decreasing $n$, we see that $H_{-\infty}$ is equal to the closure of $H^2(D^*) \times L^2(\mathbf{T})$ which is all of $L$.

Let $(A_n, B_n)$ be the element which represents the linear functional $\lambda$ in the subspace $H_n$. It is easy to see that

$$(A_\infty, B_\infty) = (1, 0)$$

$$(A_{-\infty}, B_{-\infty}) = a(\infty)(a, b)$$

Observe that the operation $(A, B) \to (B^* z^n, A^* z^n)$ is a bijection on the space $H^2(D^*) \times z^n H^2(D)$, and extends to a bijective isometry of $H_n$.

We claim

$$(A_{n+1}, B_{n+1}) = (A_n, B_n) - F_n(B_n^* z^n, A_n^* z^n) \tag{33}$$

for a certain complex number $F_n \in D$.

Indeed, since $\lambda$ is non-zero (it is so on $H_\infty$), we have

$$\operatorname{Re} A_n(\infty) = \langle (A_n, B_n), (A_n, B_n) \rangle > 0 \tag{34}$$

and there is a unique $F_n \in \mathbf{C}$ such that the off-diagonal entry on the right-hand side of (33) vanishes to order $n+1$ at 0. For later reference we pause to argue that taking real part on the left-hand side of (34) is superfluous since $A_n(\infty)$ itself is positive. Namely, $(A, B) \to (cA, cB)$ is an isometry of $H_n$ for $|c| = 1$ and since

$$\operatorname{Re}[cA_n] = \langle (A_n, B_n), (cA_n, cB_n) \rangle$$

is maximized for $c = 1$, we have that $A_n(\infty) > 0$.

Now we observe that

$$(B_n^* z^n, A_n^* z^n)$$

is orthogonal to $H_{n+1}$. Namely, let $(A,B) \in H_{n+1}$, then

$$\langle (B_n^* z^n, A_n^* z^n), (A,B)\rangle = \langle (A_n, B_n), (B^* z^n, A^* z^n)\rangle = \lambda(B^* z^n, A^* z^n) = 0$$

Thus the right-hand side of (33) is the orthogonal projection of $(A_n, B_n)$ onto $H_{n+1}$ and thus indeed equal to $(A_{n+1}, B_{n+1})$.

We now verify that $|F_n| < 1$. This simply follows from the fact that

$$\|(B_n^* z^n, A_n^* z^n)\| = \|(A_n, B_n)\|$$

and the fact that the terms in (33) form a Pythagorean triple and $(A_{n+1}, B_{n+1}) \neq 0$. Another consequence is that

$$\|(A_{n+1}, B_{n+1})\| = (1 - |F_n|^2)^{1/2} \|(A_n, B_n)\| \tag{35}$$

Each vector $(A_n, B_n)$ is the orthogonal projection of $(A_{-\infty}, B_{-\infty})$ onto $H_n$, and the projection of $(A_n, B_n)$ onto $H_\infty$ is $(A_\infty, B_\infty)$. Thus the length of each vector $(A_n, B_n)$ is squeezed between two finite numbers

$$\|(A_{-\infty}, B_{-\infty})\| \geq \|(A_n, B_n)\| \geq \|(A_\infty, B_\infty)\|$$

By an inductive argument using (35) we see that

$$\prod_n (1 - |F_n|^2)^{1/2}$$

is a convergent product and thus the sequence $F = (F_n)$ is in $l^2(\mathbf{Z}, D)$.

Let $(\tilde{a}, \tilde{b})$ be the nonlinear Fourier transform of $F$. We need to show that $(a,b) = (\tilde{a}, \tilde{b})$. We claim that

$$(A_n, B_n) \prod_{k \geq n} (1 - |F_n|^2)^{1/2}$$

is the nonlinear Fourier transform $(\tilde{a}^{(\geq n)}, \tilde{b}^{(\geq n)})$ of the truncated sequence $F^{(\geq n)}$.

Consider

$$((\tilde{a}^{(\geq n)})^*, -\tilde{b}^{(\geq n)})(A_n,, B_n) \prod_{k \geq n} (1 - |F_k|^2)^{1/2} \tag{36}$$

where we have used the convention to read each vector as the first row of a positive scalar multiple of a $SU(1,1)$ matrix.

Observe that (33) reads as

$$(A_{n+1}, B_{n+1}) = (1, -F_n z^n)(A_n, B_n)$$

Therefore, by the recursion equation for $(a^{(\geq n)}, b^{(\geq n)})$, the quantity (36) is independent of the parameter $n$.

There is an increasing sequence $n_k$ of integers such that pointwise almost everywhere on $\mathbf{T}$ we have

$$(\tilde{a}^{(\geq n_k)}, \tilde{b}^{(\geq n_k)}) \to (1,0)$$

$$(A_{n_k}, B_{n_k}) \to (A_\infty, B_\infty) = (1,0)$$

as $k \to \infty$. For the first limit this follows from convergence of the sequence $F^{(\geq n)}$ to 0 in $l^2$ and thus convergence of $\log|\tilde{a}^{(\geq n)}|$ in $L^1$, convergence of the phase $\tilde{a}^{(\geq n)}/|\tilde{a}^{(\geq n)}|$ in $L^2$, and convergence of $\tilde{b}^{(\geq n)}/\tilde{a}^{(\geq n)}$ in $L^2$. For the second limit this follows from convergence of $A_n/a$ and $B_n/a^*$ in $L^2$. Thus (36) is equal to $(1,0)$ almost everywhere. Taking a similar limit as $n \to -\infty$ we observe

$$(\tilde{a}^*, -\tilde{b})(A_{-\infty}, B_{-\infty}) \prod_k (1 - |F_k|^2)^{1/2} = (1,0)$$

This proves that $(a, b)$ is a positive scalar multiple of $(\tilde{a}, \tilde{b})$, but since both are $SU(1,1)$ valued they are indeed equal.

Finally, we observe from splitting the sequence $F$ as $F^{(<0)}+F^{(\geq 0)}$ that $(\tilde{a}^{(\geq 0)}, \tilde{b}^{(\geq 0)})$ is the right factor of a Riemann-Hilbert factorization of $(a, b)$, and that we have

$$(A_{\min}, B_{\min}) = \tilde{a}_0(\infty)(\tilde{a}_0, \tilde{b}_0)$$

This completes the proof that (30) produces a Riemann-Hilbert factorization.

The proof for $(A_{\max}, B\max)$ is similar with changes as indicated above. □

The above construction of a Riemann-Hilbert factorization easily provides the following strengthening. Let $(a, b) \in \mathbf{L}$ and

$$(a, b) = (a_-, b_-)(a_+, b_+)$$

be any Riemann-Hilbert factorization. Then the vector $(A_{\min}, B_{\min})$ constructed as above with respect to $(a, b)$ is identical to the vector $(A_{+,\min}, B_{+,\min})$ constructed as above with respect to $(a_+, b_+)$. To see this, it suffices to show that the inner products

$$\langle (A', B'), (A, B)\rangle := \int_{\mathbf{T}} \mathrm{Re}\,[A'(A^* - \frac{b}{a}B^*) + (B')^*(B - \frac{b}{a}A)]$$

$$\langle (A', B'), (A, B)\rangle_+ := \int_{\mathbf{T}} \mathrm{Re}\,[A'(A^* - \frac{b_+}{a_+}B^*) + (B')^*(B - \frac{b_+}{a_+}A)]$$

coincide on the space $H^2(D^*) \times H^2(D)$. Indeed, by polarization it suffices to show this for $(A, B) = (A', B')$. However, the difference of these inner products is then given by

$$\mathrm{Re} \int_{\mathbf{T}} 2[\frac{b}{a} - \frac{b_+}{a_+}]AB^*$$

We have

$$b_- = -ab_+ + ba_+$$

$$\frac{b_-}{aa_+} = \frac{b}{a} - \frac{b_+}{a_+}$$

Observe that here the right-hand side is bounded on $\mathbf{T}$, while the left-hand side is in $H_0^2(D^*)$ and thus in $H_0^\infty(D^*)$. The difference of the inner products is then given by

$$\mathrm{Re} \int_{\mathbf{T}} 2\frac{b_-}{aa_+}AB^*$$

Since $A$ and $B^*$ are in $H^2(D^*)$, this difference is equal to 0.

**Theorem 8.** *Let $(a, b) \in \mathbf{L}$. Then there is a unique factorization*

$$(a, b) = (a_{--}, b_{--})(a_{\mathbf{o}}, b_{\mathbf{o}})(a_{++}, b_{++})$$

*such that*

$$(a_{--}, b_{--}) \in \mathbf{H}_0^*$$

$$(a_{\mathbf{o}}, b_{\mathbf{o}}) \in \mathbf{H}_0^* \cap \mathbf{H}$$

$$(a_{++}, b_{++}) \in \mathbf{H}$$

*and $(a_{--}, b_{--})$ and $(a_{++}, b_{++})$ do not have any Riemann-Hilbert factorizations other than*

$$(a_{++}, b_{++}) = (1, 0)(a_{++}, b_{++})$$

*and*

$$(a_{--}, b_{--}) = (a_{--}, b_{--})(1, 0)$$

*Moreover, we have the sub-factorization property: Any Riemann-Hilbert factorization*

$$(a,b) = (a_-, b_-)(a_+, b_+) \tag{37}$$

*comes with further (obviously unique) Riemann-Hilbert factorizations*

$$(a_-, b_-) = (a_{--}, b_{--})(a_{-\mathbf{o}}b_{-\mathbf{o}})$$

$$(a_+, b_+) = (a_{\mathbf{o}+}, b_{\mathbf{o}+})(a_{++}b_{++})$$

*where* $(a_{--}, b_{--})$ *and* $(a_{++}, b_{++})$ *are as above and*

$$(a_{-\mathbf{o}}b_{-\mathbf{o}}), (a_{\mathbf{o}+}, b_{\mathbf{o}+}) \in \mathbf{H}_0^* \cap \mathbf{H}$$

*We have*

$$(a_{\mathbf{o}}, b_{\mathbf{o}}) = (a_{-\mathbf{o}}b_{-\mathbf{o}})(a_{\mathbf{o}+}, b_{\mathbf{o}+}) \tag{38}$$

*The passage from the Riemann-Hilbert factorization (37) of* $(a,b)$ *to the Riemann-Hilbert factorization (38) of* $(a_{\mathbf{o}}, b_{\mathbf{o}})$ *constitutes a bijective correspondence between the Riemann-Hilbert factorizations of* $(a,b)$ *and the Riemann-Hilbert factorizations of* $(a_{\mathbf{o}}, b_{\mathbf{o}})$.

Before proving the theorem, we prove the following lemma.

**Lemma 20.** *If* $(a,b) \in \mathbf{H}$ *and*

$$(a,b) = (a_-, b_-)(a_+, b_+)$$

*is a Riemann-Hilbert factorization, then* $(a_-, b_-) \in \mathbf{H}$. *Conversely, if*

$$(a_-, b_-) \in \mathbf{H}_0^* \cap \mathbf{H}$$

$$(a_+, b_+) \in \mathbf{H}$$

*then*

$$(a,b) := (a_-, b_-)(a_+, b_+) \in \mathbf{H}$$

Observe that by reflection there is an analogous lemma with $(a,b) \in \mathbf{H}_0^*$ and $(a_+, b_+) \in \mathbf{H}_0^*$.

**Proof.** Assume we have a Riemann-Hilbert factorization of $(a,b) \in \mathbf{H}$. Then

$$b_+ = a_-^* b - b_- a^*$$

$$\frac{b_-}{a_-^*} = \frac{b}{a^*} - \frac{b_+}{a_-^* a^*} \tag{39}$$

Every summand on the right-hand side is in $H^2(D)$, hence so is the expression on the left-hand side. This proves the first statement of the lemma.

Next, assume

$$(a_-, b_-) \in \mathbf{H}_0^* \cap \mathbf{H}$$

$$(a_+, b_+) \in \mathbf{H}$$

Clearly the product $(a,b)$ is in $\mathbf{L}$ by Lemma 16. Then $b/a^* \in \mathbf{H}(D)$ follows again from (39). This proves the second statement of the lemma. □

Now we can prove Theorem 8.

**Proof.** With the notation of Theorem 7, set

$$(a_{++}, b_{++}) = A_{\min}(\infty)^{-1/2}(A_{\min}, B_{\min})$$

As we have observed in the discussion prior to the statement of Theorem 8, for any Riemann-Hilbert factorization

$$(a, b) = (a_-, b_-)(a_+, b_+)$$

we have a Riemann-Hilbert factorization

$$(a_+, b_+) = (a_{\mathbf{o}+}, b_{\mathbf{o}+})(a_{++}, b_{++})$$

The lemma just shown implies that

$$(a_{\mathbf{o}+}, b_{\mathbf{o}+}) \in \mathbf{H}_0^* \cap \mathbf{H}$$

Thus we have shown the sub-factorization property for $(a_{++}, b_{++})$.

By the sub-factorization property, $(a_{++}, b_{++})$ is the only possible right factor in a Riemann-Hilbert factorization of $(a, b)$ which does not have a Riemann-Hilbert factorization other than identity times itself. We claim that $(a_{++}, b_{++})$ indeed does not have any further Riemann-Hilbert factorization. Assume we have a Riemann-Hilbert factorization

$$(a_{++}, b_{++}) = (\tilde{a}, \tilde{b})(a_{+++}, b_{+++})$$

Then by an application of Lemma 20 we observe that $(a_{+++}, b_{+++})$ is a right factor of a Riemann-Hilbert factorization of $(a, b)$, and thus by the sub-factorization property

$$(a_{+++}, b_{+++}) = (\tilde{a}^*, -\tilde{b})(a_{++}, b_{++})$$

is also a Riemann-Hilbert factorization. Thus both $(\tilde{a}, \tilde{b})$ and its inverse are in $\mathbf{H} \cap \mathbf{H}_0^*$. Thus $\tilde{b}/\tilde{a}$ is in $H^2(D) \cap H_0^2(D^*)$. Therefore $\tilde{b}/\tilde{a} = 0$ and consequently $(\tilde{a}, \tilde{b}) = (1, 0)$. This proves that $(a_{++}, b_{++})$ does not have any nontrivial Riemann-Hilbert factorization.

Since by Lemma 20 any triple factorization of $(a, b)$ as in the Theorem gives a Riemann-Hilbert factorization

$$(a, b) = [(a_{--}, b_{--})(a_{\mathbf{o}}, b_{\mathbf{o}})][(a_{++}, b_{++})]$$

we have uniquely identified the factor $(a_{++})$ as the only possible in such a triple factorization. Similarly we can uniquely identify the factor $(a_{--}, b_{--})$, and by an application of the sub-factorization property we actually obtain the triple factorization from knowledge of these two factors.

Finally, we observe that any Riemann-Hilbert factorization of $(a, b)$ gives a Riemann-Hilbert factorization of $(a_{\mathbf{o}}, b_{\mathbf{o}})$ as described in the theorem, and vice versa every Riemann-Hilbert factorization of $(a_{\mathbf{o}}, b_{\mathbf{o}})$ necessarily has two factors in $\mathbf{H}_0^* \cap \mathbf{H}$ by Lemma 20 and therefore comes from a Riemann-Hilbert factorization of $(a, b)$. This proves the theorem. □

For any $(a, b) \in \mathbf{L}$ call $\log|a(\infty)|$ the energy of $(a, b)$. The energy of $(a, b)$ is a nonnegative real number. Indeed, it is positive unless $(a, b) = (1, 0)$. If

$$(a, b) = (a_-, b_-)(a_+, b_+)$$

is a Riemann-Hilbert factorization, then we have additivity of the energies

$$\log|a(\infty)| = \log|a_+(\infty)| + \log|a_-(\infty)|$$

(evaluate $a = a_-a_+ + b_-b_+^*$ at $\infty$).

The sub-factorization property shows that the right factor $(a_{++}, b_{++})$ of the triple factorization minimizes the energy among all right factors of Riemann-Hilbert factorizations of $(a, b)$, and indeed is a unique minimizer. Since $(a_{++}, b_{++})$ was constructed through a minimal Hilbert space $H_{\min}$ in Theorem 7, it is natural to guess that the solution constructed from the space $H_{\max}$ in Theorem 7 maximizes the energy. This is the main content of the following lemma:

**Lemma 21.** *With the notation of Theorems 7 and 8 we have*

$$(a_{\mathbf{o}}, b_{\mathbf{o}})(a_{++}, b_{++}) = (A_{\max}, B_{\max})A_{\max}(\infty)^{-1/2}$$

**Proof.** Consider the notation of Theorem 7 and the space $H_{\max}$. We claim that for every Riemann-Hilbert factorization

$$(a, b) = (a_-, b_-)(a_+, b_+)$$

we have

$$(a_+, b_+) \in H_{\max}$$

By Lemma 19, $a_+/a$ and $b_+^*/a$ are in $H^2(D^*)$. Therefore it suffices to show that $\|(a_+, b_+)\|$ is finite. However, we have

$$\begin{aligned} \|(a_+, b_+)\| &= \int \operatorname{Re}\,[a_+(a_+^* - \frac{b}{a}b_+^*) + b_+^*(b_+ - \frac{b}{a}a_+^*)] \\ &= \int \operatorname{Re}\,[\frac{a_+a_-}{a} - \frac{b_+^*b_-}{a}] \end{aligned}$$

and the right-hand side has been shown to be finite in the proof of Lemma 19. Therefore, $(a_+, b_+) \in H_{\max}$.

We have the quantitative estimate

$$\begin{aligned} \|(a_+, b_+)\| &= -1 + 2\int \operatorname{Re}\,[\frac{a_+a_-}{a}] \\ &\leq -1 + 2a_+(\infty)a_-(\infty)a(\infty)^{-1} = 1 \end{aligned}$$

In the inequality we have used the observation in the proof of Lemma 19 that $a_-a_+/a$ has positive real part, and thus the real part of its value at $\infty$ is the total mass of the positive measure given by the Herglotz representation theorem. The total mass of this measure dominates the total mass of its absolutely continuous part $\operatorname{Re}(a_+a_-/a)$.

Applying the linear functional $\lambda$ to $(a_+, b_+) \in H_{\max}$ gives

$$\begin{aligned} a_+(\infty) &= \langle (a_+, b_+), (A_{\max}, B_{\max}) \rangle \\ &\leq \|(a_+, b_+)\|\|A_{\max}, B_{\max})\| \\ &\leq A_{\max}(\infty)^{1/2} \end{aligned}$$

This proves that $(a_+, b_+)$ has smaller energy than the solution $(A_{\max}, B_{\max})A_{\max}(\infty)^{-1/2}$ of the Riemann-Hilbert factorization theorem. Thus $(A_{\max}, B_{\max})A_{\max}(\infty)^{-1/2}$ maximizes the energy of the right factor of a Riemann-Hilbert factorization of $(a, b)$. On the other hand, by symmetry $(a_{--}, b_{--})$ (uniquely) minimizes the energy of the left factor of a Riemann-Hilbert factorization, just as $(a_{++}, b_{++})$ minimizes the energy of the right factor. By additivity of the energy this proves the lemma. □

# LECTURE 4
## Rational functions as Fourier transform data

### 1. The Riemann-Hilbert problem for rational functions

The class $\mathbf{L}$ of nonlinear Fourier transform data of $l^2$ sequences contains elements $(a, b)$ such that $a$ and $b$ are rational functions in $z$. We call $(a, b)$ rational if $a$ and $b$ are rational.

Indeed, any pair $(a, b)$ of rational functions is an element of $\mathbf{L}$ if $aa^* = 1 + bb^*$, $a(\infty) > 0$, and the function $a$ has no zeros and poles in $D^*$. These pairs can easily be parameterized by the function $b$, as the following lemma states:

**Lemma 22.** *For each rational function $b$ there is precisely one rational function $a$ such that $aa^* = 1 + bb^*$, $a$ has no zeros and poles in $D^*$, and $a(\infty) > 0$. This is the unique function $a$ such that $(a, b) \in \mathbf{L}$.*

*For rational $(a, b) \in \mathbf{L}$, we have $(a, b) \in \mathbf{H}$ if and only if $b$ has no poles in $D$, and $(a, b) \in \mathbf{H}_0$ if and only if in addition $b(0) = 0$. Likewise, we have $(a, b) \in \mathbf{H}^*$ if and only if $b$ has no poles in $D^*$ and we have $(a, b) \in \mathbf{H}_0^*$ if and only if in addition $b(\infty) = 0$.*

**Proof.** Let $b$ be a rational function. Consider the rational function $g = 1 + bb^*$. Then $g(z) = g^*(z)$ and the zeros and poles of $g$ are symmetric about $\mathbf{T}$: if $z$ is a pole of order $n$, then so is $z^*$ and likewise for the zeros. Moreover, there are no zeros of $g$ on $\mathbf{T}$ and the poles of $g$ on $\mathbf{T}$ are of even order.

Let $a$ be a rational function whose zeros and poles in $D$ are precisely the zeros and poles of $g$ in $D$ with the same order, whose poles on $\mathbf{T}$ are precisely the poles of $g$ but with half the order, and which has no zeros on $\mathbf{T} \cup D^*$ and no poles on $D^*$. Thus poles and zeros of $a$ are completely specified and $a$ is determined up to a scalar factor. We assume $a$ to be positive at $\infty$, which determines the phase of this scalar factor.

Consider

$$f = (aa^*)^{-1}(1 + bb^*)$$

Then this rational function evidently has no zeros and no poles and therefore it is constant. Since it is positive on $\mathbf{T}$, we may normalize $a$ with a positive factor such that $f = 1$.

We claim that $(a, b) \in \mathbf{L}$. Certainly $aa^* = 1 + bb^*$ by construction. The function $a$ is holomorphic in $D^*$ with no zeros in $D^*$. Any rational function with these properties is outer on $D^*$. To see this, it suffices by multiplicativity of outer functions to show that functions of the form $(1/z - 1/z_0)$ with $z_0 \in D \cup \mathbf{T}$ are outer, which is easy to verify.

This proves that there exists an $a$ with $(a, b) \in \mathbf{L}$. Uniqueness follows very generally from the fact that the normalized outer function $a$ is determined by $|a|$ almost everywhere on $\mathbf{T}$, and the latter is determined by $b$. This proves the first statement of the lemma.

Clearly holomorphicity of $b$ in $D$ is necessary for $(a, b) \in \mathbf{H}(D)$. However, if $b$ is holomorphic in $D$, which is the same as saying $b$ has no poles in $D$, then $b/a$ is a rational function holomorphic in $D$ and bounded by 1 almost everywhere on $\mathbf{T}$, and thus holomoprhic in a neighborhood of $D \cup \mathbf{T}$ and thus in $H^2(D)$. This proves $(a, b) \in \mathbf{H}(D)$. The statement about $\mathbf{H}_0(D)$ is clear. This together with the symmetric statement for $D^*$ proves the remaining statements of the lemma. □

The next lemma states that solving the problem of Riemann-Hilbert factorization does not leave the class of rational functions.

**Lemma 23.** *Assume $(a, b) \in \mathbf{L}$ is rational. Given any factorization*

$$(a_-, b_-)(a_+ b_+) = (a, b)$$

*with $(a_-, b_-) \in \mathbf{H}^2(D^*)$ and $(a, b) \in \mathbf{H}(D)$, then $(a_-, b_-)$ and $(a_+, b_+)$ are also rational.*

**Proof.** Recall from Lemma 19 that

$$\frac{a_-}{a}, \frac{b_-}{a}, \frac{a_+}{a}, \frac{b_+^*}{a} \in H^2(D^*)$$

In particular

$$\int_{\mathbf{T}} \left|\frac{a_+}{a}\right|^2 (r\cdot) \leq C$$

for $r \geq 1$. Now $aa^* = 1 + bb^*$ implies that the poles of $b$ on $\mathbf{T}$ have the same order as the poles of $a$ on $\mathbf{T}$. Thus $b/a$, which is also a rational function, is actually holomorphic on a neighborhood of $\mathbf{T}$.

Using

$$a_-^* = a_+ a^* - b^* b_+$$

we obtain

$$\frac{a_+}{a^*} = \frac{1}{a^*}\frac{a_-^*}{a^*} + \frac{b}{a^*}\frac{b_+}{a^*}$$

On the right-hand side, the functions $a_-^*/a^*$ and $b_+/a^*$ are in the Hardy space $H^2(D)$, while the rational functions $1/a^*$ and $b^*/a^*$ are holomorphic in a neighborhood of $\mathbf{T}$.

Therefore, $a_+$ has a meromorphic extension to $D$ which is holomorphic in an annulus $1 - \epsilon < |z| < 1$ for some small $\epsilon$ and satisfies

$$\int_{\mathbf{T}} \left|\frac{a_+}{a^*}\right|^2 (r\cdot) \leq C$$

for $1 - \epsilon < r \leq 1$

Observe that $a$ and $a^*$ have comparable moduli in a small neighborhood of $\mathbf{T}$ since the quotients $a/a^*$ and $a^*/a$ are holomorphic near $\mathbf{T}$.

Thus

$$\int_{\mathbf{T}} \left|\frac{a_+}{a}\right|^2 (r\cdot) \leq C$$

for $1-\epsilon < r < 1$ for some small $\epsilon$, and the same estimate has been observed previously for $r \geq 1$.

We claim that holomorphicity of $a_+/a$ in a neighborhood of $\mathbf{T}$ with possible exception on $\mathbf{T}$ together with the above estimates implies that $a_+/a$ is indeed holomorphic across $\mathbf{T}$.

In the current situation that $a_+/a$ is in addition meromorphic in $D$ and $D^*$ with finitely many poles we can argue as follows.

We may remove the poles of $a_+/a$ in $D$ by the following recursive procedure. If $a_+/a$ has a pole at $z_\infty \in D$, then we subtract a constant from $a_+/a$ so that the new function has a zero at a distinct point $z_0 \in D$, and then we multiply the function by $(z-z_\infty)/(z-z_0)$. This reduces the order of the pole at $z_\infty$ and leaves the order of all other poles unchanged. Iterating this procedure we obtain a function $g$ which is holomorphic in $D$ and $D^*$. The above $L^2$ estimates prevail throughout this iteration, possibly with different constants $C$, so $g$ is in $H^2(D) \cap H^2(D^*)$. Therefore $g$ is constant, and we conclude that $a_+/a$ is rational. The estimates near $\mathbf{T}$ then imply that it has no poles on $\mathbf{T}$.

More generally, the claim can be proved using the theorem of Morera: a function is holomorphic in a disc if the Cauchy integral over each triangle vanishes. For triangles which avoid $\mathbf{T}$ this is obvious for $a_+/a$, and for triangles which intersect $\mathbf{T}$ one obtains vanishing of the Cauchy integral by approximating the triangle by shapes avoiding $\mathbf{T}$ and then using maximal function estimates to pass to the limit.

This proves that $a_+$ is rational, and one can argue similarly that $b_+, a_-, b_-$ are rational.

This proves Lemma 23. □

The lemma just proved reduces the Riemann-Hilbert problem for rational $(a,b)$ to a purely algebraic problem in the class of rational functions. Even better, the following lemma states that the solution functions $a_-, b_-, a_+, b_+$ are in a sense subordinate to $a, b$. This reduces the Riemann-Hilbert problem to a finite dimensional problem.

If $f$ is meromorphic near $z \in \mathbf{C}$, denote by ord(f, z) the order of the pole of $f$ at $z$. Thus

$$f(\zeta)(z-\zeta)^{\mathrm{ord(f,z)}}$$

is holomorphic at $z$ and does not vanish at $z$. For rational functions we define the order at $\infty$ in the usual manner using a change of coordinates on the Riemann sphere. Observe that the order is a negative number if $f$ vanishes at $z$.

Call a rational function $g$ subordinate to another rational function $f$ on a certain domain if for all points $z$ in the domain such that $\mathrm{ord(g,z)} > 0$ we have $\mathrm{ord(f,z)} \geq \mathrm{ord(g,z)}$. We call $g$ subordinate to $f$ if $g$ is subordinate to $f$ on the whole Riemann sphere. Clearly, if we fix $f$, the set of rational functions $g$ subordinate to $f$ is a finite dimensional vector space.

**Lemma 24.** *Let $(a,b) \in \mathbf{L}$ be rational. Then $a$ is subordinate to $bb^*$.*

*If $(a_-, b_-) \in \mathbf{H}_0$ and $(a_+, b_+) \in \mathbf{H}$ such that*

$$(a_-, b_-)(a_+, b_+) = (a, b)$$

*then the rational functions $b_-$ and $b_+$ are subordinate to $b$.*

**Proof.** Since $(a,b) \in \mathbf{L}$, we have $aa^* = 1 + bb^*$. Since $a^*$ does not vanish in $D \cup \mathbf{T}$, we see that $a$ is subordinate to $bb^*$ on a neighborhood of $D \cup \mathbf{T}$. Since $a$ has no poles on $D^*$, it is subordinate to $bb^*$ on $D^*$ and thus on the whole Riemann sphere.

Now let $(a_-, b_-)$ and $(a_+, b_+)$ be a Riemann-Hilbert factorization as in the lemma. Then

$$b_- = -ab_+ + ba_+$$
$$b_+ = -\frac{b_-}{a} + b\frac{a_+}{a}$$

On $D^* \cup \mathbf{T}$, the functions $b_-/a$ and $a_+/a$ have no poles, since they are in $H^2(D^*)$. The last display then implies that $b_+$ is subordinate to $b$ on a neighborhood of $D^* \cup \mathbf{T}$. Since $b_+$ has no poles on $D$, it is subordinate to $b$ on the whole Riemann sphere.

Similarly one proves that $b_-$ is subordinate to $b$. □

**Lemma 25.** *Let $(a, b) \in \mathbf{L}$ be rational. Then there exists a unique Riemann-Hilbert factorization*

$$(a, b) = (a_-, b_-)(a_+, b_+) \tag{40}$$

*such that $b_+$ does not have any poles on $\mathbf{T}$. The factor $(a_+, b_+)$ coincides with the factor $(a_{++}, b_{++})$ in the triple factorization of Theorem 8. Similarly, there exists a unique Riemann-Hilbert factorization (40) such that $b_-$ does not have any poles on $\mathbf{T}$. For this factorization, the factor $(a_-, b_-)$ coincides with the factor $(a_{--}, b_{--})$ in the triple factorization of Theorem 8.*

**Proof.** If there exists a Riemann-Hilbert factorization (40) such that $b_+$ has no pole on $\mathbf{T}$, then the factor $(a_+, b_+)$ has to coincide with $(a_{++}, b_{++})$. Namely, it is clear that $(a_+, b_+)$ has no further nontrivial Riemann-Hilbert factorization by Lemma 18. This implies that $(a_+, b_+)$ is equal to $(a_{++}, b_{++})$.

In particular we have proved that the requirement that $b_+$ has no poles on $\mathbf{T}$ makes the Riemann-Hilbert factorization unique.

It remains to show that such a Riemann-Hilbert factorization exists.

We set up a Banach fixed point argument as in the proof of Lemma 18. Recall from the proof of that lemma that, for every constant $c$, the affine linear mapping

$$T : (A, B) \to (c + P_{D^*}(\frac{b^*}{a^*}B), P_D(\frac{b}{a}A))$$

is a weak contraction mapping in the sense that

$$\|T(A, B) - T(A', B')\| \leq \|(A, B) - (A', B')\|$$

where the norms are with respect to $L^2(\mathbf{T}) \oplus L^2(\mathbf{T})$. Indeed, unless $(A', B') = (A, B)$, the inequality is strict since multiplication by $b/a$ strictly lowers the $L^2$ norm of any nonzero element. This implies that on any invariant finite dimensional subspace of $L^2(\mathbf{T})$, the mapping is a strict contraction in the sense

$$\|T(A, B) - T(A', B')\| \leq (1 - \epsilon)\|(A, B) - (A', B')\|$$

for some $\epsilon$ depending on the subspace. This can be seen by a compactness argument.

Consider the finite dimensional space $V$ of all rational $(A, B)$ such that $B$ is subordinate to $b$, $A$ is subordinate to $b^*$, and $A$ and $B$ have no poles on $\mathbf{T}$. This is clearly a subspace of $L^2(\mathbf{T}) \oplus L^2(\mathbf{T})$. For any rational function $f$ without poles on $\mathbf{T}$ the projection $P_D f$ is up to an additive constant the sum of the principal parts of the poles of $f$ in $D^*$, while $P_{D^*} f$ is up to an additive constant the sum of the principal parts of the poles of $f$ in $D$. Therefore, $P_D f$ is a rational function subordinate to $f$ with no poles in $D$ and $P_{D^*} f$ is a rational function subordinate to $f$ with no poles in $D^*$.

We observe that for $(a,b) \in \mathbf{H}$ the quotient $b/a$ has no poles on $\mathbf{T}$ and is subordinate to $b$ on $D^*$. Thus for any $(A,B) \in V$ we have that $P_{D^*}(b^*B/a^*)$ is subordinate to $b^*$ and $P_D(bA/a)$ is subordinate to $b$. Thus $V$ is invariant under the mapping $T$ for any $c$. Since $T$ is a strict contraction mapping on $V$, there exists a fixed point in $V$ under this mapping. Using this fixed point $(A,B)$ for $c=1$, we can as in the proof of Lemma 18 produce a right factor

$$(a_+, b_+) = A(\infty)^{-1}(A,B)$$

to the Riemann-Hilbert factorization problem for $(a,b)$. Clearly $b_+$ has no poles on $\mathbf{T}$, so this is the desired right factor.

The symmetric statement concerning left factors is proved similarly. This completes the proof of Lemma 25. □

The above shows that there is a very satisfactory description of the set of rational elements in $\mathbf{L}$ which qualify to be left, middle, or right factors in a triple factorization as in Theorem 8. Namely, possible left (middle, right) factors are exactly those rational $(a,b) \in \mathbf{L}$ for which $b$ has only poles in $D$, ($\mathbf{T}$, $D^*$).

It remains to study the possible factorizations of a rational middle factor in the triple factorization. Thus we are reduced to study the Riemann-Hilbert problem for rational $(a,b) \in \mathbf{H}_0^* \cap \mathbf{H}$. Any factorization consists again of rational factors in $\mathbf{H}_0^* \cap \mathbf{H}$.

This problem too has a very satisfactory answer, though the formulation of the answer is a little more involved.

Before we proceed further, we shall briefly digress on the maximum principle. The maximum principle says that any nonconstant holomorphic function on $D$ which is continuous on $D \cup \mathbf{T}$ attains its maximum only on $\mathbf{T}$. If the function is actually differentiable on $D \cup \mathbf{T}$, then the following lemma gives more precise information. It is a version of the maximum principle which may be less well known.

**Lemma 26.** *Let $f$ be a nonconstant holomorphic map from $D$ to itself and assume that $f$ and $f'$ have continuous extensions to the boundary $\mathbf{T}$. Thus $f$ and $f'$ map $D \cup \mathbf{T}$ to $D \cup \mathbf{T}$.*

*If $f$ attains its maximum at $z \in \mathbf{T}$, then $f'(z) = z^*\omega f(z)$ for some strictly positive $\omega$.*

**Proof.** Multiplying $f$ by a constant phase factor, if necessary, we may assume $f(z) = 1$. Consider the real part $u$ of $f$. It has a maximum at $z$. In particular, $u$ has zero derivative in the direction tangential to $\mathbf{T}$. Therefore, the gradient of $u$ has to be radial and is either 0 or outward pointing. This proves

$$\frac{\partial u}{\partial x} + i\frac{\partial u}{\partial y} = \omega z$$

for some $\omega \geq 0$ Thus, by the Cauchy Riemann equations,

$$f'(z) = \frac{\partial u}{\partial x} + i\frac{\partial v}{\partial x} = \frac{\partial u}{\partial x} - i\frac{\partial u}{\partial y} = \omega z^* = \omega z^* f(z)$$

It remains to show that $\omega$ is not zero, i.e., that the harmonic function $u$ does not have vanishing derivative at $z$. Assume by a rotation that $z=1$. It will suffice to find some function $\tilde{u}$ which dominates $u$ in the intersection of $D$ with a neighborhood of 1, such that $\tilde{u}$ is differentiable at 1 with nonvanishing derivative.

Since $u$ is not constant, we find two points on $\mathbf{T}$ where $u$ is strictly less than 1. The two points divide the circle into wo arcs $C_1$ and $C_2$. Let $L$ be the line connecting the two points. Assume w.l.o.g. that $C_1$ contains 1 and let $z_0$ be a point of $C_2$. Define

$$\tilde{u}(\zeta) = 1 + \epsilon \mathrm{Re} \frac{\zeta + z_0}{\zeta - z_0}$$

for some small $\epsilon > 0$. Since $\tilde{u}$ is 1 on the arc $C_1$, it dominates $u$ there. Since $u$ is strictly less than 1 in the (compact) line $L$, we can choose $\epsilon$ small enough so that $\tilde{u}$ dominates $u$ on the line. By the (easy) maximum principle, $\tilde{u}$ dominates $u$ inside the enterior of $C_2 \cup L$. It remains to prove that $\tilde{u}$ has nonvanishing derivative at 1. This however can be done easily by direct inspection. □

We continue to study Riemann-Hilbert factorizations

$$(a, b) = (a_-, b_-)(a_+, b_+)$$

for rational $(a, b) \in \mathbf{H} \cap \mathbf{H}_0^*$. Thus $b$ and $a$ have only poles on $\mathbf{T}$. Indeed, they have the same poles as $1 + bb^* = aa^*$ shows.

By Lemma 24, the functions $a_+$ and $a_-$ can only have poles where $a$ has poles. Consider the identity

$$a = a_- a_+ (1 + \frac{b_-}{a_-} \frac{\overline{b_+}}{a_+}) \tag{41}$$

The function

$$\frac{b_-}{a_-} \frac{\overline{b_+}}{a_+} \tag{42}$$

maps $D \cup \mathbf{T}$ to itself. Therefore, the last factor on the right-hand side of (41) can only vanish at point $z \in \mathbf{T}$ when $z$ is a maximum of the function (42) on $D \cup \mathbf{T}$. By Lemma 26, the last factor in (41) can only have a simple zero at $z$.

Therefore, for every pole $z$ of $a$, we have either

$$ord(a, z) = ord(a_-, z) + ord(a_+, z) \tag{43}$$

or

$$ord(a, z) = ord(a_-, z) + ord(a_+, z) - 1 \tag{44}$$

We say that the pole $z$ is *split* if (43) holds, and we say that it is shared if (44) holds. If $z$ is a shared pole, then both functions $b_-/a_-$ and $b_+^*/a_+$ have modulus one at $z$. Therefore, both functions $a_+$ and $a_-$ have a pole at $z$ and by (44) both poles have order at most $\mathrm{ord(a, z)}$.

For each pole $z$ of $a$ we define

$$n := \mathrm{ord(a, z)}, \quad \mathrm{n^-} := \mathrm{ord(a_-, z)}, \quad \mathrm{n_j^+} := \mathrm{ord(a_+, z)}$$

Define the functions

$$A_+ := 1 - \frac{b_+ b^*}{a_+ a^*} = \frac{a_-^*}{a_+ a^*} = \frac{1}{a_+ a_+^*} \frac{1}{1 + \frac{b_-^*}{a_-^*} \frac{b_+}{a_+^*}}$$

$$A_- := 1 - \frac{b_- b^*}{a_-^* a} = \frac{a_+}{a_-^* a} = \frac{1}{a_- a_-^*} \frac{1}{1 + \frac{b_-}{a_-} \frac{b_+^*}{a_+}}$$

On $\mathbf{T}$, the functions $A_+$ and $A_-$ have positive real part except possibly where $a$ has a pole (use the first representation for $A_+, A_-$). There, $A_+$ vanishes of order $n + n_- - n_+$ and $A_-$ vanishes of order $n + n_+ - n_-$ (use the second representation for $A_+, A_-$). In particular, $A_+$ vanishes of order $2n^+$ if the pole is split or $2n^+ - 1$ if the pole is shared.

For each shared pole $z$, we define $\mu^+$ and $\mu^-$ by the asymptotic expansions

$$A_+(\zeta) = -\mu^+ z(\zeta - z)^{n^+-1}(\frac{1}{\zeta} - \frac{1}{z})^{n^+} + O(\zeta - z)^{2n^+}$$

$$A_-(\zeta) = -\mu^- z^*(\zeta - z)^{n^-}(\frac{1}{\zeta} - \frac{1}{z})^{n^- -1} + O(\zeta - z)^{2n^-})$$

We claim that $\mu^+$ and $\mu^-$ are positive.

To see this for $A_+$, we set

$$a_+(\zeta) = \gamma(\zeta - z)^{-n^+} + O(\zeta - z)^{-n^+ + 1}$$

and, by Lemma 26,

$$(1 + \frac{b_-^* b_+}{a_-^* a_+^*})'(z) = (\frac{b_-^* b_+}{a_-^* a_+^*})'(z) = (\frac{b_-^* b_+}{a_-^* a_+^*})(z) z^* \mu = -z^* \mu$$

for some positive $\mu$. Using the third representation of $A_+$ we obtain

$$\mu^+ = 1/(\mu|\gamma|^2)$$

awhich shows that $\mu^+$ is positive. The proof that $\mu^-$ is positive is similar.

Write

$$\frac{1}{aa^*}(\zeta) = \mu(\zeta - z)^n(\frac{1}{\zeta} - \frac{1}{z})^n + O(\zeta - z)^{2n+1}$$

Then the identity

$$A_+ A_- = \frac{1}{aa^*}$$

shows that

$$\mu^+ \mu^- = \mu$$

Our goal is to see that the parameters $n^+$ and $n^-$ for all poles together with the parameters $\mu^+$ and $\mu^-$ for all shared poles parameterize the Riemann-Hilbert factorizations of $(a, b)$. We shall first adress the easier statement that all Riemann-Hilbert factorizations are uniquely determined by these parameters.

**Lemma 27.** *Let $(a, b)$ be rational in $\mathbf{H} \cap \mathbf{H}_0^*$. Then any Riemann-Hilbert factorization*

$$(a, b) = (a_-, b_-)(a_+, b_+)$$

*is uniquely determined if the parameters $n^+$ and $n^-$ for all poles and the parameters $\mu^+$, and $\mu^-$ for all shared poles are specified.*

**Proof.** We assume to get a contradiction that there are two Riemann Hilbert factorizations with right factors $(a_+, b_+)$ and $(\tilde{a}_+, \tilde{b}_+)$ respectively, which have the same parameters listed in the lemma.

Define $(c, d)$ by

$$(\tilde{a}_+, \tilde{b}_+) = (c, d)(a_+, b_+)$$

Our task is to show that $(c, d) = (1, 0)$. It suffices to show that $d$ is constant. Then $d$ has to be constant 0 as one can see from evaluating the defining equation for $(c, d)$ at 0 and using that $b_+, \tilde{b}_+$ vanish at 0 while $a_+^*$ does not. Then $c = \tilde{a}_+/a_+$ is an outer function on $D^*$ and of constant modulus 1 on $\mathbf{T}$, and thus it is constant. This constant is equal to 1 as one can see from evaluating $c$ at $\infty$.

Observe that $d$ is a rational function and can only have poles where $a$ has poles. Therefore, it suffices to show that $d$ is holomorphic at all poles of $a$.

Fix a pole $z$. Define $r := b/a$ and similarly $r_+, r_-, \tilde{r_+}$. Then we have

$$r - r_+ = \frac{1}{r^*}(1 - r^* r_+ - (1 - r^* r)) = \frac{1}{r^*}(A_+ - (1 - r^* r))$$

Observe that $r^*$ has modulus one at $z$ and $1 - rr^* = (aa^*)^{-1}$ vanishes of order $2n$ at $z$. Moreover, $A_+$ vanishes at least of order $2n^+ - 1$ at $z$ and its Taylor coefficient of order $2n^+$ at $z$ is determined by $\mu^+$.

Therefore, $r - r_+$ vanishes at least of order $2n^+ - 1$ at $z$ and its Taylor coefficient of order $2n$ is determined by $r$ and $\mu^+$. The same holds for $r - \tilde{r}_+$, and by taking differences we see that $r_+ - \tilde{r}_+$ vanishes of order $2n^+$ at $z$. Since $a_+$ has a pole of order $n^+$ at $z$, we see that

$$(r_+ - \tilde{r}_+)a_+\tilde{a}_+ = b_+\tilde{a}_+ - \tilde{b}_+ a_+ = d$$

has no pole at $z$. □

**Theorem 9.** *Assume $(a, b) \in \mathbf{H} \cap \mathbf{H}_0^*$ is rational. Let $z_j \in \mathbf{T}$, $j = 1, \dots, N$ be the distinct poles of $a$ and denote the order of the pole $z_j$ by $n_j$.*

*Assume we are given numbers $0 \leq n_j^+, n_j^- \leq n_j$ for $j = 1, \dots, N$ such that for each $j$ either*

$$n_j^+ + n_j^- = n_j$$

*(split case) or*

$$n_j^+ + n_j^- - 1 = n_j$$

*(shared case). Assume further that for each $j$ in the shared case we are given positive real numbers $\mu_j^+, \mu_j^-$ with*

$$\mu_j^+ \mu_j^- = \mu_j$$

*where $\mu_j$ is defined by*

$$\frac{1}{aa^*}(\zeta) = \mu_j(\zeta - z_j)^{n_j}(\frac{1}{\zeta} - \frac{1}{z_j})^{n_j} + O(\zeta - z_j)^{2n_j+1}$$

*Then there exists a unique Riemann-Hilbert factorization*

$$(a, b) = (a_-, b_-)(a_+, b_+)$$

*such that*

$$\mathrm{ord}(a_+, z_j) = n_j^+$$

$$\mathrm{ord}(a_-, z_j) = n_j^-$$

*and, if $j$ is in the shared case,*

$$A_+(\zeta) = -\mu_j^+ z_j (\zeta - z_j)^{n_j^+ - 1}(\frac{1}{\zeta} - \frac{1}{z_j})^{n_j^+} + O((\zeta - z_j)^{2n_j^+})$$

$$A_-(\zeta) = -\mu_j^- z_j^* (\zeta - z_j)^{n_j^-}(\frac{1}{\zeta} - \frac{1}{z_j})^{n_j^- - 1} + O((\zeta - z_j)^{2n_j^-})$$

*All Riemann-Hilbert factorizations of $(a, b)$ are obtained in this way.*

**Proof.** Our previous discussion of the parameters $n_j^+$,$n_j^-$,$\mu_j^+$ and $\mu_j^-$ already implies that every Riemann-Hilbert factorization of $(a, b)$ comes with parameters as described in the theorem and the parameters determine the factorization uniquely.

It thus remains to show that for a given set of parameters such a Riemann-Hilbert factorization exists.

It is enough to consider the case when $b$ is nonzero at $\infty$, because

$$(a, b) = (a_-, b_-)(a_+, b_+)$$

is equivalent to

$$(a, bz^n) = (a_-, b_-z^n)(a_+, b_+z^n)$$

and thus one can reduce the case of $b$ vanishing at $\infty$ to the case of $b$ not vanishing at $\infty$.

We first prove existence of a Riemann-Hilbert factorization in the easier case when all poles are split. We write

$$b(z) = b(\infty) \left[\prod(z - y_k)\right] \left[\prod(z - z_j)^{-n_j}\right]$$

where $y_k$ are the zeros of $b$ counted with multiplicities.

For each $j$, consider points $z_j^+ \in D^*$ and $z_j^- \in D$ close to $z_j$. It shall be enough to consider $z_j^\pm$ such that they avoid the zeros of $b$ and are all pairwise distinct as well as distinct from all $(z_j^\pm)^*$. Consider the perturbation $\tilde{b}$ of $b$ defined by

$$\tilde{b}(z) = b(\infty) \left[\prod(z - y_k)\right] \left[\prod(z - z_j^+)^{-n_j^+}\right] \left[\prod(z - z_j^-)^{-n_j^-}\right]$$

This function has the same zeros with multiplicities as $b$, but the poles are at perturbed locations.

Since the zeros $y_k$ are fixed, we have an upper bound on

$$\int_{\mathbf{T}} \log_+ |\tilde{b}|$$

uniformly in the choice of the points $z_j^\pm$.

By Lemma 22, there is a unique rational $\tilde{a}$ such that $(\tilde{a}, \tilde{b}) \in \mathbf{L}$. The equation $1 + \tilde{b}\tilde{b}^* = \tilde{a}\tilde{a}^*$ implies that there is a uniform upper bound on

$$\int_{\mathbf{T}} \log_+ |\tilde{a}|$$

and since $\tilde{a}$ is outer on $D^*$ we have a uniform upper bound on $\tilde{a}(\infty)$. Trivially, we also have the lower bound $1 \le |\tilde{a}(\infty)|$.

Since $\tilde{a}$ is bounded on $\mathbf{T}$, Lemma 18 gives a unique Riemann-Hilbert factorization

$$(\tilde{a}, \tilde{b}) = (\tilde{a}_-, \tilde{b}_-)(\tilde{a}_+, \tilde{b}_+)$$

Applying Lemma 24 thoroughly to the situation at hand, we conclude that $\tilde{a}_+$ has poles only at the points $(z_j^+)^*$ with order at most $n_j^+$, while $\tilde{a}_-$ has poles only at the points $z_j^-$ with order at most $n_j^-$. Moreover, since $|\tilde{a}_+(\infty)|$ and $|\tilde{a}_-(\infty)|$ are bounded below by one and $\tilde{a}_+(\infty)\tilde{a}_-(\infty) = \tilde{a}(\infty)$, we obtain that $\tilde{a}_+(\infty)$ and $\tilde{a}_-(\infty)$ are in a fixed compact set avoiding 0.

We can write

$$\tilde{a}_+(z) = \tilde{a}_+(\infty) \left[\prod(z - x_k^+)\right] \left[\prod(z - (z_j^+)^*)^{-\tilde{n}_j^+}\right]$$

where $x_k^+$ are the zeros of $\tilde{a}_+$ with multiplicities and $\tilde{n}_j^+ \le n_j^+$.

Now we consider a sequence of choices of $z_j^{\pm}$ such that for each $j$ both points $z_j^+$ and $z_j^-$ converge to $z_j$. For ea Since the zeros $x_k^+$ remain in the compact set $D \cup T$ and the value $\tilde{a}_+(\infty)$ remains in a compact set away from 0, there is a subsequence for which the $\tilde{n}_j^+$ are constant, each zero $x_k^+$ converges (assuming the zeros are appropriatelyenumerated), and the value $\tilde{a}_+(\infty)$ converges. For this subsequence, $\tilde{a}_+$ converges uniformly on compact sets away from the poles of $a$ to a limit $a_+$. The zeros of $a_+$ are still in $D \cup \mathbf{T}$ an the poles are on $\mathbf{T}$. Thus $a_+$ is still outer. Clearly also $a_+(\infty)$ is positive as a limit of positive numbers.

Similarly, one can choose a further subsequence so that all other terms in the identity

$$(\tilde{a}, \tilde{b}) = (\tilde{a}_-, \tilde{b}_-)(\tilde{a}_+, \tilde{b}_+)$$

converge uniformly on compact sets away from the poles of $a$. In the limit, we obtain a Riemann-Hilbert factorization

$$(a, b) = (a_-, b_-)(a_+, b_+)$$

By construction, the poles of $a_+$ and of $a_-$ at $z_j$ are at most of order $n_j^+$ and $n_j^-$ respectively. Since the sum of these orders for any Riemann-Hilbert factorization has to be at least $n_j$, the orders of $a_+$ and $a_-$ at $z_j$ are exactly $n_j^+$ and $n_j^-$. Thus we have proved existence of a Riemann-Hilbert factorization for the given parameters in the completely split case.

Now we modify the above argument so that it works in the case when there are shared poles.

For each split pole $z_j$, we choose again $z_j^+$ and $z_j^-$ exactly as before. For every shared pole, we choose $z_j^+$ and $z_j^-$ as before but with the additional constraint that

$$\mu_j^+ |(z_j^+)^* - z_j|^{2n_j^+} = \mu_j^- |z_j^- - z_j|^{2n_j^-} \tag{45}$$

We define

$$\tilde{b}(z) = b(\infty) \left[\prod (z - y_k)\right] \left[\prod_{z_j \text{ shared}} (z - z_j)\right] \left[\prod (z - z_j^+)^{-n_j^+}\right] \left[\prod (z - z_j^-)^{-n_j^-}\right]$$

Compared to the completely split case, we have defined $\tilde{b}$ to have an additional zero at each shared pole $z_j$. Since for each shared pole we have $n_j^+ + n_j^- = n_j + 1$, the numerator and denominator of the rational function defining $\tilde{b}$ have the same degree and $\tilde{b}(\infty)$ is again finite.

As before, we obtain a unique Riemann-Hilbert factorization

$$(\tilde{a}, \tilde{b}) = (\tilde{a}_-, \tilde{b}_-)(\tilde{a}_+, \tilde{b}_+)$$

Then we let $z_j^{\pm}$ tend to $z_j$ respecting the additional constraint (45) for each shared pole. As before, we can choose a subsequence so that all quantities in the Riemann-Hilbert factorization converge uniformly on compact sets away from the poles $z_j$. In the limit, we obtain a Riemann-Hilbert factorization

$$(a, b) = (a_-, b_-)(a_+, b_+)$$

We need to show that this factorization has the given parameters $n_j^{\pm}$ and $\mu_j^{\pm}$. As before, for each split pole $z_j$ the order of poles of the limits $a_+$ and $a_-$ are at most $n_j^+$ and $n_j^-$ and thus have to be exactly $n_j^+$ and $n_j^-$.

We consider a shared pole $z_j$. We calculate

$$\tilde{b} = \tilde{a}_- \tilde{b}_+ + \tilde{b}_- \tilde{a}_+^*$$

$$\frac{\tilde{b}}{\tilde{a}_- \tilde{a}_+^*} = \frac{\tilde{b}_+}{\tilde{a}_+^*} + \frac{\tilde{b}_-}{\tilde{a}_-}$$

At the shared pole $z_j$, the left-hand side vanishes:

$$0 = \frac{\tilde{b}_+(z_j)}{\tilde{a}_+^*(z_j)} + \frac{\tilde{b}_-(z_j)}{\tilde{a}_-(z_j)} \tag{46}$$

Since for every element $(a', b') \in \mathbf{L}$, the modulus of $b'/a'$ on $\mathbf{T}$ can be expressed in terms of the modulus of $a$, we conclude from (46) that

$$|\tilde{a}_+(z_j)| = |\tilde{a}_-(z_j)| \tag{47}$$

In a compact neighborhood of $z_j$ avoiding the poles $z_k$ with $k \neq j$, we can write

$$\tilde{a}_+(z) = (z - (z_j^+)^*)^{-n_j^+} \tilde{h}^+(z)$$

$$\tilde{a}_-(z) = (z - z_j^-)^{-n_j^-} \tilde{h}^-(z)$$

where $\tilde{h}^+$ and $\tilde{h}^-$ converge uniformly on the neighborhood to functions $h^+$ and $h^-$. The equation (47) then becomes

$$|z_j - (z_j^+)^*|^{-n_j^+} |\tilde{h}^+(z_j)| = |z_j - z_j^-|^{-n_j^-} |\tilde{h}^+(z_j)|$$

By choice of the $z_j^\pm$ we thus have for some constant $\tilde{c}$:

$$|\tilde{h}^+(z_j)| = \tilde{c}(\mu_j^+)^{-1/2}$$

$$|\tilde{h}^-(z_j)| = \tilde{c}(\mu_j^-)^{-1/2}$$

Taking the limit, we obtain

$$|h^+(z_j)| = c(\mu_j^+)^{-1/2}$$

$$|h^-(z_j)| = c(\mu_j^-)^{-1/2}$$

for some constant $c$. We claim that $c$ is not zero. Assume to get a contradiction that it is zero, then $h^+$ and $h^-$ vanish at $z_j$. From the equations

$$a_+(z) = (z - z_j)^{-n_j^+} h^+(z)$$

$$a_-(z) = (z - z_j)^{-n_j^-} h^-(z)$$

we see that $a_+$ and $a_-$ have order of pole at most $n_j^+ - 1$ and $n_j^- - 1$ at $z_j$. The sum of these orders is less than $n_j$, a contradiction. Therefore $c$ is not zero.

Then $h^+$ and $h^-$ do not vanish at $z_j$ and $a_+$ and $a_-$ have poles of exact order $n_j^+$ and $n_j^-$ at $z_j$. Thus we are indeed in the case of a shared pole.

Moreover, we calculate in the limit

$$\frac{|a_+|^2}{|a_-|^2} = \frac{\mu_j^-}{\mu_j^+} |z - z_j|^{-2n_j^+ + 2n_j^-} + O(|z - z_j|^{-2n_j^+ + 2n_j^- + 1})$$

Now we use $\mu_j^+ \mu_j^- = \mu_j$ and

$$\frac{1}{|a|^2} = \mu_j^2 |z - z_j|^{2n_j} + O(|z - z_j|^{2n_j + 1})$$

to obtain

$$\frac{|a_+|^2}{|a_-|^2|a|^2} = |\mu_j^-|^2|z - z_j|^{-2n_j^+ + 2n_j^- + 2n_j} + O(|z - z_j|^{-2n_j^+ + 2n_j^- + 2n_j + 1}$$

Comparing with the asymptotics for $A_-$ and doing the analogue calculation for $A_+$ we conclude that the shared pole $z_j$ has indeed the parameters $\mu_j^+$ and $\mu_j^-$. □

# LECTURE 5
# Orthogonal polynomials

## 1. Orthogonal polynomials

In this lecture, we describe how the nonlinear Fourier transform on the half-line relates to orthogonal polynomials. The material on orthogonal polynomials is folklore, a standard reference is [**21**] and a more recent introduction with interesting applications is [**6**].

Let $\mu$ be any compactly supported positive measure on the plane $\mathbf{C}$ with the normalization $\|\mu\| = 1$.

Let $H$ be the Hilbert space completion of the linear span of the set of functions $z^0 = 1, z^1, z^2, \dots$ (monomials) under the inner product

$$\langle f, g \rangle = \int f\overline{g}\, d\mu$$

If a finite set $z^0, \dots, z^n$ of monomials is linearly dependent in this Hilbert space, then necessarily $\mu$ has finite support. Namely, linear dependence means that some linear combination of these monomials, which is nothing but a polynomial $P(z)$, is equivalent to 0 in the Hilbert space. This means

$$\|P\| = \int |P(z)|^2\, d\mu = 0$$

As $|P|^2$ and $d\mu$ are positive, this can only happen if the support of $\mu$ is contained in the null set of $P$, which is finite. Conversely, if $\mu$ has finite support, then it is easy to find a polynomial which vanishes on the support and thus is equivalent to 0 in the Hilbert space.

From now on the standing assumption is that $\mu$ has infinite support.

We can apply the Gram-Schmidt orthogonalization process to the sequence of vectors $z^n$. Let

$$\phi_0 = z^0 = 1$$

and let $\phi_n$ be the unique polynomial of degree $n$ which has unit length in the Hilbert space $H$, has positive highest coefficient, and is orthogonal to all polynomials of degree less than $n$. Since the space of polynomials of degree less than $n$ has exact codimension 1 in the space of polynomials of degree less than or equal to $n$, such a polynomial $\phi_n$ exists and clearly has exact degree $n$. The orthogonality condition determines $\phi_n$ up to a scalar factor. The modulus of this scalar factor is determined by the requirement that $\phi_n$ has unit length, and the phase is determined by the requirement that the highest coefficient of $\phi_n$, which is the coefficient in front of $z^n$, is positive.

The polynomials $\phi_n$ form an orthonormal set. Indeed, they form an orthonormal basis of $H$ since they span the same subspace as the monomials $z^n$, which by definition span the full space $H$.

The linear operator

$$T : f \to zf$$

originally defined on all polynomials $f$, is bounded with respect to the norm on $H$ because the function $z$ is bounded on the support of $\mu$. Therefore, the operator $T$ extends to a unique bounded operator on $H$.

We can express the map $T$ in the basis $\phi_n$. This means we represent elements in $H$ as infinite linear combinations $\sum a_j \phi_j$ and let $T$ act on the column vector $(a_j)_{j\geq 0}$ by matrix multiplication from the left by a matrix $(J_{ij})_{i,j\geq 0}$:

$$(a_j) \to (\sum_j J_{ij} a_j)$$

As $z\phi_j$ is a polynomial of degree $j+1$, we have $J_{ij} = 0$ if $i > j+1$. Moreover, $z\phi_j$ has positive highest coefficient, and thus $J_{j+1,j} > 0$. We call a matrix a Hessenberg matrix if it satisfies these two constraints, namely positivity on the subdiagonal and vanishing below the subdiagonal.

The requirement that Hessenberg matrices be positive on the subdiagonal is a matter of convenient normalization to produce uniqueness results. One can conjugate a Hessenberg matrix by a diagonal unitary matrix and thus obtain an equivalent matrix which vanishes below the subdiagonal but is merely nonzero on the subdiagonal. In the literature, one often calls these more general matrices Hessenberg matrices.

We will specialize to two cases.

Case 1: The measure $\mu$ is supported on the real line $\mathbf{R}$. In this case $T : f \to zf$ is selfadjoint since

$$\langle Tf, g\rangle = \int zf\overline{g}\, d\mu = \int f\overline{zg}\, d\mu = \langle f, Tg\rangle$$

where we have used that $z$ is real on the support of $\mu$. Thus $J$ is a selfadjoint matrix, and since it is Hessenberg, it is tridiagonal ($J_{ij} = 0$ if $|i-j| > 1$). A Hessenberg matrix which is self adjoint is called a Jacobi matrix.

Observe that the $j$-th column vector of a Jacobi matrix has two real parameters not determined by the previous columns: the subdiagonal is a positive number by definition of Hessenberg matrices and independent of the previous columns, the diagonal entry has to be a real number by selfadjointness and is independent of the previous columns, while the superdiagonal entry is determined by the previous column.

We know that $T$ and thus $J$ have to be bounded operators. For tridiagonal matrices $J_{ij}$, boundedness is equivalent to $\sup_{ij} |J_{ij}| < \infty$.

Case 2: The measure is supported on the circle $\mathbf{T}$. In this case we have $T^*T = \mathrm{id}$ since

$$\langle T^*Tf, g\rangle = \langle Tf, Tg\rangle = \int zf\overline{zg}\, d\mu = \int f\overline{g}\, d\mu = \langle f, g\rangle$$

Observe that $T^*T = \mathrm{id}$ does not imply that $T$ is unitary, since $TT^* = \mathrm{id}$ may fail. For example, if the measure $\mu$ is a normalized Lebesgue measure on $\mathbf{T}$, then the monomials $z^n$ form an orthonormal set. The operator $T$ in this case is the standard shift operator and is obviously not surjective since the image of $T$ contains only functions with zero constant term.

Unlike in the selfadjoint case, the matrix $J$ is not sparse above the diagonal. However, we still have two new real parameters per column. The condition $J^*J = \mathrm{id}$ implies that the column vectors are orthonormal. Thus the first $n$ entries of the $n$-th column have to be orthogonal to the previous columns, and thus there is only a complex parameter in $D$ (the vector of the first $n$ entries must have norm less than 1) as degree of freedom, and then the $(n+1)$-st entry is determined since it is positive and makes the column have unit length.

Observe that boundedness of a matrix satisfying $J^*J = \mathrm{id}$ is automatic.

**Theorem 10.** *The above construction of a Jacobi matrix provides a bijective correspondence between*

1. *A compactly supported positive measure $\mu$ supported on the real line with $||\mu|| = 1$ and infinite support.*
2. *A tridiagonal selfadjoint matrix $J = (J_{ij})_{i,j\geq 0}$ that has strictly positive elements on the subdiagonal and has finite operator norm.*

**Proof.** We need to show that the map from measures to Jacobi matrices described above is injective and surjective.

The matrix $J$ determines the moments of $d\mu$ because

$$\int z^n d\mu = \langle T^n \phi_0, \phi_0 \rangle = \langle J^n e_0, e_0 \rangle$$

where $e_n$ denotes the $n$-th standard unit vector whose $n$-th entry is 1 and all other entries are 0. But then by the Stone-Weierstrass theorem, this determines $d\mu$. Hence the map $\mu \mapsto J$ is injective.

To show surjectivity, let $J$ be a Jacobi matrix and write down the Neumann series

$$(J-z)^{-1} = -\sum_{n=0}^{\infty} \frac{J^n}{z^{(n+1)}}$$

Since $J$ is bounded, for $|z|$ greater than the operator norm of $J$ the Neumann series converges and provides a proper meaning to the left-hand side of the last equation.

We can form the function

$$m(z) = \frac{1}{\pi} \left\langle (J-z)^{-1} e_0, e_0 \right\rangle$$

This function is holomorphic in a neighborhood of $\infty$ on the Riemann sphere since the Neumann series converges there. Moreover,

$$m(z) = -1/z + O(z^{-2})$$

as $z \to \infty$.

We would like to show that $m$ can be extended to the open upper half-plane and has positive imaginary part there. Let $z$ have positive imaginary part and let $z_0$ be purely imaginary with sufficiently large (in modulus) negative imaginary part. Then the Neumann series

$$(J-z)^{-1} = ((J-z_0) - (z-z_0))^{-1} = \sum_{n=0}^{\infty} -\frac{(J-z_0)^n}{(z-z_0)^{n+1}}$$

converges because

$$|z_0 - z|^2 = |z_0 - \mathrm{Im}\,(z)|^2 + |\mathrm{Re}\,(z)|^2 = |z_0|^2 + 2|z_0||\mathrm{Im}\,(z)| + O(1)$$

while

$$\|z_0 - J\|^2 = \sup_{\|f\|=1} \langle (J - z_0)f, (J - z_0)f \rangle$$

$$= \sup_{\|f\|=1} \langle z_0 f, z_0 f \rangle + \langle Jf, Jf \rangle = |z_0|^2 + O(1)$$

This shows that both $(J - z)^{-1}$ and $m$ extend holomorphically to the upper half-plane.

For $\mathrm{Im}\,(z) > 0$, let $\phi_z := (J - z)^{-1} e_0$. We observe

$$\mathrm{Im}\,(m(z)) = \frac{1}{\pi} \langle \phi_z, (J - z)\phi_z \rangle = \frac{1}{\pi} \mathrm{Im}\,(z) \|\phi_z\|^2$$

Thus $m$ has positive imaginary part in the upper half plane.

By the Herglotz representation theorem there is a positive measure $\mu$ on the real line such that

$$m(z) = \frac{1}{\pi} \int \frac{1}{\zeta - z}\, d\mu(\zeta)$$

The measure is compactly supported on $\mathbf{R}$ since $m$ is holomorphic in a neighborhood of $\infty$.

As $m(z)$ has the asymptotics

$$z^{-1} + O(z^{-2})$$

near $\infty$, we have $\|\mu\| = 1$.

We prove that $\mu$ is not supported on a finite set. The most important ingredient in the following argument is that $J$ has no zeros on the subdiagonal $i = j + 1$.

Assume to get a contradiction that the measure is supported on a finite set. Then the integral representation for $m$ shows that $m$ is a rational function and therefore there is a polynomial

$$p(z) = \sum_{k=0}^{m} a_k z^k$$

such that $p(z)m(z)$ is a polynomial. The latter polynomial has the form

$$-\frac{1}{\pi} \sum_{k=0}^{m} a_k z^k \sum_{n=0}^{\infty} z^{-(n+1)} \langle J^n e_0, e_0 \rangle$$

Since this is a polynomial, for sufficiently large $n$ the coefficient in the Laurent series has to be zero. This leads to the identity

$$\sum_{k=0}^{m} a_k \left\langle J^{n+k} e_0, e_0 \right\rangle = 0$$

for large $n$. By selfadjointness of $J$, we also obtain for large $n$ and all $m \geq 0$

$$\sum_{k=0}^{m} a_k \left\langle J^{n+k} e_0, J^m e_0 \right\rangle = 0$$

Since the vectors $J^m e_0$ with $m \geq 0$ span the full Hilbert space $l^2(\mathbf{Z}_{\geq 0})$, this gives

$$\sum_{k=0}^{m} a_k J^{n+k} e_0 = 0$$

This is however absurd if $J$ is a Hessenberg matrix with positive elements on the subdiagonal.

To complete the proof of surjectivity, it remains to show that the measure $\mu$ that we have constructed gives rise to the matrix $J$. By construction we have for large $|z|$

$$\int (\zeta - z)^{-1} \, d\mu(\zeta) = \langle (J-z)^{-1} e_0, e_0 \rangle$$

Comparing coefficients in the Neumann series gives for $n \geq 0$

$$\int \zeta^n \, d\mu(\zeta) = \langle J^n e_0, e_0 \rangle$$

By self adjointness of $J$ we also have

$$\int \zeta^n \, d\mu(\zeta) = \langle J^n e_0, J^m e_0 \rangle$$

The same identity holds with $J$ replaced by the matrix $\tilde{J}$ constructed from the measure $\mu$ by the Gram-Schmidt orthogonalization process. Therefore, it suffices to show that the numbers

$$\langle J^n e_0, J^m e_0 \rangle \tag{48}$$

determine $J$. We claim that knowing these numbers, we can express all of the standard basis vectors $e_n$ as linear combination of $J^0 e_0, \ldots J^n e_0$. This is clear for $n = 0$. Assume by induction we can express $e_0, \ldots, e_n$ in terms of $J^0 e_0, \ldots, J^n e_0$. We can calculate the coefficients $\langle J^{n+1} e_0, e_j \rangle$ for $j \leq n$ by expressing $e_j$ in terms of $J^0 e_0, \ldots, J^j e_j$ and using (48). As $J^{n+1} e_0$ is a linear combination of $e_0, \ldots, e_{n+1}$, it remains to determine the coefficient in front of $e_{n+1}$. This coefficient is positive, hence we can then determine this coefficient by determining the length of $J^{n+1} e_0$, which we can do by (48). This poves the claim.

Using the claim, we can calculate $\langle J e_n, e_m \rangle$ for all $n, m$ and thus obtain all coefficients of $J$. This completes the proof of surjectivity.

□

In the case of measures on the circle we obtain

**Theorem 11.** *The above construction of a Hessenberg matrix provides a bijective correspondence between*

1. *A positive measure $\mu$ supported on the circle $\mathbf{T}$ with $||\mu|| = 1$ and infinite support.*
2. *A matrix $J = (J_{ij})_{i,j \geq 0}$ which satisfies $J_{ij} = 0$ for $i > j+1$, $J_{(j+1),j} > 0$ for all $j$, and $J^* J = \mathrm{id}$.*

**Proof.** The proof is a reprise of the arguments in the proof of the previous theorem. We shall only describe the inverse map and leave the other details as an exercise.

We can use Neumann series to define

$$\frac{z+J}{z-J} = (z+J) \sum_{n=0}^{\infty} J^n z^{-(n+1)}$$

for $|z| > 1$. Then we define

$$m(z) = \left\langle \frac{z+J}{z-J} e_0, e_0 \right\rangle$$

on $D^*$. Setting

$$\phi_z := (z - J)^{-1} e_0$$

we have

$$\operatorname{Re} m(z) = \operatorname{Re} \langle (z + J)\phi_z, (z - J)\phi_z \rangle$$
$$= \langle z\phi_z, z\phi_z \rangle - \langle J\phi_z, J\phi_z \rangle = (|z|^2 - 1)\|\phi_z\|^2$$

Thus $m$ has positive real part on $D^*$ and by the Herglotz representation theorem is the extension of a unique positive measure $\mu$:

$$m(z) = \int_{\mathbf{T}} \frac{z + \zeta}{z - \zeta} \, d\mu$$

This is the desired measure. □

Assume $J$ is a matrix with $J^*J = 1$ and $J_{ij} = 0$ for $i - 1 > j$. Then we have $JJ^*J = J$ and thus $JJ^* = 1$ on the image of $J$. We claim that the image of $J$ has codimension at most 1, and in case that the codimension is exactly 1, then $e_0$ complements the image:

$$\operatorname{image}(J) + \operatorname{span}(e_0) = l^2(\mathbf{Z}_{\geq 0})$$

The claim simply follows from the fact that

$$e_0, Je_0, J^2 e_0, \ldots, J^n e_0$$

evidently span the same space as

$$e_0, e_1, e_2, \ldots e_n$$

If the codimension of the image is 0, then $J$ is unitary. This happens if and only if $e_0$ is in the image of $J$. A criterion when this happens can be derived from Szegö's theorem stated below. Namely, $J$ is unitary if and only if both sides of the identity in Szegö's theorem are zero.

**Theorem 12.** *(Szegö) Let $\mu$ be a measure on $\mathbf{T}$ with a.c. part $w d\theta$ (we denote by $d\theta$ Lebesgue measure on $\mathbf{T}$ normalized so that $\int_{\mathbf{T}} d\theta = 1$). Then*

$$\inf_f \int |1 - f|^2 \, d\mu = \exp \int_{\mathbf{T}} \log |w| \, d\theta$$

*where $f$ runs through all polynomials in $z$ with zero constant term and $\|f\| = 1$ in $L^2(\mu)$.*

In case the left-hand side is nonzero, it has the meaning of $|\langle u, 1\rangle|^2$ where $u$ is the unit vector perpendicular to the image of $J$. The right-hand side becomes zero if $\int \log |w| = -\infty$.

We will prove this theorem in the next section in the case that both sides are finite.

## 2. Orthogonal polynomials on T and the nonlinear Fourier transform

In this section we relate orthogonal polynomials on the circle $\mathbf{T}$ to the nonlinear Fourier transform on the half-line.

The nonlinear Fourier transform $(a, b)$ of a sequence in $l^2(\mathbf{Z}_{\geq 1}, D)$ gives rise to an analytic function $b/a^*$ which maps $D$ to itself and vanishes at 0. Via a Möbius transform of the target space $D$, such an analytic functions is in unique correspondence with an analytic function $m$ on $D$ which has positive real part and is equal to 1 at 0:

$$m(z) = \frac{1 - b/a^*}{1 + b/a^*}$$

By the Herglotz representation theorem, this function is uniquely associated with a positive measure on $\mathbf{T}$ with total mass 1. By Theorem 11, assuming for the moment that this measure does not have finite support, the measure is in unique correspondence with a Hessenberg matrix. The following theorem states that the diagonal and subdiagonal entries of the Hessenberg matrix can easily be expressed in terms of the sequence $F$.

**Theorem 13.** *Let $F \in l^2(\mathbf{Z}_{\geq 1}, D)$ and let $(a, b)$ be the nonlinear Fourier transform of $F$. Let $\mu$ be the positive measure on $\mathbf{T}$ whose harmonic extension to $D$ is equal to the real part of*

$$m(z) = \frac{1 - b/a^*}{1 + b/a^*}$$

*Then the Hessenberg matrix associated to $\mu$ via Theorem (11) satisfies*

$$J_{i(i-1)} = (1 - |F_i|^2)^{1/2} \tag{49}$$

$$J_{ii} = -F_i \overline{F_{i+1}} \tag{50}$$

*for $i \geq 1$ and*

$$J_{00} = -\overline{F_1} \tag{51}$$

The main point of this theorem is that taking the forward nonlinear Fourier transform is equivalent to calculating the spectral data (the measure $\mu$) of a Hessenberg matrix, while taking the inverse nonlinear Fourier transform or "layer stripping method" is equivalent to the Gram-Schmidt orthogonalization process.

**Proof.** For $F$ a sequence in $l^2(\mathbf{Z}_{\geq 1}, D)$ with nonlinear Fourier transform $(a, b)$, we consider the truncated sequences $F_{\leq n}$ and their nonlinear Fourier transforms

$$(a_{\leq n}, b_{\leq n}) = (a_n, b_n)$$

Define for $n \geq 0$

$$\phi_n(z) = z^n[a_n(z) + b_n^*(z)] \tag{52}$$

By Lemma 2, $\phi_n$ is a polynomial in $z$ of exact degree $n$, and the highest coefficient of $\phi_n$ is equal to the constant coefficient of $a_n$ and therefore it is positive. In particular, $\phi_0 = 1$.

We can write (52) in matrix form (with our standing convention to complete the second row of a matrix by applying the $*$ operation to and reversing the order of the entries of the first row) as

$$(z^{-n}\phi_n, z^n\phi_n^*) = (1, 1)(a_n, b_n) \tag{53}$$

Thus we obtain a recursion formula by multiplying a transfer matrix from the right:

$$(z^{-(n+1)}\phi_{n+1}, z^{n+1}\phi_{n+1}^*) = (z^{-n}\phi_n, z^n\phi_n^*)(1 - |F_{n+1}|^2)^{-1/2}(1, F_{n+1}z^{n+1})$$

In particular, we have the identity

$$z^{-(n+1)}\phi_{n+1} = (1 - |F_{n+1}|^2)^{-1/2}[z^{-n}\phi_n + \overline{F_{n+1}}z^{-1}\phi_n^*]$$

which can be rewritten as

$$z\phi_n = (1 - |F_{n+1}|^2)^{1/2}\phi_{n+1} - \overline{F_{n+1}} z^n \phi_n^* \tag{54}$$

This recursion formula (54) together with $\phi_0 = 1$ contains all information on the sequence $\phi_n$, but it will be convenient to rewrite the recursion in a different form.

In (53), we may also multiply from the right by an inverse transfer matrix to obtain

$$(z^{-(n-1)}\phi_{n-1}, z^{n-1}\phi_{n-1}^*) = (z^{-n}\phi_n, z^n\phi_n^*)(1 - |F_n|^2)^{-1/2}(1, -F_n z^n)$$

$$z^{-(n-1)}\phi_{n-1} = (1 - |F_n|^2)^{-1/2}[z^{-n}\phi_n - \overline{F_n}\phi_n^*]$$

$$\phi_n = (1 - |F_n|^2)^{1/2} z\phi_{n-1} + \overline{F_n} z^n \phi_n^* \tag{55}$$

Let the collection of polynomials $\phi_n$ be fomally an orthonormal basis of a Hilbert space $H$. Thus abstractly the Hilbert space is the set of all linear combinations of the $\phi_n$ with square summable coefficients, and the inner product is given by the standard inner product on the space of square summable sequences. Further below we will identify the measure on $\mathbf{T}$ with respect to which the $\phi_n$ are the Gram-Schmidt orthogonal polynomials, but for now we only need $H$ abstractly defined.

We claim that $T : f \to zf$, originally defined on the set of polynomials, extends to an isometry on $H$. This will follow once we have shown that the $\phi_n$ are the Gram-Schmidt orthogonal polynomials of a measure supported on $\mathbf{T}$. However, we find the following proof instructive, and the proof will also help to calculate the diagonal and subdiagonal entries of the Hessenberg matrix representing $T$.

We shall prove by induction the following two statements

$$\|z\phi_{n-1}\| = 1 \tag{56}$$

$$\langle z\phi_{n-1}, \phi_n\rangle = (1 - |F_n|^2)^{1/2} \tag{57}$$

Observe that equation (52) applied for $n = -1$ gives $\phi_{-1} = z^{-1}$. In this sense the above two statements for $n = 0$ reduce to the fact $\|\phi_0\| = 1$.

Assume by induction that (56) and (57) are true for some $n$. Then by pairing (55) with $z\phi_{n-1}$ we see that the two summands on the right-hand side of (55) are orthogonal. As the coefficients in (55) form a Pythagorean triple, we conclude

$$\|z^n\phi_n^*\| = 1 \tag{58}$$

Using this in (54) together with the obvious fact that $\phi_{n+1}$ is orthogonal to $z^n\phi_n^*$, we obtain $\|z\phi_n\| = 1$ and $\langle z\phi_n, \phi_{n+1}\rangle = (1 - |F_{n+1}|^2)^{1/2}$. This concludes the induction step.

Now we prove by induction on $n$ that

$$\langle z\phi_n, z\phi_m\rangle = 0 \tag{59}$$

for $m < n$. Fix $n$ and assume that (59) is true for all indices smaller than $n$. By eliminating $z^n\phi_n^*$ from (54) and (55), $z\phi_n$ becomes a linear combination of terms manifestly orthogonal to $z\phi_m$ if $m < n - 1$. If $m = n - 1$, we use (55) and (57) to observe that $z^n\phi_n^*$ is perpendicular to $z\phi_{n-1}$ and we insert this into (54). This proves (59) for the index $n$.

Therefore, the orthonormal basis $\phi_n$ is mapped to an orthonormal set via the operation $T : f \to zf$. This proves that $T$ is an isometry.

The operator $T$ in the basis $\phi_n$ is given by a Hessenberg matrix $J$, i.e., $J_{ij} = 0$ for $i > j+1$ and $J_{j+1,j} > 0$ (recall the highest coefficient of $\phi_n$ is positive).

Indeed, we read from the above recursions that the Hessenberg matrix satisfies (49), (50), and (51). Only the proof of (50) is slightly more involved. It follows from adding $\overline{F_n}$ times (54) and $\overline{F_{n+1}}$ times (55) to obtain

$$\overline{F_n} z\phi_n + \overline{F_{n+1}}\phi_n = \overline{F_n}(1 - |F_{n+1}|^2)^{1/2}\phi_{n+1} + \overline{F_{n+1}}(1 - |F_n|^2)^{1/2} z\phi_{n-1}$$

Taking the $\phi_n$ component everywhere gives

$$\overline{F_n}\,\langle z\phi_n, \phi_n\rangle + \overline{F_{n+1}} = \overline{F_{n+1}}(1 - |F_n|^2)$$

Which proves (50).

We know from the abstract theory discussed in the previous section that there is a positive measure $\mu$ on the circle with $\|\mu\| = 1$ such that the $\phi_n$ are the Gram-Schmidt orthogonal polynomials with respect to this measure. We need to show that the Herglotz function $m$ of this measure satisfies

$$m(z) = \frac{1 - s(z)}{1 + s(z)}$$

where we have set $s = b/a^*$.

We first argue that it suffices to prove the claim under the additional assumption that the sequence $F$ is compactly supported. To pass to the case of general $F$, we approximate $s$ by $s_n$. Clearly $s_n$ converges on the disc $D$ pointwise to $s$ since we know convergence of $(a_n, b_n)$ to $(a, b)$ in $\mathbf{H}$. Thus the measures $\mu_n$ (defined as above by the function $s_n = b_n/a_n^*$) converge weakly to the measure $\mu$ (defined as above by the function $s$). Observe that the polynomials $\phi_m$ defined as $z^m(a_m + b_m^*)$ do only depend on the first $m$ coefficients of $F$ and thus are the same as if defined by the truncated sequences as long as the truncation parameter $n$ is larger than or equal to $m$.

Suppose we can prove that the inner product $\langle \phi_m, \phi_{m'}\rangle$ with respect to the measure $\mu_n$ for all sufficiently large $n$ is equal to the Kronecker delta of $m$ and $m'$. By passing to the weak limit, the inner product with respect to $\mu$ is also the Kronecker delta. Thus the $\phi_n$ are an orthonormal basis for $\mu$.

We now prove the orthonormality property of $\phi_n$ under the additional assumption that $F_n$ is a finite sequence. Then $a$, $b$ and $m$ are analytic across the circle $\mathbf{T}$, and the measure $\mu$ associated to the Herglotz function $m$ is absolutely continuous with respect to normalized Lebesgue measure and has density $\mathrm{Re}\,(m)$. We can write

$$\begin{aligned}\mathrm{Re}\,(m) &= \frac{1}{2}\left[\frac{1-s}{1+s} + \frac{1-s^*}{1+s^*}\right]\\ &= \frac{1 - ss^*}{(1+s)(1+s^*)}\\ &= \frac{1}{(a^*+b)(a+b^*)}\end{aligned}$$

We need to show that with respect to this measure,

$$\phi_n = z^n(a_n + b_n^*)$$

is orthogonal to all $z^k$ with $k < n$ and has length 1.

Consider the decomposition (we write $a_+$ for $a_{>n}$ etc.)

$$\begin{pmatrix} a & b \end{pmatrix} = \begin{pmatrix} a_n & b_n \end{pmatrix}\begin{pmatrix} a_+ & b_+ \end{pmatrix}$$

or equivalently

$$\begin{pmatrix} a & b \end{pmatrix} \begin{pmatrix} a_+^* & -b_+ \end{pmatrix} = \begin{pmatrix} a_n & b_n \end{pmatrix}$$

We obtain

$$a_n + b_n^* = aa_+^* - bb_+^* + b^* a_+^* - a^* b_+^*$$
$$= (a + b^*) a_+^* - (a^* + b) b_+^*$$

And thus

$$\int_{\mathbf{T}} \frac{z^{-k} z^n (a_n + b_n^*)}{(a^* + b)(a + b^*)}$$

$$(60) \qquad = \int_{\mathbf{T}} z^{n-k} \left[ \frac{a_+^*}{a^* + b} - \frac{b_+^*}{a + b^*} \right]$$

However,

$$\frac{a_+^*}{a^* + b}$$

is holomorphic on $D$ and

$$\frac{b_+^*}{a + b^*}$$

is holomorphic on $D^*$ with a zero of order $n+1$ at $\infty$. Thus the above integral (60) is zero for $0 \leq k < n$.

For $n = k$ we write for (60)

$$\int_{\mathbf{T}} \frac{a_+^*}{a^* + b} - \frac{b_+^*}{a + b^*}$$

$$= \frac{a_+^*(0)}{a^*(0) + b(0)} = \frac{1}{a_n(\infty)}$$

As the highest order coefficient of $\phi_n$ is $a_n(\infty)$, we obtain

$$\|\phi_n\|^2 = a_n(\infty) \langle z^n, \phi_n \rangle = 1$$

Thus we have verified that $\mu$ is the measure with respect to which the $\phi_n$ are orthonormal.

We remark that if $F$ is a finite sequence and $n$ is the order of the largest nonzero element of $F$, then the density of the measure is given by

$$\frac{1}{(a_n^* + b_n)(a_n + b_n^*)} = \frac{1}{|\phi_n|^2}$$

□

We are now ready to prove Szegö's theorem under the assumption that $F_n$ is a square integrable sequence and thus $\log |a|$ is integrable.

The measure $\mu$ splits as $\mu = \mu_s + \mu_{ac}$ where $\mu_s$ is singular with respect to the Lebesgue measure and $\mu_{ac}$ is absolutely continuous with density $w \in L^1(\mathbf{T})$.

Since $\mathrm{Re}\,(m)$ is equal to $w$ almost everywhere on $\mathbf{T}$, we have

$$\int \log w = \int \log\left(\frac{1}{(a^* + b)(a + b^*)}\right)$$

$$= \int \log \frac{1}{|a|^2} + \int \log \frac{1}{1 + s} + \int \log \frac{1}{1 + s^*}$$

The last two integrals vanish because the function $1/(1+s)$ is holomorphic and outer on $D$ and equal to 1 at 0. Thus

$$\int \log w = -2 \int \log |a|$$

On the other hand, $\inf_f \|1-f\|$ on the left-hand side of Szegö's theorem is equal to the modulus squared of the inner product $\langle u, 1\rangle$ where $u$ is a unit vector perpendicular to the image of $T : f \to zf$. (We shall momentarily establish that under the assumption of square integrable $F_n$ this image has indeed codimension 1.) We claim that the vector $u$ is the strong limit as $n \to \infty$ of the vectors

$$z^n \phi_n^*$$

These vectors are unit vectors by (58) and perpendicular to $z\phi_0, \ldots, z\phi_{n-1}$ by (55). Thus the strong limit of these vectors has to be a unit vector perpendicular to the image of $T : f \to zf$.

To show the existence of the strong limit, we observe that applying the star operation to (55) and multiplying by $z^n$ we obtain

$$z^n \phi_n^* = (1-|F_n|^2)^{1/2} z^{n-1} \phi_{n-1}^* - F_n \phi_n \tag{61}$$

This recursion implies that

$$z_n \phi_n^* = \sum_{k=0}^{n} c_{k,n} F_k \phi_k$$

for some constants $c_{k,n}$ bounded by 1. Using orthogonality of the vectors $\phi_n$ we obtain that the sequence $(z_n \phi_n^*)$ has a limit $u$.

Using (61) and orthogonality of $\phi_n$ to $\phi_0$ for $n > 1$, we obtain

$$\langle u, 1\rangle = \prod_{i=1}^{\infty} (1-|F_n|^2)^{1/2}$$

Thus the left-hand side of Szegö's theorem is $\prod_{i=1}^{\infty}(1-|F_n|^2)$ and the identity of Szegö's theorem follows by the nonlinear Plancherel identity.

This proves Szegös theorem in the setting of square summable $F$. It can be shown that in case $F$ is not square summable, orthogonal polynomials can still be defined using the Hessenberg matrix associated to the sequence $F$, and both sides of Szegö's theorem vanish.

## 3. Jacobi matrices and the nonlinear Fourier transform

Orthogonal polynomials on $\mathbf{R}$ and Jacobi matrices are related to the nonlinear Fourier transform on $l^2(\mathbf{Z}_{\geq 1}, [-1,1])$, i.e., the space of bounded real valued sequences $F_n$ on the the half-line. Since Jacobi matrices are a discrete model for Schrödinger operators, the material discussed here also relates to the spectral theory of Schrödinger operators. We only touch upon the subject, entry points to some related current literature are [**13**], [**5**], [**18**].

Our main concern is again to show that on the one hand one can parameterize a class of Jacobi matrices easily by sequences $F \in l^2(\mathbf{Z}_{\geq 1}, [-1,1])$, and on the other hand one can determine the measure $\mu$ associated to such a Jacobi matrix easily from the nonlinear Fourier transform $(a, b)$ of $F$.

We first observe that if $F$ is a sequence in $l^2(\mathbf{Z}, [-1, 1])$, then the nonlinear Fourier transform $(a, b)$ of $F$ satisfies

$$a(z^*) = a^*(z), \quad b(z^*) = b^*(z) \tag{62}$$

by property (9) stated in Lemma 1. Conversely, any element $(a, b) \in \mathbf{H}$ satisfying (62) has a real sequence $F \in l^2(\mathbf{Z}_{\geq 0}, D)$ as preimage under the nonlinear Fourier transform.

**Theorem 14.** *Let $F \in l^2(\mathbf{Z}_{\geq 1}, [-1, 1])$ and let $(a, b)$ be the nonlinear Fourier transform of $F$. Let $\mu$ be the positive measure on $[-2, 2]$ whose harmonic extension to the upper half plane is the imaginary part of the function $m$ defined by*

$$m(w + w^*) = \frac{1}{w - w^*} \frac{1 - b(w)/a^*(w)}{1 + b(w)/a^*(w)}$$

*Then the Jacobi matrix associated to $\mu$ via Theorem 10 satisfies*

$$J_{n,n} = (F_{2n+1}(1 + F_{2n}) - F_{2n-1}(1 - F_{2n})) \tag{63}$$

$$J_{n+1,n} = (1 + F_{2n})^{1/2}(1 - |F_{2n+1}|^2)^{1/2}(1 - F_{2n+2})^{1/2} \tag{64}$$

*for $n \geq 1$ and*

$$J_{1,0} = 2^{1/2}(1 - |F_1|^2)^{1/2}(1 - F_2)^{1/2} \tag{65}$$

$$J_{0,0} = 2F_1 \tag{66}$$

We remark that the elements of the sequence $F_n$ with even $n$ enter into the formulas for $J_{ij}$ in a different manner from the elements with odd index. Interesting special cases occur when either all $F_n$ with even index $n$ vanish, or all $F_n$ with odd index $n$ vanish, but we do not further elaborate on this here.

**Proof.** We study orthogonal polynomials for measures supported on the interval $[-2, 2]$. By a simple scaling argument, it is no restriction if we fix the length of this interval. One can relate these polynomials to orthogonal polynomials on the circle $\mathbf{T}$ using the conformal map

$$w \to y = w + w^*$$

from $D$ to the Riemann sphere slit at $[-2, 2]$. This map is sometimes called the Joukowski map. The Joukowski map extends to the boundary $\mathbf{T}$ of $D$ and maps $\mathbf{T} \setminus \{-1, 1\}$ two-to-one onto the interval $(-2, 2)$ and maps $\{-1, 1\}$ to one to one to the endpoints $\{-2, 2\}$ of that interval.

Using this map on $\mathbf{T}$, one can push forward measures on $\mathbf{T}$ to measures on $[-2, 2]$. This provides a bijection from measures $\mu'$ on $\mathbf{T}$ with the symmetry $\mu'(w) = \mu'(w^*)$ to measures $\mu$ on $[-2, 2]$. The relation between $\mu'$ and $\mu$ is given by the formula

$$\int f(w + w^*)\, d\mu'(w) = \int f(y)\, d\mu(y)$$

We shall assume that the measure $\mu'$ is positive and normalized to $\|\mu'\| = 1$. Then the same normalization holds for $\mu$.

By the Herglotz representation theorem, there is a holomorphic function $m'$ on $D$ whose real part is the harmonic extension of $\mu'$. By the Herglotz representation theorem for the upper half-plane, there is a function $m$ on the upper half plane whose imaginary part is the harmonic extension of the compactly supported measure $\mu$ on $\mathbf{R}$. We claim

$$m'(w) = (w - w^*)m(w + w^*)$$

Namely, using the Poisson kernel for the disc, we have

$$m'(w) = \int \frac{v + w}{v - w}\, d\mu'(v)$$

Using the symmetry of $\mu'$, we obtain

$$\begin{aligned} m'(w) &= \frac{1}{2}\int_{\mathbf{T}} \frac{v + w}{v - w} + \frac{v^* + w}{v^* - w}\, d\mu'(v) \\ &= \int_{\mathbf{T}} \frac{1 - w^2}{(v - w)(v^* - w)}\, d\mu'(v) \\ &= \int_{\mathbf{T}} \frac{w^* - w}{(vw^* - 1)(v^* - w)}\, d\mu'(v) \\ &= (w - w^*)\int_{\mathbf{T}} [(v + v^*) - (w + w^*)]^{-1}\, d\mu'(v) \\ &= (w - w^*)\int_{\mathbf{R}} [y - (w + w^*)]^{-1}\, d\mu(y) \\ &= (w - w^*)m(w + w^*) \end{aligned}$$

This proves the claim.

We now discuss how the orthogonal polynomials in the variable $y = w + w^*$ can be identified with orthogonal polynomials in the variable $w$ on $\mathbf{T}$ with respect to $\mu'$. As in the previous section, let $\phi_n$ denote the orthogonal polynomials in the variable $w$ on $\mathbf{T}$ with respect to the measure $\mu'$. Thus

$$\phi_n(w) = w^n[a_n(w) + b_n^*(w)]$$

where $(a_n, b_n)$ are the truncated nonlinear Fourier transforms of a sequence $F_n$ supported on $\mathbf{Z}_{\geq 1}$. Due to the symmetry of $\mu'$ we have

$$\phi^*(w) = \phi(w^*), \quad a^*(w) = a(w^*), \quad b^*(w) = b(w^*)$$

and the sequence $F_n$ is real.

As $\phi_{2n}$ is orthogonal to all monomials of degree up to $2n - 1$, the function

$$\psi_n(w) := (w^*)^n\phi_{2n}(w)$$

is orthogonal to all functions $w^k$ with $-n \leq k \leq n - 1$. Consequently, $\psi_n$ is also orthogonal to all

$$(w + w^*)^k \quad , \quad 0 \leq k \leq n - 1 \tag{67}$$

By symmetry under $w \to w^*$, also $\psi_n^*(w)$ is orthogonal to all (67), and so is

$$\Psi_n = \psi_n + \psi_n^* = [w^n(a + b^*) + w^{-n}(a^* + b)] \tag{68}$$

However, $\Psi_n$ is itself symmetric under $w \to w^*$, and thus a polynomial in $y = w + w^*$ of degree $n$. Therefore, up to an as of yet unspecified scalar factor, $\Psi_n$ is the $n$-th orthogonal polynomial with respect to $\mu$.

To determine the scalar factor we calculate

$$\int \Psi_n(w)\Psi_n^*(w)\, d\mu'$$

$$= \int (\psi_n + \psi_n^*)(\psi_n + \psi_n^*)$$
$$= 2 + 2\mathrm{Re} \int \psi_n^2 \, d\nu$$
$$= 2 + 2\mathrm{Re} \int w^{-2n} \phi_{2n} \phi_{2n}$$

If $n = 0$, this is simply equal to 4.

Assume $n > 0$. Since $\phi_{2n}$ is an orthogonal polynomial in the variable $w$ with respect to $\mu'$, we see that the last display is equal to

$$= 2 + 2\mathrm{Re} \int \phi_{2n}(w) w^{-2n} c_0 \tag{69}$$

where $c_0$ is the constant coefficient of $\phi_{2n}$, which can be obtained from (55) and (61)

$$c_0 = \overline{F_{2n}} \prod_{k=1}^{2n} (1 - |F_k|^2)^{-1/2}$$

However, again by the fact that $\phi_{2n}$ is an orthogonal polynomial, (69) is equal to

$$= 2 + 2\mathrm{Re} \int \phi_{2n}(w) \phi_{2n}^*(w) c_0 (\overline{c_{2n}})^{-1}$$

where $c_{2n}^{-1}$ is the highest coefficient of $\phi_{2n}$. From (61) and the value of $c_0$ stated above we obtain

$$c_n = \prod_{k=1}^{2n} (1 - |F_k|^2)^{-1/2}$$

Thus, since $F_n$ is real,

$$\|\Psi_n\|^2 = 2(1 + F_n)$$

Define

$$\Phi_0 = 1$$

and, for $n \geq 1$,

$$\Phi_n = 2^{-1/2}(1 + F_{2n})^{-1/2} \Psi_n$$

Then $\Phi_n$ is the $n$-th orthogonal polynomial with respect to $\mu$. Observe that the expression for $\Phi_n$ in the case of generic $n$ remains correct if we set $F_0 = 1$.

We calculate the Jacobi matrix associated to the polynomials $\Phi_n$. Assume first $n \geq 1$. We can write for (68)

$$(\Psi_n, \Psi_n) = (1, 1)(a, b)(w^n, 0)(1, 1)$$

Then we have the recursion equations

$$(1 - |F_{2n+1}|^2)^{1/2}(1 - |F_{2n+2}|^2)^{1/2}(\Psi_{n+1}, \Psi_{n+1})$$
$$= (1, 1)(a, b)(1, F_{2n+1} w^{2n+1})(1, F_{2n+2} w^{2n+2})(w^{n+1}, 0)(1, 1)$$
$$= (1, 1)(a, b)(w^n, 0)(w, F_{2n+1})(1, F_{2n+2})(1, 1)$$
$$= (1, 1)(a, b)(w^n, 0)([w + F_{2n+1}][(1 + F_{2n+2}], 0)(1, 1) \tag{70}$$

In the last step we have used that for any real number $\gamma$ and any $(a, b)$ we have

$$(a, b)(\gamma, \gamma) = (a + b, 0)(\gamma, \gamma)$$

Similarly,

$$(1 - |F_{2n-1}|^2)^{1/2}(1 - |F_{2n}|^2)^{1/2}(\Psi_{n-1}, \Psi_{n-1})$$

$$= (1,1)(a,b)(1,-F_{2n}w^{2n})(1,-F_{2n-1}w^{2n-1})(w^{n-1},0)(1,1)$$
$$= (1,1)(a,b)(w^n,0)(1,-F_{2n})(w^*,-F_{2n-1})(1,1)$$
$$= (1,1)(a,b)(w^n,0)(w^* - F_{2n}w - F_{2n-1} + F_{2n}F_{2n-1},0)(1,1) \tag{71}$$

Multiplying (70) by $(1+F_{2n})/(1+F_{2n+2})$ and adding to (71) we obtain on the right-hand side

$$(1,1)(a,b)(w^n,0)(w+w^*+F_{2n+1}(1+F_{2n}) - F_{2n-1}(1-F_{2n}),0)(1,1)$$
$$= (w+w^*+F_{2n+1}(1+F_{2n}) - F_{2n-1}(1-F_{2n}))(\Psi_n,\Psi_n)$$

where we have pulled a real scalar matrix out of the product. Collecting the terms and expressing $\Psi_n$ in terms of $\Phi_n$ gives for $n \geq 2$:

$$(1+F_{2n})^{1/2}(1-|F_{2n+1}|^2)^{1/2}(1-F_{2n+2})^{1/2}\Phi_{n+1}$$
$$+(1+F_{2n-2})^{1/2}(1-|F_{2n-1}|^2)^{1/2}(1-F_{2n})^{1/2}\Phi_{n-1}$$
$$= (w+w^*)\Phi_n + (F_{2n+1}(1+F_{2n}) - F_{2n-1}(1-F_{2n}))\Phi_n$$

This identity shows that we can express multiplication by $y = w + w^*$ in the basis $\Phi_n$ by a matrix $J$ with

$$J_{n,n} = (F_{2n+1}(1+F_{2n}) - F_{2n-1}(1-F_{2n}))\Phi_n$$

and

$$J_{n+1,n} = (1+F_{2n})^{1/2}(1-|F_{2n+1}|^2)^{1/2}(1-F_{2n+2})^{1/2}$$

for $n \geq 1$. To obtain the value for $J_{1,0}$, we review the above calculation for $n = 1$ and observe that it remains correct if $\Phi_0$ is replaced by $2^{1/2}\Phi_0$ in view of the special normalization of $\Phi_0$. Thus

$$J_{21} = 2^{1/2}(1-|F_1|^2)^{1/2}(1-F_2)^{1/2}$$

To calculate $J_{0,0}$, we specialize (70) to $n = 0$:

$$2^{1/2}(1-|F_1|^2)^{1/2}(1-F_2)^{1/2}(\Phi_1,\Phi_1)$$
$$= (1,1)(w+F_1,0)(1,1)$$

Or

$$2^{1/2}(1-|F_1|^2)^{1/2}(1-F_2)^{1/2}\Phi_1 = (w+w^*)\Phi_0 + 2F_1\Phi_0$$

Thus

$$J_{0,0} = 2F_1$$

□

# LECTURE 6
# Further applications

## 1. Integrable systems

The linear Fourier transform takes partial differential operators into multiplication operators by coordinate functions. Hence the linear Fourier transform takes simple partial differential equations such as linear equations with constant coefficients into algebraic equations. The latter can often be solved by explicit expressions modulo the task of taking the linear Fourier transform and inverting it. In this section we discuss that the nonlinear Fourier transform can be used similarly to obtain explicit solutions to certain nonlinear partial differential equations, again modulo the task of taking the nonlinear Fourier transform and inverting it. We will present the calculations on a purely formal and thus expository level. They can be made rigorous in appropriate function spaces, e.g., by the precise calculus in [**1**]. The formal calculations for the particular example of the modified Korteweg de Vries equation that we choose for the exposition can be found in [**23**]. A more general discussion can be found in [**12**].

We need the nonlinear Fourier transform of functions on $\mathbf{R}$ and we shall briefly introduce it.

Recall that the linear Fourier transform of sequences is defined by

$$\widehat{F}(\theta) = \sum_{n \in \mathbf{Z}} F_n e^{-in\theta}$$

Here $\theta$ lives on the interval $[-\pi, \pi] \subset \mathbf{R}$. To pass to the continuous Fourier transform on $\mathbf{R}$, one can do a limiting process by letting $\theta \in [-\pi/\epsilon, \pi/\epsilon]$ and $n \in \epsilon\mathbf{Z}$ with $\epsilon$ approaching 0. Taking an appropriate limit, one obtains the Fourier transform of a function $F$ on $\mathbf{R}$

$$\widehat{F(k)} = \int_{\mathbf{R}} F(x) e^{2ikx}\, dx \tag{72}$$

Here we have used a special normalization of $2ikx$ in the exponent of the exponential function, which is maybe unusual but convenient for the discussions to follow.

Now consider $F \in l^2(\mathbf{Z}, D)$, its truncations $F_{\leq n}$ and their nonlinear Fourier transforms $\widehat{F_{\leq n}} = (a_n, b_n)$. Then we have the recursion equation

$$(a_n(z), b_n(z)) = (a_{n-1}(z), b_{n-1}(z)) \frac{1}{1 - |F_n|^2}(1, F_n z^n)$$

Subtracting $(a_{n-1}, b_{n-1})$ on both sides we obtain

$$(a_n(z), b_n(z)) - (a_{n-1}(z), b_{n-1}(z))$$
$$= (a_{n-1}(z), b_{n-1}(z)) \left[ \frac{1}{1 - |F_n|^2}(1, F_n z^n) - (1, 0) \right]$$

A similar type of limiting process as in the linear case, leads to an expression for the nonlinear Fourier transform on $\mathbf{R}$. The discrete variable $n$ becomes a continuous variable $x \in \mathbf{R}$, and the variable $z$ becomes $e^{2ik}$ for some real $k$, and we obtain

$$\frac{\partial}{\partial x}(a(k,x), b(k,x)) = (a(k,x), b(k,x))(0, F(x)e^{2ikx})$$

In case of compactly supported $F$, solutions $(a, b)$ to this ordinary differential equation are constant to the left and to the right of the support of $F$. We denote these constant values by $(a(k, -\infty), b(k, -\infty))$ and $(a(k, \infty), b(k, \infty))$. To obtain the nonlinear Fourier transform, we set the initial value condition

$$(a(k, -\infty), b(k, -\infty)) = (1, 0)$$

and then define

$$\widehat{F}(k) = (a(k, \infty), b(k, \infty))$$

We review how the linear Fourier transform is used to solve linear constant coefficient PDE. Our example is the Cauchy problem for the Airy equation. Thus the problem is to find a solution $F(x, t)$ to the Airy equation

$$F_t = F_{xxx}$$

(where a variable in the index denotes a partial derivative) with the initial condition

$$F(0, x) = F_0(x)$$

for some given function $F_0$. Taking formally the linear Fourier transform of $F$ in the $x$ variable one obtains a function $\widehat{F}(t, k)$ satisfying

$$\widehat{F}_t = (-2ik)^3 \widehat{F} = 8ik^3 \widehat{F}$$
$$\widehat{F}(0, k) = \widehat{F_0}(k)$$

For fixed $k$ this is an ordinary differential equation in $t$ which has the solution

$$\widehat{F}(t, k) = e^{8ik^3 t} \widehat{F_0}(k)$$

Taking the inverse Fourier transform, one obtains

$$F(t, x) = (e^{8ik^3 t} \widehat{F_0}(k))^{\vee}$$

Aanlogously, the nonlinear Fourier transform can be used to solve certain nonlinear partial differential equations. As an example we discuss the Cauchy problem for the modified Korteweg-de Vries (mKdV)

(73) $$F_t = F_{xxx} + 6F^2 F_x$$
$$F(0, x) = F_0(x)$$

Observe that the mKdV equation is a perturbation of the Airy equation by the nonlinear term $6F^2F_x$.

We take the nonlinear Fourier transform of the initial data:

$$\widehat{F_0}(k) = (a(k), b(k))$$

Then the solution to the mKdV equation is formally given by

$$\widehat{F}(t,k) = (a(k), e^{8ik^3t}b(k)) \tag{74}$$

Therefore, safe for the task of taking a nonlinear Fourier transform and an inverse nonlinear Fourier transform, this is an explicit solution.

We outline a proof of (74) on a formal level. The argument can be made rigorous in approoriate function spaces.

A Lax pair is a pair of time dependent differential operators $L(t), P(t)$ in spatial variables such that $L$ is selfadjoint, $P$ is anti-selfadjoint, and

$$\frac{d}{dt}L(t) = [P(t), L(t)] = P(t)L(t) - L(t)P(t)$$

This Lax pair equation implies that eigenvectors of $L$ are preserved under the flow of $P$. More precisely, this means that if we have a solution $\phi$ to the evolution equation

$$\frac{d}{dt}\phi(t) = P(t)\phi(t)$$

and at time $t = t_0$ we have

$$L(t_0)\phi(t_0) = \lambda\phi(t_0)$$

then we also have

$$L(t)\phi(t) = \lambda\phi(t)$$

for all $t$. Namely,

$$\begin{aligned}\frac{d}{dt}[L\phi - \lambda\phi] &= [P, L]\phi + LP\phi - \lambda P\phi \\ &= P[L\phi - \lambda\phi]\end{aligned}$$

Thus if $L\phi - \lambda\phi$ vanishes for some $t_0$, then it vanishes for all time.

We introduce the Lax pair which is useful for the mKdV equation. The operators $L$ and $P$ are two by two matrices of differential operators. For some real function $F(t,x)$ define the selfadjoint operator

$$L(t) = \begin{pmatrix} 0 & -i(\frac{\partial}{\partial x} + F(t,x) \\ i(-\frac{\partial}{\partial x} + F(t,x) & 0 \end{pmatrix}$$

where $F$ denotes the operator of multiplication by $F$. We remark that this operator is called a Dirac operator, since it is a square root of a Schrödinger operator:

$$L^2(t) = \begin{pmatrix} -\frac{\partial^2}{\partial x^2} + F_x + F^2 & 0 \\ 0 & -\frac{\partial^2}{\partial x^2} - F_x + F^2 \end{pmatrix}$$

Thus $L^2$ separates into two operators that are of Schrödinger type.

We consider the eigenfunction equation for $L$:

$$L\phi = k\phi$$

If we make the ansatz

$$\phi(t,k,x) = \begin{pmatrix} a(t,k,x)e^{ikx} + b(t,k,x)e^{-ikx} \\ a(t,k,x)e^{ikx} - b(t,k,x)e^{-ikx} \end{pmatrix} \tag{75}$$

then the eigenfunction equation for $L$ turns into the ordinary differential equation

$$\frac{\partial}{\partial x}(a,b) = (a,b)(0, Fe^{2ikx})$$

used to define the nonlinear Fourier transform.

Define the anti-selfadjoint operator $P = P(t)$ by

$$P = \begin{pmatrix} 4\frac{\partial^3}{\partial x^3} + 3\left\{\frac{\partial}{dx}, F_x + F^2\right\} - 4(ik)^3 & 0 \\ 0 & 4\frac{d^3}{dx^3} + 3\left\{\frac{\partial}{\partial x}, -F_x + F^2\right\} - 4(ik)^3 \end{pmatrix}$$

where $\{A, B\} = AB + BA$ denotes the anti-commutator of $A$ and $B$. The operator $P(t)$ depends on the parameter $k$, but only through an additive multiple of the identity matrix which vanishes upon taking a commutator.

Using some elementary algebraic manipulations, the Lax pair equation

$$\frac{d}{dt}L = [P, L]$$

turns into the mKdV equation (73) for $F$. Thus $L$ and $P$ as above are a Lax pair precisely when $F$ satisfies the mKdV equation.

With the ansatz (75), the evolution equation for $\phi$,

$$\frac{d}{dt}\phi = P\phi$$

becomes a partial differential equation for $a$ and $b$.

We shall now assume that $F(t,x)$ is compactly supported in $x$ and remains in a fixed compact support as $t$ evolves. This assumption is only good for the purpose of an exposition of the main ideas. In reality no solution to the mKdV equation remains supported in a compact set under the time evolution, so in a rigorous argument one needs to discuss asymptotic behaviour of the solutions for large $|x|$.

To the right and to the left of the support of $F$, the functions $a$ and $b$ are constant and we write $a(t,k,\pm\infty)$ and $b(t,k,\pm\infty)$ for the values to the right and left of the support of $F$. Outside the support of $F$, the partial differential equations for for $a(t,k,\pm\infty)$ and $b(t,k,\pm\infty)$ become

$$\frac{d}{dt}a(t,k,\pm\infty)e^{ikx} = \left[4\frac{d^3}{dx^3} - 4(ik)^3\right] a(t,k,\pm\infty)e^{ikx}$$

$$\frac{d}{dt}b(t,k,\pm\infty)e^{-ikx} = \left[4\frac{d^3}{dx^3} - 4(ik)^3\right] b(t,k,\pm\infty)e^{-ikx}$$

The right-hand side of the equation for $a$ vanishes and thus we obtain that $a$ remains constant:

$$a(t,k,\pm\infty) = a(0,k,\pm\infty)$$

The equation for $b$ reduces to

$$\frac{d}{dt}b(t,k,\pm\infty) = 8ik^3 b(t,k,\pm\infty)$$

which has the solution

$$b(t,k,\pm\infty) = b(0,k,\pm\infty)e^{8ik^3 t}$$

We now specialize to the solutions to the eigenfunction equation for $L$ such that

$$a(t,k,-\infty) = 1$$
$$b(t,k,-\infty) = 0$$

This is consistent with the evolution calculated above. Thus we have

$$\widehat{F}(t,k) = (a(t,k,\infty), b(t,k,\infty)) = (a(0,k,\infty), b(0,k,\infty)e^{8ik^3t})$$

This is the explicit form of the solution $F$ to the mKdV equation that we claimed.

## 2. Gaussian processes

We shall very briefly discuss the link between the nonlinear Fourier transform and stationary Gaussian processes. A detailed account on the subject of Gaussian processes can be found in [**8**].

Probability theory in the view of an analyst is a theory of measure and integration where the underlying measure spaces are hidden as much as possible in language and notation. For probabilists' intuition, the underlying measure spaces are uninteresting. Many statements in probability theory are fairly independent of the special structure of the underlying measure space.

A random variable $f$ is a measurable function on some measure space of total measure 1. Analysts would write $f(x)$ referring to an element $x$ of the underlying measure space. The integral of this function over the measure space is called the mean or the expectation $E(f)$ of $f$. Analysts would write $\int f(x)d\mu(x)$ referring to the measure $\mu$. A collection of random variables living on the same measure space is called a family of random variables. The measure of a set in the measure space is called the probability $P$ of the set. Inded, since the set is usually discribed by conditions on one or several random variables, the set is called the "event" that the random variables satisfy these conditions.

For example, one writes

$$P(f > \lambda)$$

for the measure $\mu(\{x : f(x) > \lambda\})$ and calls it the probability of the event $f > \lambda$. The function $P(f > \lambda)$ in $\lambda$ is called the distribution function of the random variable $f$.

A real valued random variable is called Gaussian, if

$$P(f > \lambda) = c\int_\lambda^\infty e^{-(s-s_0)^2/2\sigma}\,ds$$

where $c$ is normalized so that $P(f > -\infty) = 1$. Observe that for a Gaussian variable, the distribution function is determined by the mean value $E(f) = s_0$ and the variance

$$E((f - E(f))^2) = E(f^2) - E(f)^2 = \sigma$$

A Gaussian family indexed by $\mathbf{Z}$ consists of a family of random variables

$$f_n, \quad n \in \mathbf{Z}$$

such that each $f_n$ is Gaussian distributed with mean zero and also each finite linear combination

$$f = \sum_{n\in\mathbf{Z}} \gamma_n f_n$$

is Gaussian distributed.

In particular, $E(f_n) = 0$ for all $n$. By linearity of the expectation, $E(f) = 0$ for all finite linear combinations $f$ as above. Therefore, the distribution of each linear combination $f$ is determined by the variance $E(f^2)$. By linearity of the expectation, we have the formula

$$E(f^2) = \sum_{n,m\in\mathbf{Z}} \gamma_n\gamma_m E(f_n f_m)$$

On the space of finite sequences $(\gamma_n)$, identified as random variables $f$ as above, we can define an inner product of two elements $f$ and $f'$ by $E(ff')$. It is positive definite since $E(f^2) > 0$ for all $f$. Let $H$ be the Hilbert space closure of this inner product space. The elements of this Hilbert space are again random variables.

Orthogonality in this space translates to the probabilistic notion of independence. Independence means a factorization of the distribution functions:

$$P(f > \lambda, f' > \lambda') = P(f > \lambda)P(f' > \lambda')$$

If two random variables are orthogonal in $H$, then

$$E(\eta_1\eta_2) = 0 = E(\eta_1)E(\eta_2)$$

and $E(\eta_1\eta_2) = E(\eta_1)E(\eta_2)$ is a necessary condition for independence. In the setting of Gaussian variables, it is also sufficient for independence.

If $H'$ is a subspace of $H$, then every $f \in H$ can be split as

$$f = f' + f''$$

where $f'$ is in $H$ and $f''$ is in the orthogonal complement of $H'$. The probabilistic interpretation is that $f'$ is known if all elements in $H$ are known, while $f''$ is independent of any knowledge about $H'$.

A stationary Gaussian process is one for which

$$E(f_n f_m) = Q(n-m)$$

for some sequence of numbers $Q \in l^\infty(\mathbf{Z})$. Equivalently, a Gaussian process is stationary if it is equal to the shifted proces $\tilde{f}_n = f_{n-1}$. In particular, the length of the vectors $f_n$ in a stationary Gaussian process is independent of $n$.

We observe that $Q(n) = E(f_n, f_0) \leq E(f_0, f_0) = Q(0)$ for all $n$. This is a special case of the property

$$\sum_{n,m} Q(n-m)\gamma_n\gamma_m \geq 0$$

for all finite sequences $\gamma_n$. Indeed, the left-hand side has the meaning of $E(f^2)$ for some $f$ and is therefore non-negative. This property is called positive definiteness of the sequence $Q$. The following theorem by Bochner characterizes all positive definite sequences.

**Theorem 15.** *A nonzero sequence $Q_n$ satisfies*

$$\sum_{n,m} Q(n-m)\gamma_n\gamma_m \geq 0$$

*for all real valued finite sequences $\gamma_n$ if and only if it is the Fourier series of a positive measure on $\mathbf{T}$:*

$$Q_n = \int_{\mathbf{T}} z^n \, d\mu$$

We only prove one direction of the theorem. Assume $Q_n$ is the Fourier series of a positive measure. Then

$$\sum_{n,m} Q(n-m)\gamma_n\gamma_m = \sum_{m,n} \int z^{n-m}\gamma_n\gamma_m \, d\mu$$
$$= \int |\sum_n \gamma_n z^n|^2 \, d\mu \geq 0$$

This proves one direction of Bochner's theorem.

Now there is an evident isometric isomorphism from $H$ to $L^2(\mu)$. It maps the elements $f_n$ to the function $z^n$. Isometry is seen as follows:

$$E(f_n, f_m) = Q(n-m) = \int z^{n-m} \, d\mu = \langle z^n, z^m \rangle_{L^2(\mu)}$$

Surjectivity follows from the definition of $L^2(\mu)$ as the closure of the linear span of the monomials $z^n$.

The space $L^2(\mu)$ takes us into the setting of orthogonal polynomials on the circle $\mathbf{T}$. Indeed, by reflection we consider polynomials in $z^{-1}$ and have the following immediate consequence of Szegö's theorem:

**Corollary 2.** *The past, namely the span of $f_{-1}, f_{-2}, \ldots$ determines the present $f_0$, i.e., $f_0$ is in the closed span of $f_{-1}, f_{-2}, \ldots$ if and only if*

$$\int_{\mathbf{T}} \log|w| = -\infty$$

*where $w$ is the absolutely continuous part of $\mu$.*

# LECTURE 7
## Appendix: Some Background material

This lecture has two sections. The first section gives some background material on boundary regularity of harmonic and holomorphic functions on the unit disc. This section contains several theorems that are important for a rigorous development of the nonlinear Fourier transform.

The second section recalls some facts about the group $Sl_2(\mathbf{R})$ and the isomorphic group $SU(1,1)$. This section is meant to help understand some of the algebraic manipulations done in this lecture series on a group theoretical level, but otherwise is somewhat irrelevant for the overall understanding of the lecture series.

## 1. The boundary behaviour of holomorphic functions

We shall study classes of holomorphic functions on the unit disc defined by some size control.

For example, for any monotone increasing function $\phi : \mathbf{R}_{\geq 0} \to \mathbf{R}_{\geq 0}$ we can consider the space of holomorphic functions

$$\{f : \sup_{r<1} \int_{\mathbf{T}} \phi(|f(r\cdot)|) < \infty\} \tag{76}$$

Here $\int_{\mathbf{T}}$ denotes the integral over the circle $\mathbf{T}$, i.e., the set of all $z \in \mathbf{C}$ with $|z| = 1$, with the usual Lebesgue measure on $\mathbf{T}$ normalized to have total mass 1.

$$\int_{\mathbf{T}} f := \frac{1}{2\pi} \int_0^{2\pi} f(e^{i\theta})\, d\theta$$

Thus for fixed $r$ the integral in (76) has the meaning of an average over the circle of radius $r$ about the origin.

We shall mainly be interested in $\phi(x) = x^p$ (producing Hardy spaces) and $\phi(x) = \log_+(x) = \max(0, \log|x|)$ (producing Nevanlinna class).

The largest space we consider and the space containing all functions we shall be concerned with is the Nevanlinna class:

$$N = \{f : \sup_{r<1} \int_{\mathbf{T}} \log_+ |f(r\cdot)| < \infty\}$$

This space can be identified with a space of almost everywhere defined functions on the circle, as the following theorem describes:

**Theorem 16.** *If $f \in N$, then $f$ has radial limits*

$$\lim_{r \to 1} f(rz)$$

*for almost every $z$ in $\mathbf{T}$. (Here $r$ is real and less than 1.) If the radial limits of two functions $f_1, f_2 \in N$ coincide on a subset of $\mathbf{T}$ of positive measure, then the functions are equal.*

The proof of this theorem will be discussed below.

The uniqueness result is not special to Nevanlinna class, provided one passes to the notion of nontangential limits. A precise version of this statement is given in the theorem below. The previous theorem also holds with "nontangential" in place of "radial", but we shall not need this.

**Theorem 17.** *Assume that two holomorphic functions on the disc each have nontangential limits on sets of positive measure, and the limits coincide on a set of positive measure, then the two functions are equal.*

*This statement is false if "nontangential" is replaced by "radial".*

For a proof of this theorem and a discussion of nontangential limits see Garnett's book [**9**].

Thus we can talk about "the holomorphic extension" of a function defined on a positive set on $\mathbf{T}$, provided such an extension exists.

One of the reasons for us to try to identify holomorphic functions with their limits on the boundary is that on the Riemann sphere, the circle $\mathbf{T}$ is the boundary of the unit disc about 0 and the boundary of the unit disc about $\infty$. We would like to study simultaneously and compare the spaces of holomorphic functions on both discs. The boundary values of the functions on $\mathbf{T}$ are the only link between the two spaces.

While knowledge of the real part of a holomorphic function on the disc is sufficient to determine the imaginary part up to an additive constant, the above uniqueness result heavily relies on the fact that we know the limits of both the real and imaginary part. In particular, the analogue statement fails for harmonic functions, as the example $\operatorname{Re}\left(\frac{z+1}{z-1}\right)$ shows, which has limit 0 almost everywhere on the circle $\mathbf{T}$ but clearly is not zero.

Thus the notion of harmonic extension cannot be used as freely as that of holomorphic extension, indeed it cannot be used easily in the context of functions defined almost everywhere. However, harmonic extensions of measures are well defined. This is described in the following theorem:

**Theorem 18.** *Given a complex Borel measure $\mu$ on $\mathbf{T}$ (an element in the dual space of the space of continuous functions with the supremum norm), the function*

$$f(z) = \int P_z \, d\mu$$

*where $P_z$ is the Poisson kernel,*

$$P_z(\zeta) = \operatorname{Re}\left(\frac{\zeta + z}{\zeta - z}\right)$$

*is a harmonic function on the disc. Radial limits of this function exist almost everywhere and coincide almost everywhere with the density $f \in L^1$ of the absolutely continuous part of $\mu$, i.e.*

$$\mu = \frac{1}{2\pi} f d\theta + \mu_s$$

*where $\mu_s$ is singular with respect to Lebesgue measure $d\theta$*

**Proof.** The kernel $P_z$ is harmonic in $z \in D$, and thus the superposition $\int P_z \, d\mu$ of harmonic functions is again harmonic (use the mean value characterization of harmonicity and Fubini's theorem). To study the existence and behaviour of radial limits, it suffices to consider separately the cases of $\mu$ absolutely continuous and $\mu$ singular with respect to Lebesgue measure. If $\mu$ is absolutely continuous, $\mu = \frac{1}{2\pi} f d\theta$, we estimate

$$|f(z) - \int_{\mathbf{T}} P_{rz} f| \leq |f(z) - g(z)| + |g(z) - \int_{\mathbf{T}} P_{rz} g| + |\int_{\mathbf{T}} P_{rz}(f-g)|$$

$$\leq |f(z) - g(z)| + |g(z) - \int_{\mathbf{T}} P_{rz} g| + C|M(f-g)(z)|$$

where $M$ is the Hardy - Littlewood maximal function, which up to a constant factor $C$ dominates the integration against the Poisson kernel independently of $r$, and $g$ is some appropriate smooth function. If $g$ is sufficiently close to $f$ in $L^1$ norm, then outside a small set the difference $|f-g|$ is small. Outside a possibly different small set, $M(f-g)$ is small by the Hardy - Littlewood maximal theorem. The term $g(z) - \int_{\mathbf{T}} P_{rz} g$ can be made small by choosing $r$ close to 1 depending on the choice of the smooth function $g$. Making this argument rigorous using the correct quantifyers one proves convergence of $\int_{\mathbf{T}} P_{rz} f$ to $f(z)$ outside a set of arbitrarily small measure, i.e., almost everywhere.

If $\mu$ is singular, one proves using a Vitali - type covering lemma that $\mu$ has vanishing density almost everywhere,

$$\lim_{h \to 0} \mu([t-h, t+h])/2h = 0$$

for almost every $t$. For points $t$ of zero density one then observes that $P_r * f(t)$ tends to 0. □

Observe that if $\mu$ is absolutely continuous with a continuous density function, then one can prove that the harmonic function $f$ defined in the above theorem has continuous extension to $D \cup \mathbf{T}$. On $\mathbf{T}$ the function $f$ coincides with the density of $\mu$. Moreover, by the maximum principle, this extension is the unique harmonic function which has continuous extension to $\mathbf{T}$ coinciding with the density function of $\mu$.

The following is a variant of the above theorem:

**Theorem 19.** *Given a real measure $\mu$ on $\mathbf{T}$, the function*

$$f(z) = \int \operatorname{Im}\left(\frac{\zeta + z}{\zeta - z}\right) d\mu(\zeta)$$

*is a harmonic conjugate to the function defined in the previous theorem in the unit disc. Its radial limits exist almost everywhere and are equal almost everywhere to the Hilbert transform of $\mu$.*

**Proof.** Harmonicity (holomorphicity) follows again by characterizing harmonic (holomorphic) functions by mean value (Cauchy) integrals and then using Fubini's theorem and the fact that the kernel $\frac{\zeta+z}{\zeta-z}$ is holomorphic. To see that radial limits exist, one has to study convergence of the conjugate Poisson kernels, which amounts to estimating maximal truncated singular integrals. We leave details as an exercise. □

There is a nice intrinsic characterization of those harmonic functions which are extensions of positive measures, due to Herglotz [**10**].

**Theorem 20.** *Any positive harmonic function is the harmonic extension of a unique positive measure.*

We shall sometimes call a holomorphic function whose real or imaginary part is positive a Herglotz function.

A real harmonic function is the extension of a measure if it can be written as $f_1 - f_2$ with two positive harmonic functions.

**Proof.** If $f$ is a positive harmonic function, then for each radius $r$ there is a positive measure $\mu_r$ on $T$ with density $f(rz)$. Let $\mu$ be a weak-$*$ limit as $r \to 1$ of an appropriate subsequence of this collection of measures (each of them has total mass equal to $f(0)$ by the mean value property). By weak convergence, the Poisson extension of $\mu$ is equal to the pointwise limit of the Poisson extensions of $\mu_r$, and thus equal to $f$. This proves existence of a measure whose Poisson extension is $f$. If $\mu_1$ and $\mu_2$ are two measures with the same harmonic extension, then the difference measure has vanishing Fourier coefficients (this follows from calculating the Taylor coefficients of the harmonic extension at 0), and thus is zero. □

Just as the real and imaginary parts of boundary values of a holomorphic function do not individually determine the function, neither does the absolute value of the boundary values. Indeed, the absolute value of the boundary value function a.e. does not even in general determine membership in a space. For Nevanlinna class, observe that Fatou's theorem implies (let $f$ also denote the boundary value function on the circle)

$$\int_{\mathbf{T}} \log_+ |f| \leq \sup_{r<1} \int_{\mathbf{T}} \log_+ |f(r\cdot)|$$

Thus

$$\int_{\mathbf{T}} \log_+ |f| < \infty$$

is a necessary condition for $f$ to belong to $N$. But it is not a sufficient condition since Fatou's inequality can be strict and finiteness of

$$\int \log_+ |f|$$

does not imply membership in Nevanlinna class. An example is the function $\exp(\frac{1+z}{1-z})$, which has absolute value 1 almost everywhere on the circle $\mathbf{T}$ but is not Nevanlinna class.

Similar statements hold for most function spaces, in particular the Hardy spaces.

Let us note another application of Fatou's theorem to Nevanlinna class functions.

**Lemma 28.** *If $f \in N$ and $f \neq 0$, then $\log |f|$ is absolutely integrable on $\mathbf{T}$.*

The point of the lemma is that originally we control only the positive part of $\log |f|$, but via the lemma we also control the negative part.

**Proof.** By dividing by a power of $z$ if necessary, we may assume $f(0) \neq 0$ (here we use $f \neq 0$). As $\log|f(z)|$ is subharmonic in the disc, we have

$$\log|f(0)| \leq \int \log|f(r\cdot)| = \int \log_+ |f(r\cdot)| + \int \log_- |f(r\cdot)|$$

The first inequality is a consequence of subharmonicity and can be shown by Green's theorem (analoguously to the proof of the mean value property of harmonic functions), using the fact that the distributional Laplacian of $\log|f|$ is a positive measure (instead of zero as for harmonic functions).

Using the Nevanlinna class assumption we observe

$$\log|f(0)| - C \leq \int \log_- |f(r\cdot)|$$

and so Fatou's lemma implies

$$\log|f(0)| - C \leq \int \log_- |f|$$

□

Observe that this lemma also proves the uniqueness part of Theorem 16. Since the difference of two Nevanlinna functions is Nevanlinna again, it suffices to prove that if the limit function vanishes on a set of positive measure, then it is constant 0. However, if the limit function vanishes on a set of positive measure, then $\log|f|$ is not integrable on the circle, so by the lemma $f$ is identically equal to 0.

When estimating the size of holomorphic functions, it is natural to consider the logarithm of $|f|$, which produces a harmonic function if $f$ is zero free and a subharmonic function if $f$ has zeros. Thus Nevanlinna class functions are related to harmonic extensions because for zero free Nevanlinna functions $\log|f|$ is the harmonic extension of a measure.

**Lemma 29.** *If $f \in N$ then $f$ can be factored as $f_1 f_2$ where $f_1$ is an analytic function in $D$ bounded by* 1 *and* $\log|f_2|$ *is the harmonic harmonic extension of a measure. If $f$ has no zeros, we may pick $f_1 = 1$.*

**Proof.** As before, we consider the measures $\mu_r$ defined by

$$\int_{\mathbf{T}} g d\mu_r = \int_{\mathbf{T}} g(\cdot) \log_+ |f(r\cdot)|$$

By the Nevanlinna property of $f$, these measures have bounded total mass (independent of $r$), and thus there is a weak-$*$ limit, $\mu$, of a subsequence of these measures. Let $f_2$ be a function in $N$ such that $\log|f_2|$ is the harmonic extension of $\mu$. Using subharmonicity of $f$ and a limiting process, one can show that $f_2$ dominates $f$, thereby proving the above theorem.

If $f$ has no zeros, then $f/f_2$ has no zeros, and $\log(f/f_2)$ is a holomoprhic function with negative real part. Thus its real part is the harmonic extension of a negative measure. Thus $\log(f)$ is the harmonic extension of a measure. □

We can now prove existence of radial limits almost everywhere for any Nevanlinna function $f$. It suffices to prove existence for $f_1$ and $f_2$ where $f = f_1 f_2$ is a splitting as in the last lemma. By adding a constant to $f_1$, we may assume that $f_1$ has positive real part. For such functions we proved existence of the radial limits almost everywhere before. The function $f_2$ has no zeros, and thus $\log(f_2)$ is holomorphic and its real part is positive. Thus radial limits of $f_2$ exist.

The bounded function $f_1$ cannot be omitted from the last theorem, because $f$ may have zeros. Using Blaschke products, one may choose $f_1$ to be a possibly infinite Blaschke product. This can be deduced from the following lemma:

**Lemma 30.** *Let $f \in N$ and let $z_n$ be the zeros of $f$ (multiple zeros appear in the sequence according to multiplicity). Then*

$$\sum_n (1 - |z_n|) < \infty$$

*Conversely, if $z_n$ is any such sequence with $0$ appearing $m$ times, then the (Blaschke) product*

$$B(z) = z^m \prod_{z_n \neq 0} \frac{-\overline{z_n}}{|z_n|} \frac{z - z_n}{1 - z\overline{z_n}}$$

*converges uniformly on compact subsets of $D$ to a bounded analytic function with exactly the zeros $z_n$. The radial limits of $B$ on $\mathbf{T}$ have modulus $1$ almost everywhere on $\mathbf{T}$.*

**Proof.** See Garnett's book [**9**]. □

We define an outer function to be a Nevanlinna function without zeros such that $\log|f|$ is the harmonic extension of an absolutely continuous measure. Thus $\log|f|$ is the Poisson integral of its a.e. radial limits.

Outer functions are very special functions: they are determined (up to a constant phase factor) by the modulus of their limits almost everywhere.

This relates to the following lemma, which is a form of an "inverse Fatou" for outer functions.

**Lemma 31.** *If the boundary value functions of an outer function is in $L^p$, then the outer function is in $H^p$.*

**Proof.** By Poisson extension we have

$$\log|f(z)| = \int P_z(.) \log|f(.)|$$

By convexity of the function $e^{rx}$ and Jensen's inequality we have

$$|f(z)|^p \leq \int P_z(.)|f(.)|^p$$

This implies that the restrictions of $f$ to smaller circles are uniformly in $L^p$. □

We will use the following criterion for outerness.

**Lemma 32.** *Let $f$ and $1/f$ be in $H^p(D)$ for some $p > 0$. Then $f$ is outer.*

**Proof.** Since $f$ and $1/f$ are hoplomoprhic in $D$, we can choose a branch of $\log(f)$ which is holomorphic in $D$. The $L^p$ estimates for $f$ and $f^{-1}$ can be used to obtain $L^2$ estimates of $\log(f)$ on circles of radius $r$ about 0, uniformly in $r$. Thus $\log(f)$ in $H^2(D)$, which implies it is the Poisson extension of its boundary values. This implies that $f$ is outer. □

## 2. The group $Sl_2(\mathbf{R})$ and friends

The general linear group $Gl_2(\mathbf{C})$ consists of all $2 \times 2$ invertible complex matrices with the usual matrix product as group multiplication.

This group acts on the complex vector space $\mathbf{C}^2$ by linear transformations

$$\begin{pmatrix} u \\ v \end{pmatrix} \to \begin{pmatrix} a & b \\ c & d \end{pmatrix} \begin{pmatrix} u \\ v \end{pmatrix}$$

Indeed, $Gl_2(\mathbf{C})$ can be identified as the group of linear automorphisms of $\mathbf{C}^2$. The projective (complex) line $P^1$ is the set of all complex lines in $\mathbf{C}^2$ of the form

$$L_{\alpha,\beta} = \{(\alpha z, \beta z), z \in \mathbf{C}\}$$

with parameters $(\alpha, \beta)$ and $\alpha\beta \neq 0$.

The projective line $P_1$ is a complex manifold with two charts $\mathbf{C} \to P^1$ given by $z \to L_{z,1}$ and $z' \to L_{1,z'}$. The first chart misses only the line $L_{1,0}$ and the second chart misses only the line $L_{0,1}$. The transition between the two charts is given by $z = 1/z'$.

Thus $P^1$ is isomorphic to the Riemann sphere. We shall mainly use the first chart described above and write $z = \infty$ for the line $L_{1,0}$.

The action of $Gl_2(\mathbf{C})$ on $P_1$ is then given by the linear fractional transformation

$$z \to \frac{az+b}{cz+d}$$

Thus the action is by biholomorphic maps (Möbius transforms) of the Riemann sphere. Indeed, every biholomorphic self map of the Riemann sphere has to be a fractional linear transformation (exercise), and thus be given by action of an element in $Gl_2(\mathbf{C})$. Thus we have a surjection of $Gl_2(\mathbf{C})$ onto the set of Möbius transforms of the Riemann sphere.

Two elements in $Gl_2(\mathbf{C})$ give the same Möbius transform, if and only if they are scalar multiples of each other. (The trivial action only comes from the scalar matrices). The group of matrices in $Gl_2(\mathbf{C})$ with determinant one is called $Sl_2(\mathbf{C})$. Since every matrix can be normalized to have determinant one (if $\det(g) = \lambda$, then $\det(\nu g) = 1$ if $\nu^2 = \lambda^{-1}$, an equation that can always be solved for $\nu$ in the complex numbers), $Sl_2(\mathbf{C})$ still covers all Möbius transforms. However, the two elements id and $-$id of $Sl_2(\mathbf{C})$ both map onto the identity Möbius transform. Thus $Sl_2(\mathbf{C})$ is a double cover of the group of Möbius transforms. The quotient of $Sl_2(\mathbf{C})$ by the central subgroup with two elements id and $-$id is called $PSL_2(\mathbf{C})$. The group of Möbius transforms of the sphere is isomorphic to $PSL_2(\mathbf{C})$.

The group of Möbius transforms which leave the real line (a great circle) on the Riemann sphere invariant, has to come from an automorphism of $\mathbf{C}^2$ which maps real vectors in $\mathbf{C}^2$ to real vectors, and thus has to come from a real linear automorphism of $\mathbf{R}^2$.

These maps precisely give the matrices in $Gl_2(\mathbf{C})$ with real entries. The group of these matrices is $Gl_2(\mathbf{R})$. This group has two connected components, the component of elements with positive determinant and the component of elements with negative determinant. The elements of the first component map the upper half plane (positive imaginary part) to the upper half plane, the elements of the other component map the upper half plane to the lower half plane. The group of matrices in $Gl_2(\mathbf{R})$ with determinant 1 is called $Sl_2(\mathbf{R})$. The group can be identified as a double cover of all biholomorphic self maps of the upper half plane. The quotient by the central subgroup of two elements is called $PSl_2(\mathbf{R})$.

By conformal equivalence, more precisely by a rotation of the Riemann sphere, the Möbius transforms of the upper half plane correspond to the Möbius transforms on the unit disc $D = \{z : |z| < 1\}$. The latter are matrices in $Sl_2(\mathbf{C})$ which leave the set of vectors of the form $(e^{i\phi}z, z)$ with $\phi \in \mathbf{R}$ invariant. These vectors are exactly the null vectors of the quadratic form

$$B(u,v) = |u|^2 - |v|^2$$

A linear map preserving the null vectors of this form has to leave the whole quadratic form invariant up to a scalar multiple. By the determinant constraint, this multiple has to be 1 or $-1$, again corresponding to the maps which map inside of the unit circle to inside, or inside to outside respectively. The matrices in $Sl_2(\mathbf{C})$ which leave $B$ invariant are precisely those that can be written in the form

$$\begin{pmatrix} a & b \\ \overline{b} & \overline{a} \end{pmatrix}$$

with $|a|^2 - |b|^2 = 1$. These matrices form the group $SU(1,1)$. It is a double cover of $PSU(1,1)$, which is $SU(1,1)/\{\mathrm{id}, -\mathrm{id}\}$. By the above discussion this group is isomorphic to $Sl_2(\mathbf{R})$, and an explicit isomorphism can be given by conjugating with a Möbius transform mapping the upper half plane to the disc. An explicit isomorphism is given by

$$Sl_2(\mathbf{R}) \to SU(1,1)$$

$$\begin{pmatrix} a & b \\ c & d \end{pmatrix} \to \begin{pmatrix} (a+d)/2 + i(b-c)/2 & (a-d)/2 + i(b+c)/2 \\ (a-d)/2 - i(b+c)/2 & (a+d)/2 - i(b-c)/2 \end{pmatrix}$$

Let us discuss eigenvalues of matrices in $Sl_2(\mathbf{R})$. The eigenvalues of a matrix $G$ in $Sl_2(\mathbf{R})$ satisfy the equation

$$\lambda^2 - \mathrm{Tr}(G)\lambda + 1 = 0$$

and the trace, $\mathrm{Tr}(G)$, is a real number. For $|\mathrm{Tr}(G)| < 2$, the two solutions are conjugates; they are distinct and have modulus one. In this case we call the group element elliptic. For $|\mathrm{Tr}(G)| = 2$, there is a double root 1 or $-1$. Such group elements (in general Jordan blocks) are called parabolic. For $|\mathrm{Tr}(G)| > 2$, we have two distinct real roots. Such group elements are called hyperbolic.

The Möbius transformation associated to a matrix in $Sl_2(\mathbf{C})$ will have two fixed points on the Riemann sphere unless it is the identity transformation (which fixes all points), or is a Jordan block and so has only one fixed point. For matrices in $SU(1,1)$ or $Sl_2(\mathbf{R})$ the classification into elliptic, parabolic, or hyperbolic points can be understood by the location of these fixed points relative to the domain (disc, half plane).

The elliptic elements have a fixed point in the interior of the domain (disc, half plane). Rotations of the unit disc are easy examples, but any fixed point is possible. Parabolic Möbius transforms have a single fixed point on the boundary of the domain. Horizontal translation of the upper half plane is an example. Hyperbolic elements have two fixed points on the boundary. Multiplication of the upper half plane by a positive scalar is such an example.

Intuitive geometric coordinates on $SU(1,1)$ are

$$(b/|a|, \arg(a)) \in D \times \mathbf{R}/2\pi\mathbf{Z}$$

Observe that the modulus of $b/|a|$ determines the modulus of $a$ and $b$ via the constraint $|a|^2 - |b|^2 = 1$. The argument of $b$ is equal to the argument of $b/|a|$, and the argument of $a$ is given. Thus the above coordinates indeed determinant $a$ and $b$. Since $b/|a|$ lies in the open unit disc, the group $SU(1,1)$ can be visualized as a solid open torus or equivalently an infinite cylinder in $\mathbf{C} \times \mathbf{R}$ where the last coordinate is taken modulo $2\pi$. The main axis $b/|a| = 0$ consists of the elements in the compact subgroup of elements

$$\begin{pmatrix} e^{i\phi} & 0 \\ 0 & e^{-i\phi} \end{pmatrix}$$

In $Sl_2(\mathbf{R})$ these elements correspond to the rotations

$$\begin{pmatrix} \cos(\phi) & \sin(\phi) \\ -\sin(\phi) & \cos(\phi) \end{pmatrix}$$

Rotating about the main axis by an angle $2\phi$ corresponds to multiplying $b$ by a phase factor $e^{2i\phi}$. This can be achieved by conjugating with the previously displayed diagonal element of $SU(1,1)$. Thus rotating the torus about the main axis is an inner automorphism.

Reflecting across the plane determined by requiring $b$ to be real corresponds to replacing $b$ by its complex conjugate. This corresponds to transposing the matrix, which is an anti-automorphism (it changes the order of multiplication). Using the previous rotation symmetry, reflecting across any plane through the main axis is an anti-automorphism.

General Lie theory puts the tangent vectors at the identity element in one-to-one correspondence with one parameter subgroups. We claim that any such subgroup lies in the plane spanned by the tangent vector and the main axis. For the compact subgroup along the main axis this is trivially evident. For any other subgroup, observe that reflecting across this plane gives another one parameter subgroup (these groups are commutative) with the same tangent vector, and thus the reflected subgroup has to be the same as the original subgroup. Thus the subgroup lies inside the plane.

Moreover, we claim that the subgroups are contained in traces of the form

$$\sin(\arg(a))|a|/|b| = const \tag{77}$$

Which is to be read projectively.

Indeed it suffices to prove this for those subgroups with real $b$. Consider two matrices in $SU(1,1)$ with parameters $a, b$ and $a', b'$ and assume $b, b'$ real. If the two elements are in the same subgroup, than the off diagonal element of the product is again real,

$$0 = \mathrm{Im}\,(ab' + b\overline{a'}) = |a|\sin(\arg(a)b' - |a|'\sin(\arg(a'))b$$

This proves the claim.

If the modulus of the constant in (77) is less than one, the entire trace determined by the above equation is contained in the solid torus. The entire trace is a subgroup consisting of elliptic elements. The subgroup is compact. We remark that the special subgroup consisting of the main axis of the infinite cylinder is only special in the chosen coordinates. It is by inner automorphisms equivalent to any of the other elliptic subgroups.

If the constant in (77) has modulus equal to one, the solution set of the equation (77) meets the boundary of the solid torus. Intersecting with the open torus, we obtain two connected components. The component containing the identity element is a non-compact subgroup, the other component is $-1$ times this subgroup. All elements in these groups and remainder class are parabolic.

If the constant has modulus larger than one, then the trace of the equation intersected with the torus has again two components. On component is a subgroup consisting of hyperbolic elements, the other component is $-1$ times the subgroup and also consists of hyperbolic elements. (The two components are mapped onto the same group in $PSU(1,1)$.)

If the constant is infinite, which means $\sin(\arg(a)) = 0$, then the trace consists of two lines, one through the origin is a hyperbolic subgroup, the other one through $-1$ is $-1$ times this subgroup.

Note that all automorphisms of the group leave the infinitesimal cone of parabolic elements near the origin invariant. The group acts naturally on its Lie algebra, which is $\mathbf{R}^3$. Since the group leaves a cone invariant, it is easy to see that it acts as $SO(2,1)$. Thus there is a map from $Sl_2(\mathbf{R})$ to $SO(2,1)$. The kernel of this map consists of the identity matrix and minus the identity matrix, therefore there is an embedding of $PSL_2(\mathbf{R})$ into $SO(2,1)$.

*Exercise:* If $K$ is the compact subgroup of $SU(1,1)$ of diagonal elements, calculate the residue classes $SU(1,1)/K$ and $K\backslash SU(1,1)$. They are spirals in the solid torus representing $SU(1,1)$.

**Lemma 33.** *If $G \in SU(1,1)$, then*

$$||G||_{op} = |a| + |b|$$
$$\log ||G||_{op} = \operatorname{arccosh}(|a|) = \operatorname{arcsinh}(|b|)$$

**Proof.** The operator norm of $G$ does not change if we multiply from the left or from the right by an element in the subgroup of diagonal elements. Thus we may assume that $a$ and $b$ are real and positive. Now the matrix is real and symmetric and thus its operator norm is equal to the maximal eigenvalue. However, a basis of eigenvectors is $(1,1)$ and $(1,-1)$ so the eigenvalues are $|a| + |b|$ and $|a| - |b|$. This proves the lemma. $\square$